%% file: ConvexCurvature_140207.tex
% With \loadreferencesfalse, the file makes its own counters and saves them to "references.tex".
% With \loadreferencestrue, the file loads the counters saved in "references.tex".
%
\newif\ifloadreferences\loadreferencestrue %
%
\input preamble %
%
% In order to use "dvips", enter:
% dvips -N0 -Z0 -K0 [whatever.dvi] -o [whatever.ps]
%
\newref{Caffarelli}{Caffrelli L., Monge-Amp\`ere equation, div-curl theorems in Lagrangian coordinates, compression and rotation, Lecture Notes, 1997}
\newref{CaffNirSprI}{Caffarelli L., Nirenberg L., Spruck J., The Dirichlet problem for nonlinear second-order elliptic equations. I. Monge Amp\`ere equation, {\sl Comm. Pure Appl. Math.} {\bf 37} (1984), no. 3, 369--402}
\newref{CaffNirSprII}{Caffarelli L., Kohn J. J., Nirenberg L., Spruck J., The Dirichlet problem for nonlinear second-order elliptic equations. II. Complex Monge Amp\`ere, and uniformly elliptic, equations, {\sl Comm. Pure Appl. Math.} {\bf 38} (1985), no. 2, 209--252}
\newref{CaffNirSprIII}{Caffarelli L., Nirenberg L., Spruck J., The Dirichlet problem for nonlinear second-order elliptic equations. III. Functions of the eigenvalues of the Hessian, {\sl Acta Math.} {\bf 155} (1985), no. 3-4, 261--301}
\newref{CaffNirSprV}{Caffarelli L., Nirenberg L., Spruck J., Nonlinear second-order elliptic equations. V. The Dirichlet problem for Weingarten hypersurfaces, {\sl Comm. Pure Appl. Math.} {\bf 41} (1988), no. 1, 47--70}
\newref{CalI}{Calabi E., Improper affine hyperspheres of convex type and a generalization of a theorem by K. J\"orgens, {\sl Michigan Math. J.} {\bf 5} (1958), 105--126}
\newref{CrandallIshiiLions}{Crandall M. G., Ishii H., Lions P. L., User's guide to viscosity solutions of second order partial differential equations, {\sl Bull. Amer. Math. Soc.} {\bf 27} (1992), no. 1, 1--67}
\newref{GilbTrud}{Gilbard D., Trudinger N. S., {\sl Elliptic Partial Differential Equations of Second Order}, }
\newref{GuanSpruckA}{Guan B., Spruck J., Boundary value problems on $S^n$ for surfaces of constant Gauss curvature, {\sl Ann. of Math.} {\bf 138} (1993), 601--624}
\newref{Guan}{Guan B., The Dirichlet problem for Monge-Amp\`ere equations in non-convex domains and spacelike hypersurfaces of constant Gauss curvature, {\sl Trans. Amer. Math. Soc.} {\bf 350} (1998), 4955--4971}
\newref{GuanSpruck}{Guan B., Spruck J., The existence of hypersurfaces of constant Gauss curvature with prescribed boundary, {\sl J. Differential Geom.} {\bf 62} (2002), no. 2, 259--287}%
\newref{GuanSpruckI}{Guan B., Spruck J., Locally convex hypersurfaces of constant curvature with boundary, {\sl Comm. Pure Appl. Math.} {\bf 57} (2004), no. 10, 1311--1331}%
\newref{GuanSpruckII}{Guan B., Spruck J., Szapiel M., Hypersurfaces of constant curvature in hyperbolic space. I., {\sl J. Geom. Anal.} {\bf 19} (2009), no. 4, 772--795}%
\newref{GuanSpruckIII}{Guan B., Spruck J., Hypersurfaces of constant curvature in hyperbolic space. II., {\sl J. Eur. Math. Soc. (JEMS)} {\bf 12} (2010), no. 3, 797--817}
\newref{GuanSpruckIV}{Guan B., Spruck J., Convex hypersurfaces of constant curvature in hyperbolic space, {\sl in preparation}}
\newref{IvochTomi}{Ivochkina N. M., Tomi F., Locally convex hypersurfaces of prescribed curvature and boundary, {\sl Calc. Var. Partial Differential Equations} {\bf 7} (1998), no. 4, 293--314}
\newref{Jorgens}{J\"orgens K., \"Uber die L\"osungen der Differentialgleichung $rt-s^2=1$  (German), {\sl Math. Ann.} {\bf 127} (1954), 130--134}%
\newref{Lab}{Labourie F., Un lemme de Morse pour les surfaces convexes (French), {\sl Invent. Math.} {\bf 141} (2000), no. 2, 239–297}
\newref{Pog}{Pogorelov A. V., On the improper convex affine hyperspheres, {\sl Geometriae Dedicata} {\bf 1} (1972), no. 1, 33--46}
\newref{SmiDTI}{Rosenberg H., Smith G., Degree Theory of Immersed Hypersurfaces, {\sl in Preparation}}
\newref{RosSpruck}{Rosenberg H., Spruck J., On the existence of convex hypersurfaces of constant Gauss curvature in hyperbolic space, {\sl J. Differential Geom.} {\bf 40} (1994), no. 2, 379–409}
\newref{ShengUrbasWang}{Sheng W., Urbas J., Wang X., Interior curvature bounds for a class of curvature equations. (English summary), {\sl Duke Math. J.} {\bf 123} (2004), no. 2, 235--264}
\newref{Smale}{}
\newref{SmiAAT}{}
\newref{SmiCGC}{Smith G., Constant Gaussian Curvature Hypersurfaces in Hadamard Manifolds, arXiv:0912.0248}
\newref{SmiSLC}{Smith G., Special Lagrangian curvature, arXiv:math/0506230}
\newref{SmiNLD}{Smith G., The non-linear Dirichlet problem in Hadamard manifolds, arXiv:0908.3590}
\newref{SmiFCS}{Smith G., Moduli of flat conformal structures of hyperbolic type, arXiv:0804.0744}
\newref{SmiPPH}{Smith G., The Plateau problem in Hadamard manifolds, arXiv:1002.2982}
\newref{SmiNLP}{Smith G., The non-linear Plateau problem in non-positively curved manifolds, arXiv:1004.0374}
\newref{SmiPPG}{Smith G., {\sl The Plateau problem for Gaussian curvature}, arXiv:1206.5544}
\newref{Spruck}{Spruck J., }
\newref{Trud}{Trudinger N. S., On the Dirichlet problem for Hessian equations, {\sl Acta Math.}, {\bf 175}, (1995), 151--164}
\newref{TrudWang}{Trudinger N. S., Wang X., On locally locally convex hypersurfaces with boundary, {\sl J. Reine Angew. Math.} {\bf 551} (2002), 11--32}
\newref{WhiteI}{White B., The space of $m$-dimensional hypersurfaces that are stationary for a parametric elliptic functional, {\sl Indiana Univ. Math. J.}, {\bf 36}, (1987), no. 3, 567--602}
\newref{WhiteII}{White B., The space of minimal submanifolds for varying Riemannian metrics, {\sl Indiana Univ. Math. J.}, {\bf 40}, (1991), no. 1, 161--200}
%
%%%%%%%%%%%%%%%%%%%%%%%%%%%%%%%%%%%%%%%%%%%%%%%%%%%%%%%%%%%%%%%%%%%%%%%%%%%%%%%%%%%%%%%%%%%%%%%%%%%%%%%%%%%%%%%%%%%%%%%
%
% 1: The Paper.
%
%%%%%%%%%%%%%%%%%%%%%%%%%%%%%%%%%%%%%%%%%%%%%%%%%%%%%%%%%%%%%%%%%%%%%%%%%%%%%%%%%%%%%%%%%%%%%%%%%%%%%%%%%%%%%%%%%%%%%%%
%
\def\Pagetitle{}
\def\Pagefooter{\hfil{\myfontdefault\folio}\hfil}
\makeop{CN}%
\makeop{N}%
\makeop{Ln}%
\makeop{Sect}%
\makeop{Inj}%
\null \vfill
\def\centre{\rightskip=0pt plus 1fil \leftskip=0pt plus 1fil \spaceskip=.3333em \xspaceskip=.5em \parfillskip=0em \parindent=0em}%
\def\textmonth#1{\ifcase#1\or January\or Febuary\or March\or April\or May\or June\or July\or August\or September\or October\or November\or December\fi}
\font\abstracttitlefont=cmr10 at 14pt {\abstracttitlefont\centre
The Plateau problem for convex curvature functions\par}
\bigskip
{\centre 13th February 2014\par}
%{\centre \the\day\ \textmonth\month\ \the\year\par}
\bigskip
{\centre Graham Smith\par}
\bigskip
{\centre Instituto de Matem\'atica,\par
UFRJ, Av. Athos da Silveira Ramos 149,\par
Centro de Tecnologia - Bloco C,\par
Cidade Universit\'aria - Ilha do Fund\~ao,\par
Caixa Postal 68530, 21941-909,\par
Rio de Janeiro, RJ - BRASIL\par}
\bigskip
\noindent{\bf Abstract:\ }We present a novel and comprehensive approach to the study of the parametric Plateau problem for locally strictly convex (LSC) hypersurfaces of prescribed curvature for general convex curvature functions inside general Riemannian manifolds. We prove existence of solutions to the Plateau problem with outer barrier for LSC hypersurfaces of constant or prescribed curvature for general curvature functions inside general Hadamard manifolds modulo a single scalar condition. In particular, convex curvature functions of bounded type are fully treated.
\medskip
\noindent (First published on arXiv on 13th October 2010, this paper has been extensively redrafted to incorporate notational improvements made to \cite{SmiPPH}).
%
% 58E12 : Global analysis, analysis on manifolds; applications to minimal surfaces.
% 35J25 : Boundary value problems for second order elliptic equations
% 35J60 : Nonlinear elliptic equations
% 53A10 : Minimal surfaces, surfaces with prescribed mean curvature
% 53C21 : Methods of Riemannian geometry, including PDE methods; curvature restrictions
% 53C42 : Immersions (minimal, prescribed curvature, tight, etc.)
%
\par
\vfill
\nextoddpage
\global\pageno=1
\myfontdefault
\def\Pagetitle{\sl The Plateau problem for convex curvature functions}
\newhead{Introduction}
%
\newsubhead{Non-linear curvature functions} In the classical theory of hypersurfaces one usually studies functions of the principal curvatures, such as their mean or product which yield respectively the mean curvature and extrinsic curvature of the hypersurface. However, many other such ``curvature functions'' have also been formulated and studied. In order to correctly appreciate the results of this paper and their historical context, it is worthwhile to spend some time recalling the general framework of curvature functions as laid down elegantly by Caffarelli, Nirenberg \& Spruck in \cite{CaffNirSprIII} and \cite{CaffNirSprV}.
\medskip
\noindent Let $\Lambda_+:=\Lambda_+^n\subseteq\Bbb{R}^n$ be the open cone of vectors all of whose components are positive. Let $\Lambda:=\Lambda^n\subseteq\Bbb{R}^n$ be another open cone such that:
\subheadlabel{SubheadNonLinearCurvatureFunctions}
\medskip
\myitem{$(1)$} for all $x:=(x_1,...,x_n)\in\Lambda$ and for every permutation $\sigma$, $x_\sigma:=(x_{\sigma(1)},...,x_{\sigma(n)})$ is also an element of $\Lambda$;
\medskip
\myitem{$(2)$} for all $x\in\Lambda$ and for all $t>0$, $t x$ is also an element of $\Lambda$ (in particular, the vertex of $\Lambda$ lies on the origin);
\medskip
\myitem{$(3)$} for all $x\in\Lambda$ and for all $y\in\Lambda_+$, $x+y$ is also an element of $\Lambda$; and
\medskip
\myitem{$(4)$} $\Lambda$ is convex.
\medskip
\noindent A {\bf non-linear curvature function} is a non-negative function $K\in C^\infty(\Lambda)\minter C^0(\overline{\Lambda})$ such that:
\medskip
\myitem{$(1)$} $K$ is {\bf invariant} in the sense that for all $x\in\Lambda$ and for every permutation $\sigma$, $K(x_\sigma)=K(x)$;
\medskip
\myitem{$(2)$} $K$ is {\bf homogeneous} of order $1$;
\medskip
\myitem{$(3)$} $K$ is {\bf normalised} in the sense that $K(1,...,1)=1$; and
\medskip
\myitem{$(4)$} $K$ is {\bf compatible} with $\Lambda$ in the sense that it is strictly positive over $\Lambda$ and vanishes over $\partial\Lambda$.
\medskip
\noindent Scalar notions of curvature for suitable classes of hypersurfaces are defined using non-linear curvature functions. Indeed, let $M:=M^{n+1}$ be an $(n+1)$-dimensional Riemannian manifold and let $\Sigma:=(i,(S,\partial S))$ be a smooth {\bf immersed hypersurface} in $M$. We recall that this means that $(S,\partial S)$ is an $n$-dimensional manifold (possibly with boundary) and $i:S\rightarrow M$ is a smooth immersion. Let $\kappa:=(\kappa_1,...,\kappa_n)$ be its vector of principal curvatures. We say that $\Sigma$ is {\bf strictly} $K$-{\bf convex} whenever $\kappa(p)$ is an element of $\Lambda$ for all $p$ in $S$. When $\Sigma$ is strictly $K$-convex, we define the $K$-{\bf curvature} of $\Sigma$, $K(\Sigma):S\rightarrow\Bbb]0,\infty[$, by:
$$
K(\Sigma) = K(\kappa_1,...,\kappa_n).
$$
\noindent We now see that invariance ensures that both strict $K$-convexity and $K$-curvature are well-defined; homogeneity ensures that $K$-curvature transforms in a familiar manner under rescaling of the metric over $M$; and the normalisation condition ensures that the $K$-curvature of a unit sphere in Euclidean space is equal to $1$. Finally, compatibility ensures that no sequence $(x_n)_\ninn$ of points in $\Lambda$ with $(K(x_n))_\ninn$ uniformly bounded below may have a limit point on $\partial\Lambda$. In geometric terms, this means that in order for a smooth limit of strictly $K$-convex hypersurfaces to be also strictly $K$-convex, it is sufficient to ensure that the $K$-curvatures of all hypersurfaces in the sequence are uniformly bounded below. In other words, strict $K$-convexity is a closed property, modulo a scalar condition that is easy to verify.
\medskip
\noindent For PDE reasons, it is usual to require that a curvature function, $K$, satisfy, in addition, the following two conditions:
\medskip
\myitem{$(5)$} $K$ is {\bf strictly elliptic} in the sense that for all $x\in\Lambda$ and for all $1\leqslant i\leqslant n$, $\partial_ i K(x) > 0$; and
\medskip
\myitem{$(6)$} $K$ is a {\bf concave} function over $\Lambda$.
\medskip
\noindent We recall that the Jacobi operator of $K$-curvature over any strictly $K$-convex hypersurface measures the infinitesimal variation of the $K$-curvature arising from an infinitesimal normal perturbation of the immersion.\footnote*{More precisely, let $\Sigma:=(i,S)$ be an immersed hypersurface in $M$. Let $N$ be the upward-pointing unit normal vector field over $M$. For $\phi\in C_0^\infty(S)$, we define the family $(i_t)_{t\in]-\epsilon,\epsilon[}$ of immersions from $S$ into $M$ by $i_t(p)=\opExp(t\phi(p)N(p))$, where $\opExp$ is the exponential map of $M$. For all $t$, let $K_t:S\rightarrow\Bbb{R}$ be such that $K_t(p)$ is the $K$-curvature of $i_t$ at $p$. The Jacobi operator $J$ of $K$-curvature over $\Sigma$ is defined such that for all such $\phi$, $J\phi=\partial_t K_t|_{t=0}$.} The Jacobi operator is always a second order partial differential operator. Ellipticity of $K$ ensures that the Jacobi operator is elliptic for every smooth strictly $K$-convex immersed hypersurface. Finally, the condition of concavity is natural from the perspective of non-linear elliptic PDEs, and unavoidable at the current level of technology (c.f. \cite{CaffNirSprII}, \cite{CaffNirSprIII}).
\medskip
{\bf\noindent Examples:\ }Setting $\Lambda=\left\{x|x_1+...+x_n>0\right\}$ and $K(x) = (x_1 + ... + x_n)/n$, we recover the notions of local strict mean convexity and mean curvature. Setting $\Lambda=\Lambda_+=\left\{x|x_i>0\ \forall\ i\right\}$ and $K(x) = (x_1\cdot...\cdot x_n)^{1/n}$, we recover the notions of local strict convexity and extrinsic curvature. More generally, for $1\leqslant k< n$, define $f_k:\Bbb{R}^n\rightarrow\Bbb{R}$ by:
$$
f_k(x) = \frac{1}{n!}\sum_{\text{permutations}\ \sigma}x_{\sigma(1)}\cdot...\cdot x_{\sigma(k)},
$$
\noindent Define $\Lambda_k:=\left\{x\ |\ f_1(x),...,f_k(x)>0\right\}$ and define $\sigma_k:\Lambda_k\rightarrow]0,\infty[$ by $\sigma_k(x):=f_k(x)^{1/k}$. With this definition of $(\Lambda_k,\sigma_k)$, we recover the notions of local strict $k$-convexity and $\sigma_k$-curvature. In this notation, $\sigma_1$-curvature is mean curvature, $\sigma_n$-curvature is extrinsic curvature, and, when the ambient space is $\Bbb{R}^n$, $\sigma_2$-curvature coincides with the square root of the scalar curvature of the immersed hypersurface. \qed
\medskip
\noindent $K$ is said to be a {\bf convex curvature function} whenever $\Lambda=\Lambda_+$. We say that a smooth immersed hypersurface is {\bf locally strictly convex} (LSC) whenever it is strictly $\Lambda_+$-convex. Convex curvature functions are of particular interest as they are considerably more tractable to mathematical analysis than more general non-linear curvature functions. Indeed, whilst very general Plateau problems can be solved for convex curvature functions (in the literature discussed below, as well as in the current paper), comparatively little is known about Plateau problems for non-convex non-linear curvature functions when the desired hypersurface is anything other than a graph. In this paper, therefore, {\sl we will only treat smooth LSC hypersurfaces and convex curvature functions}\/ and we henceforth write $\Lambda$ instead of $\Lambda_+$.
\medskip
{\bf\noindent Examples:\ }Extrinsic curvature $K:=\sigma_n$ is a convex curvature function. For all $1\leqslant k<n$, the curvature quotient $K:=\sigma_{n,k}:=\sigma_n/\sigma_k$ is a convex curvature function (c.f. Proposition \procref{PropCurvatureQuotients}, below). Special Lagrangian curvature $K:=\rho_{(n-1)\pi/2}$ is a convex curvature function (c.f. \cite{SmiSLC} and \cite{SmiNLD}). If $K_1$ and $K_2$ are both convex curvature functions, then any {\bf convex sum} or {\bf convex product} of $K_1$ and $K_2$ given respectively by $\alpha K_1 + (1-\alpha)K_2$ and $K_1^\alpha K_2^{1-\alpha}$ for $\alpha\in[0,1]$ is also a convex curvature function. Any non-trivial convex product of a convex curvature function with any other curvature function is also a convex curvature function.\qed
\newsubhead{The non-linear Plateau problem} Let $M:=M^{n+1}$ be a complete $(n+1)$-dimensional Riemannian manifold. Let $K$ be a convex curvature function. Let $\Gamma:=(i,G)$ be a smooth compact codimension-$2$ immersed submanifold in $M$. The (convex) {\bf non-linear Plateau problem} asks for the existence of a smooth compact LSC immersed hypersurface $\Sigma:=(i,(S,\partial S))$ in $M$ with boundary equal to $\Gamma$ and with $K$-curvature, $K(\Sigma)(p)$ equal to some positive constant, $\kappa$, say, for all $p\in S$. More generally, for a given smooth function $\kappa:M\rightarrow]0,\infty[$, we ask for the existence of a smooth compact LSC immersed hypersurface $\Sigma$ such that $\partial\Sigma=\Gamma$ and $K(\Sigma)(p)=(\kappa\circ i)(p)$ for all $p\in S$. When this latter condition is satisfied, we say that the $K$-curvature of $\Sigma$ is {\bf prescribed} by $\kappa$, and we write $K(\Sigma)=\kappa$.
\subheadlabel{SubheadTheNonLinearPlateauProblem}
\medskip
\noindent The non-linear Plateau problem differs from the linear case in that $\Gamma$ is usually required to be the boundary of a $K$-convex outer barrier (which is the geometric analogue of the PDE concept of supersolution), $\hat{\Sigma}$, with $K$-curvature at all points greater than $k$. Furthermore, in treating the non-linear Plateau problem, it becomes necessary to separate the family of convex curvature functions into two classes as follows. For $K:\Lambda^n\rightarrow]0,\infty[$ a convex curvature function, following Trudinger (c.f. \cite{Trud}), we define $K_\infty:\Lambda^{n-1}\rightarrow]0,\infty]$ by:
$$
K_\infty(x_1,...,x_{n-1}) = \mlim_{t\rightarrow\infty} K(x_1,...,x_{n-1},t).
$$
\noindent Since $K$ is concave, $K_\infty$ is either everywhere infinite or everywhere finite and concave (c.f. Proposition \procref{PropLimitOfCurvFns}, below). We say that $K$ is of {\bf unbounded type} in the former case and of {\bf bounded type} in the latter. Convex curvature functions of bounded type present an extra level of analytic complexity. They are much less studied in the literature.
\medskip
{\bf\noindent Examples:\ }Extrinsic curvature $K:=\sigma_n$ is of unbounded type. The curvature quotients $K:=\sigma_{n,k}$ discussed in the preceeding section are of bounded type (c.f. Proposition \procref{PropCurvatureQuotients}, below). Special Lagrangrian curvature $K:=\rho_{(n-1)\pi/2}$ is of bounded type (c.f. \cite{SmiSLC} and \cite{SmiNLD}). Any convex sum or product of two convex curvature functions of bounded or unbounded type is also respectively of bounded or unbounded type. Any non-trivial convex sum or product of a convex curvature function of bounded type and a convex curvature function of unbounded type is of unbounded type.\qed
\medskip
\noindent An extensive literature has grown around the non-linear Plateau problem with outer barrier over the past three decades. The subject was essentially opened to modern mathematical analysis in the 1980's with the application by Caffarelli, Nirenberg \& Spruck in \cite{CaffNirSprV} of the celebrated barrier technique first introduced to the study of non-linear boundary value problems by the same authors in \cite{CaffNirSprI}. This allowed them to prove existence of solutions to the non-linear Plateau problems with outer barrier in $\Bbb{R}^{n+1}$ for (not necessarily convex) non-linear curvature functions of unbounded type. Importantly, however, they only treat the case where $\Gamma$ is the boundary of a strictly convex subset of a hyperplane and where the solution is a graph over this hyperplane. This is referred to as the {\bf non-parametric} case, as a specific parametrisation of the solution is provided by the geometry of the problem. Much subsequent work aimed to remove this restrictive geometric condition and thus treat the so-called {\bf parametric} case. The most notable result of this period was the existence theorem \cite{GuanSpruckA} of Guan \& Spruck for solutions to the non-linear Plateau problem with outer barrier in $\Bbb{R}^{n+1}$ for LSC graphs over hyperspheres. Spruck then conjectured in his presentation to the 1994 ICM (c.f. \cite{Spruck}) the existence of solutions to the parametric non-linear Plateau problem with outer barrier for hypersurfaces of constant extrinsic curvature in $\Bbb{R}^{n+1}$. This conjecture was solved simultaneously using identical techniques by Guan \& Spruck in \cite{GuanSpruck} and Trudinger \& Wang in \cite{TrudWang}. It was then extended by Guan \& Spruck in \cite{GuanSpruckI} to include all convex curvature functions of unbounded type, and by Sheng, Urbas \& Wang in \cite{ShengUrbasWang} to include the curvature quotients $K:=\sigma_{n,k}$ discussed in the preceeding section. The case of general convex curvature functions has remained to date unstudied.
\medskip
\noindent The recent results of Guan \& Spruck, Trudinger \& Wang and subsequent authors are obtained by supplementing Caffarelli, Nirenberg \& Spruck's barrier technique with a Perron-type argument. This beautiful approach suffers nonetheless from two significant limitations. First, the Perron method yields no information concerning the uniqueness of the solutions obtained, and second, and more significantly, by requiring in a fundamental manner the existence of large families of complete, totally geodesic hypersurfaces, it cannot be applied when the ambient manifold is anything other than a space-form except in certain very special cases where extra properties of the curvature function used compensate for this limitation (as is the case of special Lagrangian curvature, c.f. \cite{SmiSLC} and \cite{SmiNLP}). Nonetheless, in a different direction, using completely different techniques, Labourie showed in \cite{Lab} the existence of unique solutions to the non-linear Plateau problem with outer barrier for surfaces of constant extrinsic curvature immersed inside $3$-dimensional Hadamard manifolds of sectional curvature bounded above by $-1$. However, these techniques, involving an ingenious application of Gromov's theory of pseudo-holomorphic curves, do not extend to higher-dimensional ambient spaces.
\medskip
\noindent In summary, the existing results to date leave open the following three problems, in order of significance: firstly, to extend the result of Guan \& Spruck, Trudinger \& Wang and Labourie to general ambient spaces of arbitrary dimension; secondly, to extend these results to general convex curvatures functions, including the case of curvature functions of bounded type; and, thirdly, to prove the uniqueness or otherwise of the solutions obtained. In \cite{SmiCGC} and \cite{SmiPPH}, with the aim of addressing these open problems, we initiated a programme of developing a fully geometric approach to the original barrier technique of Caffarelli, Nirenberg \& Spruck as presented in \cite{CaffNirSprIII} and \cite{CaffNirSprV}. By developing furthermore a parametric Smale-type degree theory in the spirit of the work \cite{WhiteI} and \cite{WhiteII} of White, we extend the results of Guan \& Spruck, Trudinger \& Wang and Labourie to solve the parametric Plateau problem with outer barrer for hypersurfaces of prescribed extrinsic curvature inside general Hadamard manifolds of arbitrary dimension. Furthermore, although we do not, strictly speaking, prove uniqueness, the differential topological degree that we obtain nonetheless guarantees that under generic conditions, the number of solutions {\sl counted with algebraic sign} is equal to $1$.
\medskip
\noindent The current paper concerns the final stage of this programme which involves extending the existence result of \cite{SmiPPH} to the parametric non-linear Plateau problem with outer barrier for general convex curvature functions inside general manifolds. To this end, we prove a-priori estimates for the norms of the shape operators of hypersurfaces of prescribed $K$-curvature. The challenge here is two-fold, but in each case involves a highly non-trivial reworking of Caffarelli, Nirenberg \& Spruck's barrier arguments. The first challenge lies in the fact that the barrier functions used to date, being non-geometric in nature, become unworkably complex in general ambient spaces. In the case of extrinsic curvature, the multiplicative nature of the determinant function allows us to remove most of the problematic terms. However, for general convex curvature functions, we no longer have this luxury. We found that the best approach involves a complete reconstruction of Caffarelli, Nirenberg \& Spruck's original barrier functions in terms of natural geometric concepts. We find the resulting barrier arguments very satisfying in their generality as well as their relative simplicity. The second challenge concerns convex curvature functions of bounded type which have not hitherto been treated in full generality. However, by once again reconstructing in terms of natural geometric objects Caffarelli, Nirenberg \& Spruck's original barrier argument for the non-linear {\sl Dirichlet} problem (c.f. \cite{CaffNirSprIII}), we successfully eliminate the extra complexities that arise in this case, thus solving, in a relatively straightforward manner, the non-linear Plateau problem for all convex curvature functions inside general Hadamard manifolds, modulo a single scalar condition on the curvature function in question. It is an interesting open question to know to what extent this condition is necessary.
\medskip
\noindent Finally, the techniques used here do not restrict to Hadamard manifolds. In forthcoming work, we shall review the straightforward geometric conditions required to extend these results to Plateau problems inside suitable open subsets of general Riemannian manifolds.
\newsubhead{Main results}Given a convex curvature function, $K$, we define $\mu_\infty(K)$ by:
$$
\mu_\infty(K) = \mliminf_{\|x\|=1,x\rightarrow\partial\Gamma}DK_x(1,...,1).\eqnum{\nexteqnno}
$$
\noindent For all $K$, $\mu_\infty(K)>1$ (c.f. Proposition \procref{PropMuInfinityIsGreaterThanOne}, below). Furthermore, when $K$ is of unbounded type, $\mu_\infty(K)=\infty$ (c.f. Proposition \procref{PropMuInfinityIsInfinite}, below).
\medskip
\noindent Given a smooth compact codimension-$2$ immersed submanifold $\Gamma:=(i,G)$ in $M$, we say that $\Gamma$ is {\bf generic} whenever $T_p\Gamma\neq T_q\Gamma$ for all distinct pairs of points $p$ and $q$ in $G$. Observe that this condition is weaker than transversality. Our main result is:
\proclaim{Theorem \nextprocno}
\noindent Let $K$ be a convex curvature function; let $M$ be an $(n+1)$-dimensional Hadamard manifold; let $\kappa:M\rightarrow]0,\infty[$ be a smooth function; and let $\hat{\Sigma}$ be a smooth compact LSC immersed hypersurface in $M$. Suppose that:
\medskip
\myitem{$(1)$} $K(\hat{\Sigma})>\kappa$; and
\medskip
\myitem{$(2)$} $\partial\hat{\Sigma}$ is generic.
\medskip
\noindent Suppose, furthermore, that there exists a point $p\in M$ and $R>0$ such that:
\medskip
\myitem{$(3)$} $\hat{\Sigma}\subseteq B_R(p)$; and
\medskip
\myitem{$(4)$} $\kappa(q)<\frac{1}{R}\mu_\infty(K)$ for all $q\in B_R(p)$.
\medskip
\noindent There exists a smooth compact LSC immersed hypersurface $\Sigma$ in $M$ such that:
\medskip
\myitem{$(1)$} $\partial\Sigma=\partial\hat{\Sigma}$;
\medskip
\myitem{$(2)$} $\Sigma<\hat{\Sigma}$; and
\medskip
\myitem{$(3)$} $K(\Sigma)=\kappa$.
\endproclaim
\proclabel{ThmFirstExistenceTheorem}
\remark In particular, when $K$ is of unbounded type, Conditions $(3)$ and $(4)$ are trivially  satisfied for all $p$ and for all sufficiently large $R$.\qed
\medskip
\remark The notation $\Sigma<\hat{\Sigma}$ is explained in Section \subheadref{SubheadDefinitionOfConvexCobordisms} below. In particular, $\Sigma<\hat{\Sigma}$ implies that $\partial\Sigma=\partial\hat{\Sigma}$.\qed
\medskip
\remark The fact that $\mu_\infty(K)>1$ has the following pleasing consequence. Let $B_r(0)$ be the ball of radius $r$ about the origin in $\Bbb{R}^{n+1}$. Let $\hat{\Sigma}$ be an open subset in $\partial B_r(0)$ with smooth boundary. For all $\kappa\in]0,1/r]$ let $\Cal{S}_\kappa$ be the family of smooth compact LSC immersed hypersurfaces $\Sigma$ in $B_r(0)$ such that $K(\Sigma)=\kappa$ and $\Sigma<\hat{\Sigma}$. By Proposition \procref{PropCurvFnsII}, below, $\mu_\infty(K)\geqslant 1$. It follows that $\Cal{S}_\kappa$ varies continuously with $\kappa\in]0,1/r[$. Since $\mu_\infty(K)$ is in fact greater than $1$, it follows from the compactness results of this paper that this family also various continuously as $\kappa$ tends to $1/r$. Furthermore, using standard barrier techniques, we show that $\hat{\Sigma}$ is the only element of $\Cal{S}_{1/r}$, so that $(\Cal{S}_\kappa)_{\kappa\in]0,1/r[}$ converges to $\left\{\hat{\Sigma}\right\}$ as $\kappa$ tends to $1/r$. That is to say, the solutions approach that data as $\kappa$ tends to $1/r$. This would not necessarily hold if we did not have strict inequality in the above relation.\qed
\medskip
\noindent Using a slightly different approach, we obtain the following complementary result:
\proclaim{Theorem \nextprocno}
\noindent Let $K$ be a convex curvature function; let $M$ be an $(n+1)$-dimensional Hadamard manifold of sectional curvature bounded above by $-1$; let $\kappa:M\rightarrow]0,1[$ be a smooth function; and let $\hat{\Sigma}$ be a smooth compact LSC immersed hypersurface in $M$. Suppose that:
\medskip
\myitem{$(1)$} $K(\hat{\Sigma})>\kappa$; and
\medskip
\myitem{$(2)$} $\partial\hat{\Sigma}$ is generic.
\medskip
\noindent There exists a smooth compact LSC immersed hypersurface $\Sigma$ in $M$ such that:
\medskip
\myitem{$(1)$} $\partial\Sigma=\partial\hat{\Sigma}$;
\medskip
\myitem{$(2)$} $\Sigma<\hat{\Sigma}$; and
\medskip
\myitem{$(3)$} $K(\Sigma)=\kappa$.
\endproclaim
\proclabel{ThmSecondExistenceTheorem}
\noindent Concerning uniqueness, as indicated above, the degree theory of \cite{SmiDTI} should readily adapt to the current framework of compact hypersurfaces with boundary. It would then follow that the number of solutions (counted algebraically) would be equal to $1$ for generic $\hat{\Sigma}$. In any case, in certain situations uniqueness can be obtained, essentially because the contribution of every solution to the degree is positive, although this perspective is not necessary to obtain the result:
\proclaim{Theorem \nextprocno}
\noindent Let $K$ be a convex curvature function. Suppose that, for all $x:=(x_1,...,x_n)\in\Lambda$ such that $K(x_1,...,x_n)\leqslant 1$ the derivative of $K$ satisfies:
$$
DK_x\cdot(1,...,1)\geqslant DK_x\cdot(x_1^2,...,x_2^2).
$$
\noindent Then the solution obtained in Theorem \procref{ThmSecondExistenceTheorem} is unique.
\endproclaim
\proclabel{ThmDegreeUniqueness}
\remark This property is satisfied by the curvature quotients $K:=\sigma_{n,n-1}$ and $K:=\sigma_{n,n-2}$ (c.f. \cite{GuanSpruckIV}), by special Lagrangian curvature (c.f. Lemma $7.4$ of \cite{SmiFCS}) and by Gaussian curvature when the hypersurface is $2$-dimensional (c.f. Proposition $3.2.1$ of \cite{Lab}).\qed
\newsubhead{Conventions and acknowledgements} All immersed submanifolds will be oriented. For an immersed hypersurface $\Sigma$, we say that a unit normal vector field $N$ over $\Sigma$ points {\bf upward} whenever it is compatible with the orientation of $\Sigma$. We say it points {\bf downward} otherwise. If $\Sigma$ is LSC and lies on the boundary of a convex set, then the orientation of $\Sigma$ is chosen such that the convex set lies below $\Sigma$. If $\Gamma$ is the boundary of $\Sigma$, then the {\bf right} hand side of $\Gamma$ is the side on which $\Sigma$ lies, and the {\bf left} hand side is the other side.
\medskip
\noindent The first draft of this paper was written in response to a question asked by Joel Spruck at the conference, ``Algebraic, Geometric and Analytical Aspects of Surface Theory'', held in Buzios, Brazil in April 2010. The main results of this paper first appeared in their current form - with slightly stronger hypotheses - in the second draft which has been available on arXiv since October 2010. The author is grateful to the IMPA, Rio de Janeiro, Brazil, for providing the excellent conditions required to write the second draft. Finally, Sylvestre Gallot, Lucio Rodriguez and Harold Rosenberg are warmly thanked for helpful comments made to this and earlier drafts, as well as for their continual encouragement.
\newhead{Convex Curvature Functions}
\newsubhead{Convex curvature functions}We recall the various properties of convex curvature functions which are required throughout the sequel. Proofs of most of these relatively standard results are scattered throughout the literature. We collect most of them here for the reader's convenience.
\headlabel{HeadConvexCurvatureFunctions}
\proclaim{Proposition \nextprocno}
\noindent If $K$ is a convex curvature function, then $K_\infty$ is either everywhere infinite, or everywhere finite. Moreover, if $K_\infty$ is everywhere finite, then:
\medskip
\myitem{$(1)$} $K_\infty$ is invariant;
\medskip
\myitem{$(2)$} $K_\infty$ is homogeneous of order $1$;
\medskip
\myitem{$(3)$} $K_\infty(1,...,1)>1$;
\medskip
\myitem{$(4)$} $K_\infty$ is (possibly non-strictly) elliptic; and
\medskip
\myitem{$(5)$} $K_\infty$ is concave;
\medskip
\endproclaim
\proclabel{PropLimitOfCurvFns}
\remark Importantly, $K_\infty$ is not necessarily continuous over $\partial\Lambda^{n-1}$.\qed
\medskip
\proof The limit of an increasing family of concave functions over an open set is either everywhere infinite, or everywhere finite and concave. This proves the main assertion and $(5)$. $(1)$, $(2)$ and $(4)$ follow trivially. Finally $K_\infty(1,...,1)>K(1,...,1)=1$. $(3)$ follows, and this completes the proof.\qed
\medskip
\noindent For $1\leqslant k<n$, let $\sigma_{n,k}$ be the curvature quotient as defined in the introduction.
\proclaim{Proposition \nextprocno}
\noindent For $1\leqslant k< n$, $\sigma_{n,k}$ possesses Properties $(1)$ to $(6)$ and is of bounded type.
\endproclaim
\proclabel{PropCurvatureQuotients}
\proof Denote $K:=\sigma_{n,k}$. $(1)$ follows from the definition. $(2)$ and $(3)$ are trivial. Define $\Phi:]0,\infty[^n\rightarrow]0,\infty[^n$ by:
$$
\Phi(x_1,...,x_n) = (x^{-1},...,x^{-n}).
$$
\noindent Then:
$$
K\circ\Phi = \sigma_l^{-1/l},
$$
\noindent where $l=n-k$. Trivially, if any component of $x\in]0,\infty]^n$ is infinite, then so is $\sigma_l(x)$ and so $K$ extends to a continuous function over $\overline{\Lambda}$ which vanishes over $\partial\Lambda$. This proves $(4)$. For $x\in]0,\infty[^n$, if we denote $y=\Phi(x)$, then, for each $i$, by the chain rule:
$$
(\partial_i K)(y) = c_1 y_i^2 \sigma_l^{-(l+1)/l} \sigma_{l-1}(x_1,...,\hat{x}_i,...,x_n),
$$
\noindent for some positive constant $c_1$. This is trivially positive, and $(5)$ follows. $(6)$ is proven in \cite{CaffNirSprIII}. Choose $(y_1,..,y_n)\in]0,\infty[^{n-1}$. Trivially:
$$
\mlim_{t\rightarrow 0}\sigma_l(y_1,...,y_{n-1},t) = c_2 \sigma_l(y_1,...,y_{n-1}),
$$
\noindent for some positive constant $c_2$. It follows that $K_\infty = c_2 \sigma_{n-1,k-1}$. In particular, $K$ is of bounded type. This completes the proof.\qed
\medskip
\noindent We relate convex curvature functions to $O(n)$-invariant functions over the space of symmetric matrices. Let $\opSymm(n)$ be the space of real valued symmetric $n$-dimensional matrices and let $\Lambda:=\Lambda^n\subseteq\opSymm(n)$ now be the open convex cone of positive definite matrices. Let $f$ be a convex curvature function. We define $F\in C^\infty(\Lambda)\minter C^0(\overline{\Lambda})$ by:
$$
F(A) = f(\lambda_1,...,\lambda_n),
$$
\noindent where $\lambda_1,...,\lambda_n$ are the eigenvalues of $A$. Interpreting $F$ in this manner, yields the following important characterisation of $F_\infty$:
\proclaim{Proposition \nextprocno}
\noindent For all $A\in\Lambda^{n-1}$:
$$
F_\infty(A) = \msup\left\{F(B)\ |\ B\in\Lambda^n,\ B|_{\Bbb{R}^{n-1}}=A\right\}.
$$
\noindent Furthermore, this supremum is not attained by any $B\in\Lambda^n$.

\endproclaim
\proclabel{PropCharacterisationOfFInfinity}
\proof We denote the supremum by $F$. Trivially, $F_\infty(A)\leqslant F$. Conversely, let $B\in\Lambda^n$ be such that its restriction to $\Bbb{R}^{n-1}$ concides with $A$. Let $\lambda_1\leqslant...\leqslant\lambda_{n-1}$ be the eigenvalues of $A$. Let $\mu_1\leqslant...\leqslant\mu_n$ be the eigenvalues of $B$. By the classical minimax principle for eigenvalues of a real valued symmetric matrix (c.f., for example, Lemma $10.2$ of \cite{SmiFCS}), for all $1\leqslant\alpha\leqslant n-1$, $\mu_\alpha\leqslant\lambda_\alpha$. Thus, bearing in mind ellipticity of $f$
$$
F(B) = f(\mu_1,...,\mu_n) \leqslant f(\lambda_1,...,\lambda_{n-1},\mu_n) < f_\infty(\lambda_1,...,\lambda_{n-1}) = F_\infty(A).
$$
\noindent Taking the supremum over all $B$, it follows that $F_\infty(A)\geqslant F$ and that $F_\infty(A)$ is not attained by any element of $\Lambda^n$, as desired.\qed
\medskip
\noindent We finally list the properties of $F$:
\proclaim{Proposition \nextprocno}
\noindent Let $f$ be a convex curvature function. For all $A\in\Lambda\subseteq\opSymm(n)$, there exists a unique matrix $B\in\opSymm(n)$ such that, for all $M\in\opSymm(n)$:
$$
DF_A(M) = \opTr(BM).
$$
\noindent Moreover:
\medskip
\myitem{$(1)$} $B$ is positive definite; and
\medskip
\myitem{$(2)$} $A$ and $B$ are simultanously diagonalisable.
\medskip
\noindent In addition, if $e_1,...,e_n$ is a system of shared eigenvectors for $A$ and $B$ and if $\lambda_1,...,\lambda_n$ and $\mu_1,...,\mu_n$ are the corresponding eigenvalues of $A$ and $B$ respectively, then:
\medskip
\myitem{$(3)$} for all $i\neq j$:
$$
\lambda_i\geqslant\lambda_j \Leftrightarrow \mu_i\leqslant\mu_j;
$$
\myitem{$(4)$} $DF_A(A)=\sum_{i=1}^n\lambda_i\mu_i=F(A)$;
\medskip
\myitem{$(5)$} for all $A'\in A$, $DF_A(A')\geqslant F(A')$;
\medskip
\myitem{$(6)$} $\sum_{i=1}^n\mu_i=DF_A(\opId)\geqslant 1$; and
\medskip
\myitem{$(7)$} for all $M\in\opSymm(n)$:
$$
-(D^2f)_A(M,M) \geqslant \sum_{i\neq j}^n\frac{\mu_j-\mu_i}{\lambda_i-\lambda_j}M_{ij}^2 \geqslant 0.
$$
\noindent In particular, $F$ is concave.
\endproclaim
\proclabel{PropCurvFnsII}
\proof $DF_A:\opSymm(n)\rightarrow\opSymm(n)$ is linear. There therefore exists a unique matrix $B\in\opSymm(n)$ such that, for all $M\in\opSymm(n)$:
$$
DF_A(M) = \opTr(BM).
$$
\noindent By invariance of $f$, $F$ is $O(n)$ invariant. Thus, for all $A\in\Lambda$ and $M\in O(n)$:
$$
F(M^tAM) = F(A).
$$
\noindent Differentiating, for all antisymmetric $M$:
$$\matrix
&DF_A(MA - AM) \hfill&=0\hfill\cr
\Rightarrow\hfill&\opTr([AB]M)\hfill&=0.\hfill\cr
\endmatrix$$
\noindent However, since $A$ and $B$ are both symmetric, $[AB]$ is antisymmetric, and so, since $M$ is arbitrary:
$$
[AB] = 0.
$$
\noindent $(2)$ now follows. Let $e_1,...,e_n$ be a system of shared eigenvectors of $A$ and $B$ and let $\lambda_1,...,\lambda_n$ and $\mu_1,...,\mu_n$ be the respective corresponding eigenvalues. By strict ellipticity of $f$, for all $i$:
$$
\mu_i = DF_A(e_i\otimes e_i) > 0.
$$
\noindent $(1)$ now follows. $(3)$ follows by concavity of $f$ and $O(n)$-invariance and $(4)$ follows by homogeneity. Likewise, by concavity, for all $A'\in\Lambda$:
$$
DF_A(A'-A) \geqslant F(A') - F(A).
$$
\noindent Thus, by $(4)$, $DF_A(A')\geqslant F(A')$, and $(5)$ follows. $(6)$ is a special case of $(5)$. Finally, $(7)$ follows by concavity as in Lemma $2.3$ of \cite{ShengUrbasWang}. This completes the proof.\qed
\medskip
\newhead{Convex Cobordisms and Embedding Radii}
\newsubhead{Definition of convex cobordisms}Convex cobordisms provide a concise means of expressing the outer barrier condition. They are defined as follows. Let $M:=M^{n+1}$ be an $(n+1)$-dimensional Hadamard manifold. For $m\in\left\{1,2\right\}$, let $\Sigma_m:=(i_m,(S_m,\partial S_m))$ be a smooth compact LSC immersed hypersurface in $M$. We say that $\Sigma_2$ {\bf bounds} $\Sigma_1$, and we denote $\Sigma_1<\Sigma_2$ whenever there exists a smooth compact $(n+1)$-dimensional manifold $(N,\partial N)$ with piecewise smooth boundary and a smooth immersion $I:N\rightarrow M$ such that $\partial N$ consists of $2$ smooth components $\partial_1N$ and $\partial_2 N$ and:
\headlabel{HeadConvexCobordismsAndEmbeddingRadii}
\subheadlabel{SubheadDefinitionOfConvexCobordisms}
\medskip
\myitem{$(1)$} the boundaries of $\partial_1 N$ and $\partial_2 N$ coincide;
\medskip
\myitem{$(2)$} $\partial_1 N$ and $\partial_2 N$ meet transversally along their shared boundary;
\medskip
\myitem{$(3)$} the immersed hypersurface $(\partial_1 N,I)$ coincides with $\Sigma_1$ and furthermore $N$ lies above this LSC hypersurface;
\medskip
\myitem{$(4)$} the immersed hypersurface $(\partial_2 N,I)$ coincides with $\Sigma_2$, and furthermore $N$ lies below this LSC hypersurface; and
\medskip
\myitem{$(5)$} furnishing $N$ with the unique metric that makes $I$ into a local isometry, $N$ is foliated by the geodesic segments normal to $\partial_1 N$.
\medskip
\noindent We refer to $(N,I)$ as the {\bf convex cobordism} from $\Sigma_1$ to $\Sigma_2$. When $(N,I)$ only satisfies conditions $(1)$ to $(4)$, we call it a {\bf convex precobordism} from $\Sigma_1$ to $\Sigma_2$. This weaker concept is useful in proving technical results.
\medskip
\noindent The concept of bounding generalises the concept of graph for LSC immersions. Indeed, we recall that $\Sigma_2$ is said to be a graph over $\Sigma_1$ whenever there exists a smooth non-negative function $f:S_1\rightarrow [0,\infty[$ and a smooth diffeomorphism $\alpha:S_1\rightarrow S_2$ such that $f$ vanishes along $\partial S_1$ and for all $p\in S_1$, $(i_2\circ\alpha)(p)=\opExp(f(p)N_1(p))$, where $N_1$ is the upward-pointing unit normal vector field over $\Sigma_1$. We define the manifold $N\subseteq S_1\times\Bbb{R}$ and the smooth immersion $I:N\rightarrow M$ by:
$$
N = \left\{ (p,t)\ |\ 0\leqslant t\leqslant f(p)\right\},\qquad I(p,t) = \opExp(tN_1(p)).
$$
\noindent $(N,I)$ defines a convex cobordism from $\Sigma_1$ to $\Sigma_2$. We conclude that $\Sigma_2$ bounds $\Sigma_1$ whenever $\Sigma_2$ is a graph over $\Sigma_1$. Conversely, one may show that if $\Sigma_2$ bounds $\Sigma_1$ and makes an angle of less than $\pi/2$ with $\Sigma_1$ at every point of their common boundary, then $\Sigma_2$ is a graph over $\Sigma_1$, although we will have no need for this result.
\medskip
\noindent In order for boundedness to be used as part of a Smale-type degree theory, we require that it constitute an open and closed condition modulo suitable hypotheses. Sections \subheadref{SubheadGeometryOfConvexPrecobordisms} to \subheadref{SubheadOpennessAndClosedness} are devoted to establishing these properties.
\newsubhead{Geometry of convex precobordisms}Let $\Sigma_1$ and $\Sigma_2$ be as before. Let $N:=(N,I)$ be a convex precobordism from $\Sigma_1$ to $\Sigma_2$. The geometry of $N$ is encoded in the following:
\subheadlabel{SubheadGeometryOfConvexPrecobordisms}
\proclaim{Proposition \nextprocno}
\noindent Choose $p\in N\setminus S_1$, let $q\in S_1$ be a point minimising distance to $p$ and let $\gamma:[0,1]\rightarrow N$ be a length minimising rectifiable curve such that $\gamma(0)=q$ and $\gamma(1)=p$. Then:
\medskip
\myitem{$(1)$} $\gamma(]0,1[)$ lies in the interior of $N$;
\medskip
\myitem{$(2)$} $\gamma(]0,1[)$ is a smooth geodesic;
\medskip
\myitem{$(3)$} if $p\in S_2\setminus\partial S_2$, then $\gamma$ is transverse to $S_2$ at $p$;
\medskip
\myitem{$(4)$} if $q\in S_1\setminus\partial S_1$, then $\gamma$ makes a right angle with $S_1$ at $q$; and
\medskip
\myitem{$(5)$} if $q\in\partial S_1$, then $\gamma$ is normal to $\partial S_1$ at $q$ and:
$$
\langle \partial_t\gamma(0),N_{\partial S_1}\rangle \geqslant 0,
$$
\noindent where $N_{\partial S_1}$ is the outward-pointing normal vector to $\partial S_1$ along $S_1$.
\endproclaim
\proclabel{PropGeometryOfCobordism}
\proof Define $t_0\in [0,1]$ by:
$$
t_0 = \minf\left\{t\in[0,1]\ |\ \gamma(t)\in S_2\right\}.
$$
\noindent If $q\in S_1\setminus\partial S_1$, then trivially $t_0>0$. If $q\in\partial S_1$, then likewise, by local strict convexity of $S_2$, $t_0>0$ unless $\gamma$ is trivial. Since $\gamma([0,t_0])$ lies inside $S_2$, by local strict convexity of $S_2$, it cannot be tangent to $S_2$ at $t_0$. It follows that $\gamma$ does not minimise length unless $t_0=1$ or $t_0=+\infty$ (the latter case occuring when $\gamma$ never intersects $S_2$) and this proves $(1)$. $(2)$ is trivial. Likewise, taking $t_0=1$, we see that $\gamma$ is transverse to $S_2$ at $\gamma(1)=p$, and this proves $(3)$. $(4)$ and $(5)$ are trivial, since $\gamma$ is length minimising. This completes the proof.\qed
\medskip
\noindent We characterise convex cobordisms amongst convex precobordisms:
\proclaim{Proposition \nextprocno}
\noindent If there exists no geodesic arc $\gamma:[0,1]\rightarrow N$ such that $\gamma(0)\in S_1$; $\gamma(1)\in S_1$; and $\gamma$ is normal to $S_1$ at $\gamma(0)$, then $(N,I)$ is a convex cobordism from $\Sigma_1$ to $\Sigma_2$.
\endproclaim
\proclabel{PropCharacterisationOfConvexCobordisms}
\proof For $p\in N$, let $Q(p)$ be the set of points $q\in S_1$ such that there exists a geodesic $\gamma_q:[0,1]\rightarrow N$ such that $\gamma_q(0)=q$; $\gamma_q$ is normal to $S_1$ at $q$; and $\gamma_q(1)=p$. Since $S_1$ is LSC, and since $N$ is non-positively curved, $Q(p)$ is discrete. Since $S_1$ is compact, $Q(p)$ is therefore finite. Denote $D(p)=\left|Q(p)\right|$. It suffices to show that $D(p)=1$ for all $p$. Choose $q\in Q(p)$. As in Proposition \procref{PropGeometryOfCobordism}, $(1)$, $\gamma_q$ may only intersect $S_2$ at $p$. Furthermore, if $p\in S_2$, then $\gamma_q$ is transverse to $S_2$ at this point. Finally, by hypothesis, $\gamma_q$ only intersects $S_1$ at $q$. It follows that for all $q\in Q(p)$, and for all $p'$ sufficiently close to $p$, $\gamma(q)$ may be perturbed to another geodesic segment normal to $S_1$ and terminating at $p'$. $D(p)$ is therefore locally constant. However, by hypothesis, $D(p)$ is equal to $1$ along $S_1$. It follows by connectedness that $D(p)=1$ for all $p$, as desired.\qed
\newsubhead{Geometry of convex cobordisms}Let $\Sigma_1$ and $\Sigma_2$ continue to be as before. Let $N:=(N,I)$ be a convex cobordism from $\Sigma_1$ to $\Sigma_2$.
\subheadlabel{SubheadGeometryOfConvexCobordisms}
\proclaim{Proposition \nextprocno, {\bf Uniqueness}}
\noindent If $(N',I')$ is another convex cobordism from $\Sigma_1$ to $\Sigma_2$, then there exists a diffeomorphism $\alpha:N\rightarrow N'$ such that $I=I'\circ\alpha$.
\endproclaim
\proclabel{PropUniquenessOfConvexCobordisms}
\proof Let $\Cal{F}$ and $\Cal{F}'$ be the foliations of $N$ and $N'$ respectively by geodesic segments normal to $S_1$. Observe that $\Cal{F}$ and $\Cal{F}'$ are canonically homeomorphic. We thus identify $\Cal{F}$ and $\Cal{F}'$. Observe that, as $\Cal{F}$ is also homeomorphic to $\Sigma_2$, it has the structure of a topological manifold with boundary. For $L\in\Cal{F}$, let $h(L)$ and $h'(L)$ be the length of the leaf $L$ in $N$ and $N'$ respectively. It suffices to show that $h=h'$ throughout $\Cal{F}$. Let $\Omega\subseteq\Cal{F}$ be the set of all points where $h=h'$. Observe that $\partial\Cal{F}\subseteq\Omega$. Let $\Omega_0$ be the connected component of $\Omega$ containing $\partial\Cal{F}$. $\Omega_0$ is open and, by continuity, it is also closed. It follows by connectedness that $\Omega_0\subseteq\Cal{F}$. The result follows.\qed
\medskip
\noindent Let $d:N\rightarrow\Bbb{R}$ be the distance in $N$ to $S_1$. The main consequence of the foliation condition is the following result:
\proclaim{Proposition \nextprocno}
\noindent For all $p\in N$, there exists a unique point $q\in S_1$ minimising distance to $p$.
\endproclaim
\proclabel{PropExistsUniqueMinimiser}
\proof Let $q\in S_1$ minimise distance to $p$. Let $\gamma:[0,1]\rightarrow N$ be a length minimising curve from $p$ to $q$. By Proposition \procref{PropGeometryOfCobordism}, $(2)$, $(4)$ and $(5)$, $\gamma$ is a geodesic in $N$ which is normal to $S_1$ at $q$. However, by hypothesis, there is only one such geodesic passing through $p$. The result follows.\qed
\proclaim{Proposition \nextprocno}
\noindent The function $d$ is strictly convex.
\endproclaim
\proclabel{PropDistInCobordismIsConvex}
\proof Choose $p\in N\setminus S_1$. By Proposition \procref{PropExistsUniqueMinimiser}, there exists a unique point $q\in S_1$ minimising distance to $p$. Since $S_1$ is LSC and since $N$ is non-positively curved, the result now follows.\qed
\medskip
\noindent Furthermore, Proposition \procref{PropExistsUniqueMinimiser} yields a well-defined closest-point projection $\pi:N\rightarrow S_1$. Since $S_1$ is LSC and since $N$ is non-positively curved, $\pi$ is distance decreasing. Considering the restriction of $\pi$ to $S_2$, we therefore obtain:
\proclaim{Proposition \nextprocno}
\noindent For each $i$, let $\opDiam(\Sigma_i)$ and $\opVol(\Sigma_i)$ denote the diameter and volume of $\Sigma_i$ respectively. Then:
$$
\opDiam(\Sigma_2)\geqslant\opDiam(\Sigma_1),\qquad \opVol(\Sigma_2)\geqslant\opVol(\Sigma_1).
$$
\endproclaim
\proclabel{PropDiameterBounds}
\newsubhead{Glueing and excision}For $m\in\left\{1,2,3\right\}$, let $\Sigma_m:=(i_m,(S_m,\partial S_m))$ be smooth compact LSC immersed hypersurfaces in $M$. Suppose that $\Sigma_1<\Sigma_2$ and $\Sigma_2<\Sigma_3$. Let $(N_{12},I_{12})$ and $(N_{23},I_{23})$ be convex cobordisms from $\Sigma_1$ to $\Sigma_2$ and from $\Sigma_2$ to $\Sigma_3$ respectively. We define $N_{12}\munion N_{23}$ by joining $N_{12}$ to $N_{23}$ along $S_2$. We define $I_{12}\munion I_{23}:N_{12}\munion N_{23}\rightarrow M$ by:
\subheadlabel{SubheadGlueingAndExcision}
$$
(I_{12}\munion I_{23})(x) = \left\{\matrix I_{12}(x)\text{ if }x\in N_{12};\text{ and }\hfill\cr
I_{23}(x)\text{ if }x\in N_{23}.\hfill\cr\endmatrix\right.
$$
\noindent $N_{12}\munion N_{23}$ is trivially a convex precobordism from $\Sigma_1$ to $\Sigma_3$.
\proclaim{Proposition \nextprocno}
\noindent There exists no non-trivial geodesic arc $\gamma:[0,1]\rightarrow N_{12}\munion N_{23}$ such that $\gamma(0),\gamma(1)\in S_1$.
\endproclaim
\proclabel{PropGeosDontReturn}
\proof Suppose the contrary. Let $\gamma:[0,1]\rightarrow N_{12}\munion N_{23}$ be a geodesic arc such that $\gamma(0),\gamma(1)\in S_1$. By Proposition \procref{PropDistInCobordismIsConvex}, $\gamma$ is not contained in $N_{12}$. There therefore exists $t_0\in]0,1[$ such that $\gamma(t_0)$ lies in the interior of $N_{23}$, and there exist $t_1<t_0<t_2$ such that $\gamma([t_1,t_2])$ is contained in $N_{23}$ and $\gamma(t_1)$ and $\gamma(t_2)$ both lie in $S_2$. However, this is also absurd by Proposition \procref{PropDistInCobordismIsConvex}, and this completes the proof.\qed
\medskip
\noindent We thus obtain a glueing operation for convex cobordisms:
\proclaim{Proposition \nextprocno, {\bf Glueing}}
\noindent $(N_{12}\munion N_{23},I_{12}\munion I_{23})$ defines a convex cobordism from $\Sigma_1$ to $\Sigma_3$.
\endproclaim
\proclabel{PropCobIsTrans}
\proof This follows from Propositions \procref{PropCharacterisationOfConvexCobordisms} and \procref{PropGeosDontReturn}.\qed
\medskip
\noindent Now, for $m\in\left\{1,2\right\}$, let $\Sigma_m:=(i_m,(S_m,\partial S_m))$ be a smooth compact LSC immersed hypersurface in $M$. Suppose that $\Sigma_1<\Sigma_2$. Let $(N,I)$ be the convex cobordism from $\Sigma_1$ to $\Sigma_2$. Denote by $\Cal{F}$ the foliation of $N$ by geodesic segments normal to $S_1$. Let $(S_3,\partial S_3)$ be a smooth compact LSC embedded submanifold of $N$ such that $\partial S_3=\partial S_1=\partial S_2$; $S_3$ is transverse to every leaf of $\Cal{F}$; and $S_3$ is not tangent to $S_1$ or to $S_2$ at any point.
\proclaim{Proposition \nextprocno}
\noindent $S_3$ meets every leaf of $\Cal{F}$ at exactly $1$ point.
\endproclaim
\proclabel{PropHowTheHypersurfaceMeetsTheFoliation}
\proof For $L$ in $\Cal{F}$, since $S_3$ is transverse to $L$, $L\minter S_3$ is discrete. Since $S_3$ is compact, $L\minter S_3$ is therefore finite. Denote $D(L)=\left|L\minter S_3\right|$. By transversality and finiteness, $D(L)$ is locally constant. However, close to $\partial S_1=\partial S_2$, $D(L)=1$. It follows by connectedness that $D(L)=1$ for all $L$, as desired.\qed
\proclaim{Proposition \nextprocno}
\noindent $S_3$ divides $N$ into two connected components.
\endproclaim
\proof This follows from Proposition \procref{PropHowTheHypersurfaceMeetsTheFoliation} since each leaf divides into two subintervals lying below and above $S_3$ respectively.\qed
\medskip
\noindent Let $i_3$ be the restriction of $I$ to $S_3$. Denote $\Sigma_3:=(i_3,(S_3,\partial S_3))$. Denote by $N_-$ and $N_+$ the closures of the connected components of $N\setminus S_3$ lying below and above $S_3$ respectively.
\proclaim{Proposition \nextprocno}
\noindent $(N_-,I)$ and $(N_+,I)$ define convex precobordisms from $\Sigma_1$ to $\Sigma_3$ and from $\Sigma_3$ to $\Sigma_2$ respectively.
\endproclaim
\proof It suffices to show that $S_1$ and $S_2$ lie below and above $S_3$ respectively. However, since this holds along the boundary, the result follows.\qed
\medskip
\noindent This yields the first excision operation for convex cobordisms:
\proclaim{Proposition \nextprocno, {\bf Excision I}}
\noindent $(N_-,I)$ defines a convex cobordism from $\Sigma_1$ to $\Sigma_3$.
\endproclaim
\proclabel{PropExcisionI}
\proof Since $S_3$ meets every leaf of $\Cal{F}$ at exactly $1$ point, the geodesic segments normal to $S_1$ also foliate $N_-$, and the result follows.\qed
\proclaim{Proposition \nextprocno}
\noindent There exists no geodesic segment $\gamma:[0,1]\rightarrow N_+$ such that $\gamma(0)\in S_3$; $\gamma(1)\in S_3$; and $\gamma$ is normal to $S_3$ at $\gamma(0)$.
\endproclaim
\proclabel{PropGeoDoesNotReturnII}
\proof Suppose the contrary. Let $\gamma$ be such a geodesic. Let $d:N\rightarrow[0,\infty[$ be the distance to $S_1$ along $\Cal{F}$. Since $\gamma$ points upwards from $S_2$, $\langle\dot{\gamma}(0),(\nabla d\circ\gamma)(0)\rangle > 0$. However, by Proposition \procref{PropDistInCobordismIsConvex}, $d$ is convex. Thus, since $\gamma$ is a geodesic:
$$
\partial_t\langle\dot{\gamma}(t),(\nabla d\circ\gamma)(t)\rangle = \opHess(d)(\dot{\gamma}(t),\dot{\gamma}(t)) \geqslant 0.
$$
\noindent In particular, $\langle\dot{\gamma}(1),(\nabla d\circ\gamma)(1)\rangle>0$. That is, $\gamma$ meets $S_2$ from below at $\gamma(1)$. In other words, for $t$ sufficiently close to $1$, $\gamma(1)\in N_-\setminus S_3$. This is absurd, and the result follows.\qed
\medskip
\noindent This yields the second excision operation for convex cobordisms:
\proclaim{Proposition \nextprocno, {\bf Excision II}}
\noindent $(N_+,I)$ defines a convex cobordism from $\Sigma_3$ to $\Sigma_2$.
\endproclaim
\proclabel{PropExcisionII}
\proof This follows from Propositions \procref{PropCharacterisationOfConvexCobordisms} and \procref{PropGeoDoesNotReturnII}.\qed
\newsubhead{Openness and closedness}Let $\Sigma_1^-$,$\Sigma_1$,$\Sigma_2$ and $\Sigma_2^+$ be smooth compact LSC immersed hypersurfaces in $M$ such that $\Sigma_1^-<\Sigma_1<\Sigma_2<\Sigma_2^+$. By Proposition \procref{PropCobIsTrans}, $\Sigma_1^-<\Sigma_2^+$. Let $N$ be the convex cobordism from $\Sigma_1^-$ to $\Sigma_2^+$. By Proposition \procref{PropUniquenessOfConvexCobordisms}, $N$ is unique. Observe that $\Sigma_1$ and $\Sigma_2$ identify with smooth compact LSC embedded hypersurfaces in $N$ which we denote by $S_1$ and $S_2$ respectively.
\subheadlabel{SubheadOpennessAndClosedness}
\proclaim{Proposition \nextprocno}
\noindent $S_1$ lies below $S_2$ and both $S_1$ and $S_2$ are transverse to the foliation of $N$ by geodesic segments normal to $S_1^-$.
\endproclaim
\proclabel{PropTransversalityOfInteriorHypersurfaces}
\proof The convex cobordism $N$ may be constructed by glueing the convex cobordism from $\Sigma_1$ to $\Sigma_2$ to the convex cobordism from $\Sigma_1^-$ to $\Sigma_1$ along $S_1$ and then glueing the convex cobordism form $\Sigma_2$ to $\Sigma_2^+$ to this convex convex cobordism along $S_2$. The result follows.\qed
\medskip
\noindent We now obtain openness:
\proclaim{Proposition \nextprocno, {\bf Openness}}
\noindent Let $(\Sigma_n)_\ninn$ and $(\hat{\Sigma}_n)_\ninn$ be sequences of smooth compact LSC immersed hypersurfaces in $M$ converging to the smooth compact LSC immersed hypersurfaces $\Sigma_\infty$ and $\hat{\Sigma}_\infty$ respectively. If $\Sigma_\infty<\hat{\Sigma}_\infty$, then for all sufficiently large $n$, $\Sigma_n<\hat{\Sigma}_n$.
\endproclaim
\proclabel{PropBoundednessIsAnOpenProperty}
\proof Upon perturbing $\Sigma_\infty$ and $\hat{\Sigma}_\infty$ we obtain smooth compact LSC immersed hypersurfaces $\Sigma_\infty'$ and $\hat{\Sigma}_\infty'$ respectively such that $\Sigma_\infty'<\Sigma_\infty$ and $\hat{\Sigma}_\infty<\hat{\Sigma}_\infty'$. By Proposition \procref{PropCobIsTrans}, $\Sigma'_\infty<\hat{\Sigma}'_\infty$. Let $N$ be the convex cobordism from $\Sigma'_\infty$ to $\hat{\Sigma}'_\infty$. Denote by $S_\infty'$ and $\hat{S}_\infty'$ the lower and upper boundary components of $N$ respectively. Denote by $\Cal{F}$ the foliation of $N$ by geodesic segments normal to $S_\infty'$. $\Sigma_\infty$ and $\hat{\Sigma}_\infty$ identify with smooth compact LSC embedded hypersurfaces in $N$ which we denote by $S_\infty$ and $\hat{S}_\infty$ respectively. By Proposition \procref{PropTransversalityOfInteriorHypersurfaces}, $S_\infty$ lies below $\hat{S}_\infty$ and both $S_\infty$ and $\hat{S}_\infty$ are transverse to $\Cal{F}$. For sufficiently large $n$, $\Sigma_n$ and $\hat{\Sigma}_n$ identify with smooth compact LSC embedded hypersurfaces in $N$ which we denote by $S_n$ and $\hat{S}_n$ respectively. Upon increasing $n$ further if necessary, we may suppose that $S_n$ lies below $\hat{S}_n$ and that both $S_n$ and $\hat{S}_n$ are transverse to $\Cal{F}$. By Proposition \procref{PropExcisionI}, $\Sigma_\infty'<\hat{\Sigma}_n$, and by Proposition \procref{PropExcisionII}, $\Sigma_n<\hat{\Sigma}_n$, as desired.\qed
\medskip
\noindent We obtain closedness modulo a condition on the curvature:
\proclaim{Proposition \nextprocno, {\bf Closedness}}
\noindent Let $(\Sigma_n)_\ninn$ and $(\hat{\Sigma}_n)_\ninn$ be sequences of smooth compact LSC immersed hypersurfaces in $M$ converging to the smooth compact LSC immersed hypersurfaces $\Sigma_\infty$ and $\hat{\Sigma}_\infty$ respectively. If $\Sigma_n<\hat{\Sigma}_n$ for all $n$, and if $K(\Sigma_\infty)<K(\hat{\Sigma}_\infty)$, then $\Sigma_\infty<\hat{\Sigma}_\infty$.
\endproclaim
\proclabel{PropBoundednessIsAClosedProperty}
\proof Upon perturbing $\Sigma_\infty$ and $\hat{\Sigma}_\infty$ we obtain smooth compact LSC immersed hypersurfaces $\Sigma_\infty'$ and $\hat{\Sigma}_\infty'$ respectively such that $\Sigma_\infty'<\Sigma_\infty$ and $\hat{\Sigma}_\infty<\hat{\Sigma}_\infty'$. By Proposition \procref{PropBoundednessIsAnOpenProperty}, there exists $n_0$ such that for all $n\geqslant n_0$, $\Sigma_\infty'<\Sigma_n$ and $\hat{\Sigma}_n<\hat{\Sigma}_\infty'$. Thus, by Proposition \procref{PropCobIsTrans}, $\Sigma_\infty'<\hat{\Sigma}_\infty'$. Let $N$ be the convex cobordism from $\Sigma'_\infty$ to $\hat{\Sigma}'_\infty$. Denote by $S_\infty'$ and $\hat{S}_\infty'$ the lower and upper boundary components of $N$ respectively. Denote by $\Cal{F}$ the foliation of $N$ by geodesic segments normal to $S_\infty'$.  By uniqueness of $N$, for $n\geqslant n_0$, $\Sigma_n$ and $\hat{\Sigma}_n$ identify with smooth compact LSC embedded hypersurfaces in $N$ which we denote by $S_n$ and $\hat{S}_n$ respectively. By Proposition \procref{PropTransversalityOfInteriorHypersurfaces}, for $n\geqslant n_0$, $S_n$ lies below $\hat{S}_n$ and both $S_n$ and $\hat{S}_n$ are transverse $\Cal{F}$. Upon taking limits, it follows that $\Sigma_\infty$ and $\hat{\Sigma}_\infty$ identify with smooth compact LSC embedded hypersurfaces in $N$ which we denote by $S_\infty$ and $\hat{S}_\infty$ respectively. We claim that $S_\infty$ is transverse to $\Cal{F}$. Indeed, otherwise, taking limits, there exists a geodesic segment normal to $S_\infty'$ which is an interior tangent to $S_\infty$ at some point. This is absurd, and it follows that $S_\infty$ is transverse to $\Cal{F}$ as asserted. Likewise, $\hat{S}_\infty$ is also transverse to $\Cal{F}$. Furthermore, $S_\infty$ lies below $\hat{S}_\infty$, and since $K(S_\infty)<K(\hat{S}_\infty)$, it follows from the geometric maximum principle that $S_\infty$ is not tangent to $\hat{S}_\infty$ at any point. By Proposition \procref{PropExcisionI}, $\Sigma_\infty'<\hat{\Sigma}_\infty$, and by Proposition \procref{PropExcisionII}, $\Sigma_\infty<\hat{\Sigma}_\infty$, as desired.\qed
\newsubhead{Boundedness and minimal embedding radii}In Proposition $4.1.1$ of \cite{SmiPPH}, we obtain a-priori lower bounds of the radius about any boundary point over which a smooth compact LSC immersed hypersurface is embedded. We introduce the terminology required to use this result within the current framework. {\sl The definitions and results of this section, though highly technical, are of central importance to the sequel.}\/ We recommend the reader study them carefully.
\subheadlabel{SubheadBoundednessAndMinimalEmbeddingRadii}
\medskip
\noindent Let $M:=M^{n+1}$ be an $(n+1)$-dimensional Hadamard manifold. Let $\Sigma:=(i,(S,\partial S))$ be a smooth immersed submanifold in $M$. Let $p$ be a point in $S$ which we identify with its image in $M$. Let $U$ be a neighbourhood of $p$ in $M$. Let $S'$ be the connected component of $i^{-1}(U)$ containing $p$. We denote $\Sigma\minter_p U:=(i,(S',\partial S'))$. We refer to $\Sigma\minter_p U$ as the {\bf connected component} of $\Sigma\minter U$ containing $p$.
\medskip
\noindent Let $\Gamma:=(i,G)$ be a smooth oriented compact codimension-$2$ submanifold of $M$. Let $\opN\Gamma$ be the circle bundle of unit normal vectors over $\Gamma$. For $N\in\opN\Gamma$, we define $A_\Gamma(N)$ to be the {\bf second fundamental form} of $\Gamma$ in the direction of $N$. In other words, if $X$ and $Y$ are vector fields tangent to $\Gamma$, then:
$$
A_\Gamma(N)(X,Y) = -\langle\nabla_X Y,N\rangle.\eqnum{\nexteqnno}
$$
\noindent For all $p\in\Gamma$, we define $\opCN_p\Gamma$ by:
$$
\opCN_p\Gamma = \left\{ N\in\opN_p\Gamma\ |\ A_\Gamma(N) > 0\right\},
$$
\noindent where, for any matrix $A$, we write $A>0$ whenever $A$ is positive definite. We refer to $\opCN\Gamma$ as the bundle of {\bf convex normal vectors} over $\Gamma$ at $p$. We say that $\Gamma$ is {\bf locally strictly convex} (LSC) whenever $\opCN\Gamma$ has non-trivial fibre above every point. In this case, by Proposition $4.1.2$ of \cite{SmiPPH}, for all $p\in\Gamma$, $\opCN_p\Gamma$ is an open subinterval of $\opN_p\Gamma$ of length at most $\pi$. Furthermore, by Proposition $4.1.3$ of \cite{SmiPPH}, $\opCN_p\Gamma$ varies continuously with $p$ in the Hausdorff topology.
\medskip
\noindent Let $\Cal{G}$ be a family of smooth oriented compact codimension-$2$ LSC immersed submanifolds in $M$. Suppose that $\Cal{G}$ is compact in the $C^\infty$-sense. For all $\Gamma=(i,G)$ in $\Cal{G}$, for all points $p$ in $G$, and for all $N\in\opCN_p\Gamma$, there exists a (non-complete) smooth LSC embedded hypersurface $H$ in $M$ such that $p$ lies in $H$; $N$ is the upward-pointing normal to $H$ at $p$; and $\Gamma\minter_p B_r(p)$ is contained in $H\minter_p B_r(p)$ for some small $r$. Although $H$ is, strictly speaking, not canonical, we assume throughout the sequel that $H$ depends continuously on $\Gamma\in\Cal{G}$ and $N\in\opCN\Gamma$. We then refer to $H$ as the {\bf LSC extension} of $\Gamma$ at $p$ with normal $N$. We denote $\Gamma\minter_p B_r(p)$ and $H\minter_p B_r(p)$ by $\Gamma_{p,r}$ and $H_r(\Gamma,N)$ respectively. Observe that, upon reducing $r$ if necessary, $\Gamma_{p,r}$ divides $H_r(\Gamma,N)$ into two connected components. We denote the closure of the component lying to the left (resp. right) of $\Gamma_{p,r}$ by $H^-_r(\Gamma,N)$ (resp. $H^+_r(\Gamma,N)$).
\medskip
\noindent A subset $X$ of $M$ is said to be {\bf convex} whenever the shortest geodesic in $M$ joining any two points of $X$ is also contained in $X$. Let $\Cal{G}$ be as above. We make a uniform choice of the radius $r$ as follows. For $\theta>0$, there exists $r>0$ with the following property: for all $\Gamma=(i,G)$ in $\Cal{F}$; for all points $p$ in $G$; and for all $N\in\opCN_p\Gamma$ such that $N$ makes an angle of at least $\theta$ with each end-point of $\opCN_p\Gamma$, if $H:=H_r(\Gamma,N)$ is the LSC extension of $\Gamma$ at $p$ with normal $N$, then $H$ is closed in $B_r(p)$; $H$ meets $\partial B_r(p)$ transversally; and $H$ divides $B_r(p)$ into two connected components, one of which is convex. We denote the closure of the convex component by $X_r(\Gamma,N)$.
\medskip
\noindent This construction yields a canonical order on $\opCN_p\Gamma$. Indeed, choose $\Gamma$ in $\Cal{G}$. Let $p$ be a point in $\Gamma$. Let $N$ and $N'$ be distinct vectors in $\opCN_p\Gamma$. Let $\theta>0$ be such that both $N$ and $N'$ make an angle of at least $\theta$ with each end-point of $\opCN_p\Gamma$. Let $r$ be defined as above. Since $H_r(\Gamma,N)$ and $H_r(\Gamma,N')$ meet transversally at $p$, upon reducing $r$ if necessary, we may assume that $H_r^+(\Gamma,N')$ is wholly contained in one of the connected components of the complement of $H_r(\Gamma,N)$ in $B_r(p)$. We say that $N$ lies {\bf above} $N'$ whenever $H_r^+(\Gamma,N')$ is wholly contained in the convex component $X_r(\Gamma,N)$. We define $\partial^+\opCN_p\Gamma$ (resp. $\partial^-\opCN_p\Gamma$) to be the upper (resp. lower) end-point of $\opCN_p\Gamma$.
\medskip
\noindent Proposition $4.1.1$ of \cite{SmiPPH} yields a priori lower bounds of the radius over which a smooth compact LSC hypersurface with generic boundary is embedded. In the present context, this is expressed as follows:
\proclaim{Theorem \nextprocno}
\noindent Let $\Cal{G}$ be a family of smooth oriented compact generic codimension-$2$ LSC immersed submanifolds in $M$. Suppose that $\Cal{G}$ is compact in the $C^\infty$-sense. For all $\theta>0$, there exists $r>s>0$ with the following property. If $\Sigma=(i,(S,\partial S))$ is a smooth compact LSC immersed hypersurface with boundary $\Gamma:=\partial\Sigma$ in $\Cal{G}$; if $p$ is a boundary point of $\Sigma$; if $N$ is the upward-pointing unit normal to $\Sigma$ at $p$; if $N'\in\opCN_p\Gamma$ is a convex normal vector to $\Gamma$ lying above $N$; and if the angles between $\partial^+\opCN_p\Gamma$ and $N'$ and $N'$ and $N$ are no less than $\theta$, then $\Sigma\minter_p B_r(p)$ is embedded; meets $\partial B_r(p)$ transversally; is contained within $X_r(\Gamma,N)$; and divides $X_r(\Gamma,N)$ into two connected components, one of which is convex and contains a ball of radius $s$. Furthermore $\Gamma_{p,r}$ is the only component of $\Gamma$ in $\Sigma\minter_p B_r(p)$.
\endproclaim
\proclabel{ThmMinimalEmbeddingRadius}
\noindent We denote $\Sigma_{p,r}:=\Sigma\minter_p B_{r,p}$. We denote the closure of the convex component of the complement of $\Sigma_{p,r}$ in $X_r(\Gamma,N)$ by $X_r(\Sigma,N)$. Observe that the boundary of $X_r(\Sigma,N)$ is the union of $H_r^-(\Gamma,N)$, $\Sigma_{p,r}$ and some open subset of $\partial B_r(p)$.
\medskip
\noindent Let $\Cal{F}$ be a family of smooth compact LSC immersed hypersurfaces with generic boundaries in $M$. Suppose that $\Cal{F}$ is compact in the $C^\infty$ topology. Let $\Cal{G}$ be the family of all $\partial\hat{\Sigma}$ where $\hat{\Sigma}\in\Cal{F}$. Choose $\hat{\Sigma}\in\Cal{F}$. Denote $\Gamma=\partial\hat{\Sigma}$. Let $\hat{N}$ be the upward-pointing unit normal vector field over $\hat{\Sigma}$. Let $p$ be a boundary point of $\hat{\Sigma}$. Observe that $\hat{N}(p)\in\opCN_p\partial\Gamma$. We denote by $\opCN_p^+\hat{\Sigma}$ (resp. $\opCN_p^-\hat{\Sigma}$) the open subinterval of $\opCN_p\Gamma$ lying above (resp. below) $\hat{N}(p)$. We extend $\hat{\Sigma}$ smoothly across its boundary to obtain a larger smooth compact LSC immersed hypersurface $\hat{\Sigma}'$. Although $\hat{\Sigma}'$ is, strictly speaking, not canonical, we assume throughout the sequel that it depends continuously on $\hat{\Sigma}\in\Cal{F}$. Furthermore, we may suppose that $H_r(\Gamma,\hat{N}(p))=\hat{\Sigma}'\minter_p B_r(p)$, and we denote $H_r(\Gamma,\hat{N})$ and $X_r(\Gamma,\hat{N})$ by $\hat{\Sigma}_{p,r}$ and $X_{p,r}(\hat{\Sigma})$ respectively. Likewise, if $\Sigma$ is another smooth compact LSC immersed hypersurface in $M$ such that $\Sigma<\hat{\Sigma}$ and $\Sigma$ makes an angle of at least $\theta$ with $\hat{\Sigma}$ along their common boundary, for suitable values of $r>0$, we denote $X_{p,r}(\Sigma,\hat{\Sigma})$ instead of $X_r(\Sigma,\hat{N}(p))$.
\newsubhead{Supporting normal vectors} We recall the definition of supporting normal vectors to convex sets. Let $X$ be any subset of $M$. Let $p$ be a point in $X$ and let $N\in T_pM$ be a unit vector to $M$ at $p$. Let $P\subseteq T_pM$ be the hyperplane normal to $N$. We identify $P$ with its image under the exponential map in $M$. We orient $P$ such that $N$ points upwards. Observe that $P$ divides $M$ into two connected components. We say that $N$ is a {\bf supporting normal} to $X$ at $p$ whenever $X$ lies in the closure of the connected component lying below $P$. We recall that the set of supporting normals to $X$ at $p$ is a closed convex subset of the sphere of unit vectors in $T_pM$ (c.f. \cite{SmiPPG}). Furthermore, when the set $X$ is convex, the set of supporting normals to $X$ is non-trivial at all of its boundary points (c.f. \cite{SmiPPG}).
\subheadlabel{SubheadSupportingNormalVectors}
\medskip
\noindent We recall that the set of supporting normal vectors to a convex set varies upper-semi\-con\-tinously in the Hausdorff sense. This permits us to obtain uniform moduli of continuity as follows. Let $(\hat{\Sigma}_n)_\ninn$, $\hat{\Sigma}$ be smooth compact LSC immersed hypersurfaces with generic boundaries in $M$ such that $(\hat{\Sigma}_n)_\ninn$ converges to $\hat{\Sigma}$ in the $C^\infty$ topology. For all $n$, let $p_n$ be a boundary point of $\hat{\Sigma}_n$. Suppose that $(p_n)_\ninn$ converges to the boundary point $p$ of $\hat{\Sigma}$. Choose $\theta>0$ and for all $n$, let $\Sigma_n$ be a smooth compact LSC immersed hypersurface in $M$ such that $\Sigma_n<\hat{\Sigma}_n$ and $\Sigma_n$ makes an angle of at least $\theta$ with $\hat{\Sigma}_n$ along their common boundary. For all $n$, let $N_n$ be the upward-pointing unit normal vector field over $\Sigma_n$. Let $r>s>0$ be as in Theorem \procref{ThmMinimalEmbeddingRadius}. For all $n$, denote $X_n:=X_{p_n,r}(\Sigma_n,\hat{\Sigma}_n)$. For all $n$, by definition, $X_n$ is contained in $B_r(p_n)$. Furthermore, by Theorem \procref{ThmMinimalEmbeddingRadius}, for all $n$, $X_n$ contains a ball of radius $s$. Thus, by compactness of the family of convex sets, we may assume that there exists a convex set $X$ with non-trivial interior towards which $(X_n)_\ninn$ converges in the Hausdorff sense. Observe that $p$ is a boundary point of $X$. Let $I$ be the set of supporting normals to $X$ at $p$:
\proclaim{Proposition \nextprocno}
\noindent $I$ is contained in the closure of $\opCN_p^-\hat{\Sigma}$.
\endproclaim
\proclabel{PropSupportingNormalIsContainedBelowOuterBarrier}
\proof For all $n$, let $\Gamma_n$ be the boundary of $\hat{\Sigma}_n$. Let $\Gamma$ be the boundary of $\hat{\Sigma}$. For all $n$, $\Gamma_{n,p_n,r}$ is contained in $X_n$. Thus, taking limits, $\Gamma_{p,r}$ is contained in $X$. It follows that $I$ is contained in the closure of $\opCN_p\Gamma$. For all $n$, $\hat{\Sigma}_{n,p_n,r}^-$ is contained in $X_n$. Thus, taking limits, $\hat{\Sigma}_{p,r}^-$ is contained in $X$. It follows that $I$ is contained in the closure of $\opCN_p^-\hat{\Sigma}$, as desired.\qed
\goodbreak
\proclaim{Proposition \nextprocno}
\noindent Upon extracting a subsequence, there exist subsets $(I_n)_\ninn$ of $UM$ and a continuous function $m:[0,\infty[\rightarrow[0,\infty[$ such that:
\medskip
\myitem{$(1)$} for all $n$, $I_n$ is contained in the closure of $\opCN_{p_n}^-\hat{\Sigma}_n$;
\medskip
\myitem{$(2)$} $(I_n)_\ninn$ converges to $I$ in the Hausdorff topology;
\medskip
\myitem{$(3)$} $m(0)=0$; and
\medskip
\myitem{$(4)$} for all $n$, and for all $q\in\Sigma_{n,p_n,r}$:
$$
D(I_n,N_n(q))\leqslant m(d(p,q)),
$$
\noindent where $D$ and $d$ are the distances in the total space of $UM$ and in $M$ respectively.
\endproclaim
\proclabel{PropUniformModulusOfContinuity}
\remark Observe that, if $X,Y\in UM$ are two unit vectors in the same fibre, then $D(X,Y)$ is the angle between $X$ and $Y$.\qed
\medskip
\proof Let $l$ be the length of $I$. For all $n$, we define $I_n$ to be the shortest interval of the closure of $\opCN_{p_n}^-\hat{\Sigma}_n$ which contains both $\hat{N}_n(p_n)$ and $N_n(p_n)$ and has length at least $l$. $(I_n)_\ninn$ trivially satisfies Condition $(1)$. By upper-semicontinuity of the sets of supporting normals to convex sets (c.f. \cite{SmiPPG}), upon extracting a subsequence, we may suppose that $(N_n(p_n))_\ninn$ converges to a vector in $I$. Condition $(2)$ follows. Let $(n_k,q_k)_{k\in\Bbb{N}}$ be such that for all $k$, $n_k\in\Bbb{N}$ and $q_k\in\Sigma_{n_k,p_{n_k},r}$. Denote $p_\infty:=p$ and $I_\infty:=I$. Suppose that $(n_k,q_k)_\ninn$ converges to $(n_\infty,p_{n_\infty})$ where $n_\infty\in\Bbb{N}\munion\left\{\infty\right\}$. By upper-semincontinuity of the sets of supporting normals to convex sets, every limit point of the sequence $(N_{n_k}(q_k))_{k\in\Bbb{N}}$ is contained in $I_{n_\infty}$. The result follows.\qed
\newhead{First Order Lower Estimates}
\headlabel{HeadFirstOrderLowerEstimates}
\newsubhead{Main results} Let $M:=M^{n+1}$ be an $(n+1)$-dimensional Riemannian manifold. Let $K$ be a convex curvature function. Let $\Cal{F}$ be a family of pairs $(\hat{\Sigma},\kappa)$ where $\hat{\Sigma}$ is a smooth compact LSC immersed hypersurface with generic boundary in $M$; $\kappa$ is a smooth positive function over $M$; and $K(\hat{\Sigma})>\kappa$. We furnish $\Cal{F}$ with the product topology of the $C^\infty$ topology in the first component and the $C^\infty_\oploc$ topology in the second. We obtain lower estimates for the normals of LSC hypersurfaces of prescribed $K$-curvature.
\subheadlabel{SubheadFirstOrderLowerEstimatesMainResults}
\proclaim{Proposition \nextprocno}
\noindent If $\Cal{F}$ is compact, then there exists $\theta>0$ such that for all smooth compact LSC immersed hypersurfaces $\Sigma$ in $M$ such that $K(\Sigma)=\kappa$ and $\Sigma<\hat{\Sigma}$ for some element $(\hat{\Sigma},\kappa)\in\Cal{F}$, if $N$ is the upward pointing unit normal vector field of $\Sigma$, then $N(p)$ makes an angle of at least $\theta$ with each of the end-points of $\opCN_p\partial\Sigma$.
\endproclaim
\proclabel{PropFirstOrderEstimatesUnboundedCase}
\noindent When $K$ is of bounded type, we require the following finer estimate:
\proclaim{Proposition \nextprocno}
\noindent Suppose that $K$ is of bounded type. If $\Cal{F}$ is compact, then there exists $\delta>0$ such that for all smooth compact LSC immersed hypersurfaces $\Sigma$ in $M$ such that $K(\Sigma)=\kappa$ and $\Sigma<\hat{\Sigma}$ for some element $(\hat{\Sigma},\kappa)\in\Cal{F}$, and for all $p\in\partial\Sigma$:
$$
K_\infty(A_\Gamma(N(p)) \geqslant \kappa(p) + \delta,
$$
\noindent where $N$ is the outward pointing unit normal vector field over $\Sigma$ and $A_\Gamma$ is the second fundamental form of $\Gamma := \partial\Sigma$, as defined in Section \subheadref{SubheadBoundednessAndMinimalEmbeddingRadii}.
\endproclaim
\proclabel{PropFirstOrderLowerEstimatesBoundedCase}
\noindent Logically, Proposition \procref{PropFirstOrderEstimatesUnboundedCase} precedes Proposition \procref{PropFirstOrderLowerEstimatesBoundedCase}. However, the proofs are similar, and we prove Proposition \procref{PropFirstOrderLowerEstimatesBoundedCase} first as it is more involved. We then review the modifications required to prove Proposition \procref{PropFirstOrderEstimatesUnboundedCase}.
\medskip
\noindent The following technical result will be used repeatedly throughout the sequel:
\proclaim{Lemma \nextprocno}
\noindent Let $f:M\rightarrow\Bbb{R}$ be a smooth function. Let $\opHess(f)$ and $\opHess^\Sigma(f)$ denote the Hessians of $f$ and the restriction of $f$ to $\Sigma$ respectively. Then:
$$
\opHess^\Sigma(f) = \opHess(f)|_{T\Sigma} - \langle\nabla f,N\rangle II,
$$
\noindent where $N$ is the upward-pointing unit normal vector field over $\Sigma$ and $II$ is the corresponding second fundamental form.
\endproclaim
\proclabel{PropHessianOfRestriction}
\proof Choose $p\in\Sigma$. Let $X$ and $Y$ be tangent vector fields over $\Sigma$ which are parallel at $p$. Let $\nabla^M$ denote the Levi-Civita covariant derivative of $M$. Since $Y$ is parallel at $p$ as a vector field over $\Sigma$, $\nabla^M_X Y = -II(X,Y)N$. Thus:
$$\matrix
\opHess(f)(X,Y)(p) \hfill&= (X(Yf))(p) - df(\nabla^M_X Y)(p)\hfill\cr
&= \opHess^\Sigma(f)(X,Y)(p) + \langle \nabla f,N\rangle II(X,Y)(p).\hfill\cr
\endmatrix$$
\noindent The result follows.\qed
\newsubhead{Analytic properties of the barrier function} In this and the following section, we assume that $K$ is of bounded type. Choose $(\hat{\Sigma},\kappa)\in\Cal{F}$. Denote $\Gamma:=\partial\hat{\Sigma}$. Let $p$ be a point in $\Gamma$. Let $d_p$ be the distance to $p$ in $M$. Let $\hat{N}$ be the upward pointing unit normal vector to $\hat{\Sigma}$ at $p$. Let $\opCN_p\Gamma$ and $\opCN_p^\pm\hat{\Sigma}$ be defined as in Section \subheadref{SubheadBoundednessAndMinimalEmbeddingRadii}. Choose $N\in\opCN_p^-\hat{\Sigma}$ such that $K_\infty(A_\Gamma(N))=\kappa(p)$. We aim to show that it is not possible for a sequence of smooth compact LSC hypersurfaces satisfying the hypotheses of Proposition \procref{PropFirstOrderLowerEstimatesBoundedCase} to have as a limit a hypersurface whose normal at $p$ is equal to $N$.
\subheadlabel{SubheadAnalyticPropertiesOfTheBarrierFunction}
\medskip
\noindent Let $\Phi_0$ be a smooth function defined in a neighbourhood of $p$ such that $\nabla\Phi_0(p)=N$ and $\Phi_0$ vanishes along $\Gamma$. The required barrier function is constructed by perturbing $\Phi_0$ as follows. Choose $V\in\opCN_p^+\hat{\Sigma}$. We denote $H_0:=H_{p,r}(\Gamma,V)$ for sufficiently small $r$. For $V$ sufficiently close to $\hat{N}$, we may suppose that $K(H_0)(p)>\kappa(p)$. Upon perturbing $H_0$ slightly and reversing the orientation, we obtain a smooth locally strictly {\sl concave} embedded hypersurface $H$ in $M$ such that $H$ passes through $p$; the upward pointing normal to $H$ at $p$ is equal to $-V$; $K(-H)(p)>\kappa(p)$, where $-H$ is $H$ furnished with the reverse orientation; and $H_0$ lies above the graph of $\epsilon d_p^2$ over $H$ for some $\epsilon>0$. Observe in particular that $\Gamma_{p,r}$ also lies above the graph of $\epsilon d_p^2$ over $H$. Furthermore, since $V$ lies above $\hat{N}$, $\hat{\Sigma}_{p,r}^+$ also lies above the graph of $\epsilon d_p^2$ over $H$. We define $d_H$ near $p$ to be the signed distance to $H$ in $M$ with sign chosen so that it is positive above $H$. For appropriate functions $x$ and $h$ defined near $p$, we define $\Phi_1$ by:
$$
\Phi_1 = \Phi_0 + x(d_H - h).\eqnum{\nexteqnno}
$$
\noindent Given $h$, the function $x$ is determined as follows. For any two functions $f$ and $g$ with non-colinear derivatives at $p$, we define the $(n-2)$-dimensional distribution $E(f,g)$ near $p$ by:
\eqnlabel{EqnDefinitionOfPhiOne}
$$
E(f,g) = \langle\nabla f,\nabla g\rangle^\perp,
$$
\noindent where $\langle X,Y\rangle$ here denotes the space spanned by the vectors $X$ and $Y$. Let $e_1,...,e_{n-1}$ be an orthonormal basis for $T_p\Gamma$ with respect to which $A_\Gamma(N)$ is diagonal. Observe that $\nabla d_H=-V$ and $\nabla\Phi_0=N$ are non-colinear at $p$. We thus extend $e_1,...,e_{n-1}$ to a local frame in $TM$ such that, at $p$, for all $X$ and for all $1\leqslant\alpha\leqslant n-1$:
$$\matrix
\langle\nabla_Xe_\alpha,\nabla d_H\rangle \hfill&= -\opHess(d_H)(e_\alpha,X),\hfill\cr
\langle\nabla_Xe_\alpha,\nabla\Phi_0\rangle \hfill&= -\opHess(\Phi_0)(e_\alpha,X).\hfill\cr
\endmatrix\eqnum{\nexteqnno}$$
\noindent We define the distribution $E$ near $p$ to be the linear span of $e_1,...,e_{n-1}$. Bearing in mind Lemma \procref{PropHessianOfRestriction}, for a smooth function $f$, we define $K_{\infty,E}(f)$ by:
\eqnlabel{EqnFirstDerivativeOfTheseVectorFields}
$$
K_{\infty,E}(f) = \frac{1}{\|\nabla f\|}K_\infty(\opHess(f)|_E),
$$
\noindent where $\opHess(f)|_E$ is the restriction to $E$ of the Hessian of $f$. Given $h$, the function $x$ is now determined by the following result:
\proclaim{Proposition \nextprocno}
\noindent If $N$ is not an end-point of $\opCN_p\Gamma$, then, for $V$ sufficiently close to $\hat{N}$ and for all $h$ defined near $p$ such that:
\medskip
\myitem{$(1)$} $h(p)=0$;
\medskip
\myitem{$(2)$} $\nabla h(p)=0$;
\medskip
\myitem{$(3)$} the Hessian of the restriction of $d_H - h$ to $\Gamma$ vanishes at $p$; and
\medskip
\myitem{$(4)$} the restriction of $\opHess(h)$ to $H$ is positive definite,
\medskip
\noindent there exists a function $x$ such that $x(p)=0$, $\opHess(x)(p)=0$ and:
$$
K_{\infty,E}(\Phi_0 + x(d_H - h))\leqslant \kappa + O(d_p^2).
$$
\endproclaim
\proclabel{PropFirstBarrierEstimate}
\proof By Lemma \procref{PropHessianOfRestriction}, the restriction of $\|\nabla\Phi_0\|^{-1}\opHess(\Phi_0)$ to $E_p$ coincides with $A_\Gamma(N)$. Thus, by definition of $N$, at $p$, $K_{\infty,E}(\Phi_0) = \kappa(p)$. The gradient of $x(d_H-h)$ vanishes at $p$. The Hessian of $xh$ vanishes at $p$. The Hessian of $xd_H$ vanishes on $(\nabla d_H)^\perp$ at $p$ and thus so too does its restriction to $E$. It follows that for all $x$, at $p$, $K_{\infty,E}(\Phi_1) = \kappa(p)$.
\medskip
\noindent Let $\lambda_1,...,\lambda_{n-1}$ be the eigenvalues of $A_\Gamma(N)$. Since $N$ is not an end-point of $\opCN_p\Gamma$, $\lambda_i>0$ for all $i$. Thus, by concavity, $K_\infty$ has a finite supporting tangent $(\mu_1,...,\mu_{n-1})$ at $(\lambda_1,...,\lambda_{n-1})$. Suppose first that all the $\lambda_i$ are distinct. Define $\tilde{K}_\infty$ such that, for all $\lambda_1',...,\lambda_{n-1}'$:
$$
\tilde{K}_\infty(\lambda_1',...,\lambda_{n-1}') := K_\infty(\lambda_1,..,\lambda_{n-1}) + \sum_{i=1}^n\mu_i(\lambda_i'-\lambda_i).
$$
\noindent By concavity, for all $\lambda_1',...,\lambda_{n-1}'$:
$$
K_\infty(\lambda'_1,...,\lambda'_{n-1})\leqslant\tilde{K}_\infty(\lambda'_1,...,\lambda'_{n-1}).
$$
\noindent For a smooth function $f$, we define:
$$
\tilde{K}_{\infty,E}(f) = \frac{1}{\|\nabla f\|}\tilde{K}_\infty(\opHess(f)|_E).
$$
\noindent Denote $P=x(d_{H}-h)$. At $p$:
$$
\opHess(P) = \nabla x\otimes\nabla d_H + \nabla d_H\otimes\nabla x.
$$
\noindent At $p$, for all $1\leqslant\alpha\leqslant n-1$, by definition, $\langle e_\alpha,\nabla d_H\rangle=0$. Thus, by \eqnref{EqnFirstDerivativeOfTheseVectorFields}, for all $X$, and for all $1\leqslant\alpha,\beta\leqslant n-1$:
$$\matrix
X\opHess(P)(e_\alpha,e_\beta) \hfill&= (\nabla_X\opHess(P))(e_\alpha,e_\beta) + \opHess(P)(\nabla_X e_\alpha,e_\beta) + \opHess(P)(e_\alpha,\nabla_X e_\beta)\hfill\cr
&=(\nabla_X\opHess(P))(e_\alpha,e_\beta) - \opHess(d_H)(X,e_\alpha)x_{;\beta} - \opHess(d_H)(X,e_\beta)x_{;\alpha}.\hfill\cr
\endmatrix$$
\noindent We extend $e_1,..,e_{n-1}$ to a basis $e_0,...,e_n$ for $T_pM$. Observe, in particular, that the plane spanned by $e_0$ and $e_n$ coincides with the plane spanned by $V$ and $N$. Since all the $\lambda_i$ are distinct, they are smooth in a neighbourhood of $p$ and thus, with respect to this basis, for all $1\leqslant\alpha\leqslant n-1$ and $0\leqslant k\leqslant n$, bearing in mind those terms which vanish at $p$, we obtain:
$$\matrix
\partial_k\lambda_\alpha\hfill&=\partial_k\opHess(\Phi_0)(e_\alpha,e_\alpha) + \partial_k\opHess(P)(e_\alpha,e_\alpha)\hfill\cr
&= \partial_k\opHess(\Phi_0)(e_\alpha,e_\alpha) -2x_{;\alpha}h_{;\alpha k} - x_{;k}h_{;\alpha\alpha} + x_{;k}d_{H;\alpha\alpha}.\hfill\cr
\endmatrix$$
\noindent The Hessian of the restriction of $d_H - h$ to $\Gamma$ vanishes at $P$. Furthermore $\nabla(d_H - h)(p)=\nabla d_H(p)=-V$. Thus, by Lemma \procref{PropHessianOfRestriction}, at $p$, for all $\alpha$:
$$
d_{H;\alpha\alpha} - h_{;\alpha\alpha} = A_{\Gamma}(\nabla d_H)_{\alpha\alpha} = -A_{\Gamma}(V)_{\alpha\alpha}.
$$
\noindent Thus, for all $1\leqslant\alpha\leqslant n-1$ and for all $0\leqslant k\leqslant n$ at $p$:
$$
\partial_k\lambda_\alpha = \partial_k\opHess(\Phi_0)(e_\alpha,e_\alpha) - 2x_{;\alpha}h_{;\alpha k} - x_{;k}A_{\Gamma}(V)_{\alpha\alpha}.\eqnum{\nexteqnno}
$$
\noindent Thus, for all $k$, at $p$ bearing in mind that $\tilde{K}_{\infty,E}(\Phi_1)=\kappa$ at $p$:
\eqnlabel{EqnFirstVariationOfPrincipleCurvature}
$$
\partial_k\tilde{K}_{\infty,E}(\Phi_1)=\partial_k\tilde{K}_{\infty,E}(\Phi_0) + (M\nabla x)_k,
$$
\noindent where, for any vector $U$:
$$
(MU)_k = \kappa\langle N,V\rangle U_k + \kappa\langle N,U\rangle V_k
-2\sum_{\alpha=1}^{n-1}\mu_\alpha U_{\alpha}h_{;\alpha k} - \sum_{\alpha=1}^{n-1}\mu_\alpha A_{\Gamma}(V)_{\alpha\alpha}U_k.
$$
\noindent We claim that for $V$ sufficiently close to $\hat{N}$, $M$ is invertible. Indeed, suppose that $MU=0$ for some non-trivial $U$. Taking the inner product with $(0,\mu_1U_1,...,\mu_{n-1}U_{n-1},0)$ yields:
$$2\sum_{\alpha,\beta=1}^{n-1}(\mu_\alpha U_\alpha)(\mu_\beta U_\beta)h_{;\alpha\beta}+ (\sum_{\alpha=1}^{n-1}\mu_\alpha A_{\Gamma}(V)_{\alpha\alpha}-\kappa\langle N,V\rangle)\sum_{\beta=1}^{n-1}\mu_\beta U_\beta^2=0
$$
\noindent However, by construction, and bearing in mind Proposition \procref{PropCharacterisationOfFInfinity}:
$$
K_\infty(A_\Gamma(V)) > K(-H)(p) > \kappa(p).
$$
\noindent Thus, bearing in mind concavity and Proposition \procref{PropCurvFnsII}, $(4)$ applied to $K_\infty$:
$$
\sum_{\alpha=1}^{n-1}\mu_\alpha A_\Gamma(V)_{\alpha\alpha}=K_\infty(A_\Gamma(N)) + \sum_{\alpha=1}^{n-1}\mu_\alpha(A_\Gamma(V)_{\alpha\alpha} - \lambda_\alpha)\geqslant K_\infty(A_\Gamma(V))>\kappa(p).
$$
\noindent Hence:
$$
\sum_{\alpha=1}^{n-1}\mu_\alpha A_\Gamma(V)_{\alpha\alpha}-\kappa\langle N,V\rangle>0.\eqnum{\nexteqnno}
$$
\noindent In particular:
\eqnlabel{EqnAUsefulTermThatDoesNotVanish}
$$
\sum_{\alpha,\beta=1}^{n-1}(\mu_\alpha U_\alpha)(\mu_\beta U_\beta)h_{;\alpha\beta}=0.
$$
\noindent Since the restriction of $\opHess(h)$ to $H$ is positive definite, it follows that $\mu_\alpha U_\alpha=0$ for all $1\leqslant\alpha\leqslant n-1$. Taking the inner product of $MU$ with $(0,U_1,...,U_{n-1},0)$ now yields:
$$
\left(\sum_{\alpha=1}^{n-1}\mu_\alpha A_\Gamma(V)_{\alpha\alpha} - \kappa\langle N,V\rangle\right)\sum_{\beta=1}^{n-1}U_\beta^2 = 0.
$$
\noindent Thus by \eqnref{EqnAUsefulTermThatDoesNotVanish}, $U_\alpha=0$ for all $1\leqslant\alpha\leqslant n-1$. Now observe that $M$ preserves the plane generated by $e_0$ and $e_n$. Observe, furthermore, that this plane coincides with the plane generated by $N$ and $V$. With respect to the basis $(N,V)$, the matrix of the restriction of $M$ to this plane is given by:
$$
M|_{\langle N,V\rangle} = \pmatrix
\kappa\langle N,V\rangle - \lambda(V)\hfill& 0\hfill\cr
\kappa\hfill&2\kappa\langle N,V\rangle - \lambda(V)\hfill\cr
\endpmatrix,
$$
\noindent where, for all $W$:
$$
\lambda(W) = \sum_{\alpha=1}^{n-1}\mu_\alpha A_{\Gamma}(W)_{\alpha\alpha}.
$$
\noindent By linearity, for any constant, $c$, the set of all unit vectors $W$ in this plane such that $c\langle N,W\rangle - \lambda(W)=0$ is either empty, consists of two points, or is the entire circle. However, on the one hand, by \eqnref{EqnAUsefulTermThatDoesNotVanish}:
$$
\kappa\langle N,V\rangle - \lambda(V)\neq 0.
$$
\noindent On the other hand, by Proposition \procref{PropCurvFnsII}, $(4)$:
$$
\lambda(N) = \sum_{\alpha=1}^{n-1}\mu_\alpha A_\Gamma(N)_{\alpha\alpha} = K_\infty(A_\Gamma(N)) = \kappa(p).
$$
\noindent In particular, $2\kappa\langle N,N\rangle - \lambda(N)\neq0$. The restriction of $M$ to the plane generated by $e_0$ and $e_n$ is therefore non-invertible for at most $4$ distinct values of $V$. In particular, for $V$ sufficiently close to $\hat{N}$, it is invertible. It follows that, $U_0=U_n=0$, and so $U=0$. $M$ is therefore invertible, as asserted.
\medskip
\noindent There therefore exists $x$ such that, at $p$, for all $k$:
$$
(M\nabla x)_{;k} = -\partial_k\tilde{K}_{\infty,E}(\Phi_0) + \kappa_{;k}.
$$
\noindent Consequently:
$$\matrix
&\tilde{K}_{\infty,E}(\Phi) \hfill&= \kappa + O(d_p^2)\hfill\cr
\Rightarrow\hfill&K_{\infty,E}(\Phi) \hfill&\leqslant \kappa + O(d_p^2).\hfill\cr
\endmatrix$$
\noindent Finally, if $\lambda_i=\lambda_j$ for some $i\neq j$, then, by convexity, we may choose $\mu$ such that $\mu_i=\mu_j$, and we proceed as before. This completes the proof.\qed
\medskip
\noindent For $M>0$, we define $\Phi$ by:
$$
\Phi = \Phi_0 + x(d_{H} - h) + Md_{H}^2.\eqnum{\nexteqnno}
$$
\proclaim{Proposition \nextprocno}
\noindent If $D$ represents the Grassmannian distance between two $(n-2)$-dimensional subspaces then:
$$
D(E,E(\Phi,d_{H})) = O(d_p^2) + O(d_{H}).
$$
\endproclaim
\proclabel{PropDistributionsAreClose}
\proof Bearing in mind, \eqnref{EqnFirstDerivativeOfTheseVectorFields}, for all $i$, at $p$:
$$
X\langle e_i,\nabla d_{H}\rangle = \langle\nabla_X e_i,\nabla d_H\rangle + \langle e_i,\nabla_X\nabla d_H\rangle = 0.
$$
\noindent Likewise, $X\langle e_i,\nabla\Phi_0\rangle = 0$. It follows that $D(E,E(\Phi_0,d_H))=O(d_p^2)$. Moreover, since $xh$ is of order $3$, near $p$:
$$
\nabla\Phi = \nabla\Phi_0 + (x + 2Md_{H})\nabla d_{H} + O(d_p^2) + O(d_{H}).
$$
\noindent Thus:
$$
\langle \nabla\Phi,\nabla d_{H}\rangle = \langle \nabla\Phi_0 + O(d_p^2) + O(d_{H}),\nabla d_{H}\rangle.
$$
\noindent It follows that $D(E(\Phi_0,d_H),E(\Phi,d_H))=O(d_p^2) + O(d_H)$. The result now follows by the triangle inequality.\qed
\proclaim{Proposition \nextprocno}
\noindent For $d_p\leqslant 1$, and for $Md_H\leqslant 1$:
$$\matrix
1/\|\nabla\Phi\|\hfill&\leqslant 1 + O(d_p) + O(Md_H),\hfill\cr
\|\nabla\Phi_1\|/\|\nabla\Phi\| \hfill&\leqslant 1 + 2Md_H + O(Md_Hd_p) + O(M^2d_H^2).\hfill\cr
\endmatrix$$
\endproclaim
\proclabel{PropControlOfDerivativeOfPhi}
\proof Indeed:
$$
\nabla\Phi = \nabla\Phi_1 + 2Md_H\nabla d_H.
$$
\noindent Thus:
$$
\|\nabla\Phi\|^2 \geqslant \|\nabla\Phi_1\|^2 + 4Md_H\langle\nabla d_H,\nabla\Phi_1\rangle.
$$
\noindent Observe that:
$$
\|\nabla\Phi_1\| \leqslant 1 + O(d_p).\eqnum{\nexteqnno}
$$
\noindent Furthermore $\|\nabla d_H\|=1$. It follows by the Cauchy-Schwarz inequality that:
\eqnlabel{EqnControlOfPhiOne}
$$
\|\nabla\Phi\|^2 \geqslant \|\nabla\Phi_1\|^2 - 4Md_H + O(Md_Hd_p).
$$
\noindent Since $Md_H\leqslant 1$, by it follows by Taylor's Theorem that:
$$
1/\|\nabla\Phi\| \leqslant 1/\|\nabla\Phi_1\| + 2Md_H/\|\nabla\Phi_1\|^3 + O(Md_Hd_p) + O(M^2d_H^2).
$$
\noindent The first inequality follows by \eqnref{EqnControlOfPhiOne}. Furthermore, by \eqnref{EqnControlOfPhiOne} again:
$$
\|\nabla\Phi_1\|/\|\nabla\Phi\| \leqslant 1 + 2Md_H + O(Md_Hd_p) + O(M^2d_H^2),
$$
\noindent as desired.\qed
\proclaim{Proposition \nextprocno}
\noindent If $N$ is not an end-point of $\opCN_p\Gamma$, then there exists $\delta>0$ such that for $M\geqslant 1$, for $d_p\leqslant 1$ and for $Md_H\leqslant 1$:
$$
K_{\infty,E}(\Phi) \leqslant \kappa - \delta Md_H + O(d_p^2) + O(d_H) + O(Md_Hd_p) + O(M^2d_H^2).
$$
\endproclaim
\proclabel{PropFirstBarrierEstimateII}
\proof By Proposition \procref{PropFirstBarrierEstimate}:
$$
K_{\infty,E}(\Phi_1) \leqslant \kappa + O(d_p^2).
$$
\noindent Since $\opHess(\Phi_1) = O(1)$, by Proposition \procref{PropDistributionsAreClose}:
$$
K_{\infty,E(\Phi,d_H)}(\Phi_1) \leqslant \kappa + O(d_p^2) + O(d_H).
$$
\noindent Thus, by Proposition \procref{PropControlOfDerivativeOfPhi}:
$$
\frac{\|\nabla\Phi_1\|}{\|\nabla\Phi\|}K_{\infty,E(\Phi,d_H)}(\Phi_1) \leqslant (1 + 2Md_H)\kappa + O(d_p^2) + O(d_H) + O(Md_Hd_p) + O(M^2d_H^2).
$$
\noindent However, differentiating $Md_H^2$ yields:
$$
\opHess(Md_H^2) = 2M\nabla d_H\otimes\nabla d_H + 2Md_H\opHess(d_H).
$$
\noindent Observe that $\nabla d_H\otimes\nabla d_H$ vanishes along $(\nabla d_H)^\perp$. Furthermore, since $H$ is locally strictly concave, the restriction of $\opHess(d_H)$ to $E(\Phi,d_H)$ is negative definite. Thus, if the restriction of $\opHess(\Phi_1)$ to $E(\Phi,d_H)$ is not positive definite, then the restriction of $\opHess(\Phi)$ to $E(\Phi,d_H)$ is also not positive definite, and we are done. Otherwise, let $A$ and $B$ be the restrictions of $\|\Phi\|^{-1}\opHess(\Phi_1)$ and $2Md_H\|\Phi\|^{-1}\opHess(d_H)$ to $E(\Phi,d_H)$ respectively. By concavity of $K_\infty$:
$$
\|\nabla\Phi\|^{-1}K_\infty(\opHess(\Phi)|_{E(\Phi,d_H)}) = K_\infty(A + B) \leqslant K_\infty(A) + DK_{\infty,A}(B).
$$
\noindent However by construction, $K(-H)(p)>\kappa(p)$, where $-H$ is $H$ furnished with the reverse orientation. Thus, bearing in mind Proposition \procref{PropCharacterisationOfFInfinity}, Proposition \procref{PropCurvFnsII}, $(5)$ and Lemma \procref{PropHessianOfRestriction}, at $p$:
$$\matrix
DK_{\infty,A}(-\opHess(d_H)|_{E(\Phi,d_H)}) \hfill&\geqslant K_\infty(-\opHess(d_H|_{E(\Phi,d_H)}))\hfill\cr
&> K(-\opHess(d_H)|_{(\nabla d_H)^\perp})\hfill\cr
&=K(-H)(p)\hfill\cr
&>\kappa(p).\hfill\cr
\endmatrix$$
\noindent Thus, bearing in mind Proposition \procref{PropControlOfDerivativeOfPhi}, there exists $\delta>0$ such that:
$$
DK_{\infty,A}(B) \leqslant -Md_H(2\kappa + \delta) + O(Md_H d_p) + O(M^2d_H^2).
$$
\noindent Combining these relations yields:
$$
K_{\infty,E}(\Phi) \leqslant \kappa - \delta Md_H + O(d_p^2) + O(d_H) + O(Md_Hd_p) + O(M^2d_H^2),
$$
\noindent as desired.\qed
\newsubhead{Geometric properties of the barrier} We now suppose the contrary in Proposition \procref{PropFirstOrderLowerEstimatesBoundedCase}. Let $(\hat{\Sigma}_n,\kappa_n)_\ninn$ be a sequence in $\Cal{F}$ converging to $(\hat{\Sigma},\kappa)$, say. For convenience, we assume that $(\hat{\Sigma}_n,\kappa_n)=(\hat{\Sigma},\kappa)$ for all $n$. Let $\Gamma$ be the boundary of $\hat{\Sigma}$. Let $\hat{N}$ be the upward-pointing unit normal vector to $\hat{\Sigma}$ at $p$. For all $n$, let $\Sigma_n$ be a smooth compact LSC immersed hypersurface such that $K(\Sigma_n)=\kappa_n$ and $\Sigma_n<\hat{\Sigma}_n$. For all $n$, let $N_n$ be the upward-pointing unit normal vector field over $\Sigma_n$. For all $n$, let $p_n$ be a point in $\Gamma$. Suppose that $(p_n)_\ninn$ converges to the point $p$ in $\Gamma$ and that $(N_n(p_n))_\ninn$ converges to a limit $N$, say. For convenience, we assume that $p_n=p$ for all $n$. By Proposition \procref{PropFirstOrderEstimatesUnboundedCase}, we may suppose that there exists $\theta>0$ such that for all $n$, $N_n(p)$ makes an angle of at least $\theta$ with the lower end-point of $\opCN_p(\Gamma)$. In particular, $N$ is an element of $\opCN_p\Gamma$. $A_\Gamma(N)$ is therefore positive-definite.
\subheadlabel{SubheadGeometricPropertiesOfFirstOrderBarrier}
\medskip
\noindent We assume that $K_\infty(A_\Gamma(N))\leqslant\kappa(p)$ and obtain a contradiction. Let $V$, $H$, $d_H$, $\Phi_0$ and $\Phi_1$ be defined as in the preceeding section. Define the family $\Cal{G}$ by $\Cal{G}:=\left\{ \partial\hat{\Sigma}\ |\ (\hat{\Sigma},\kappa)\in\Cal{F}\right\}$. Let $\theta>0$ be such that both $\partial^+\opCN_p\Gamma$ and $V$ and $V$ and $\hat{N}$ make angles of at least $\theta$. Suppose furthermore that $\hat{N}$ and $N_n$ make angles of at least $\theta$ for all $n$. Let $r>s>0$ be as in Theorem \procref{ThmMinimalEmbeddingRadius} for the family $\Cal{G}$ and the angle $\theta$. For all $n$, we henceforth identify $\Sigma_n$, $\hat{\Sigma}$ and $\Gamma$ with $\Sigma_{n,p,r}$, $\hat{\Sigma}_{p,r}$ and $\Gamma_{p,r}$. By Theorem \procref{ThmMinimalEmbeddingRadius}, for all $n$, $\Sigma_n$ is embedded; meets $\partial B_r(p_n)$ transversally; and is contained within both $X_{p,r}(\hat{\Sigma})$ and $X_r(\Gamma,V)$. We denote $X_n:=X_{p,r}(\Sigma_n,\hat{\Sigma})$. By Theorem \procref{ThmMinimalEmbeddingRadius}, for all $n$, $X_n$ contains a ball of radius $s$. Thus, by compactness of the family of convex sets, we may suppose that there exists a convex subset $X$ with non-trivial interior towards which $(X_n)_\ninn$ converges in the Hausdorff sense. Observe that $p$ is a boundary point of $X$. Let $I$ be the set of supporting normals to $X$ at $p$. By Proposition \procref{PropSupportingNormalIsContainedBelowOuterBarrier}, $I$ is contained in the closure of $\opCN_p^-\hat{\Sigma}$.
\proclaim{Proposition \nextprocno}
\noindent If $N'\in I$ lies strictly below $N$, then $K_\infty(A_\Gamma(N'))<K_\infty(A_\Gamma(N))\leqslant\kappa(p)$.
\endproclaim
\proof Observe that $\hat{N}(p)$ lies strictly above $N$. Furthermore, by Proposition \procref{PropCharacterisationOfFInfinity}, $K_\infty(A_\Gamma(\hat{N}(p)))>\kappa(p)$. Thus, if $K_\infty(A_\Gamma(N'))\geqslant\kappa(p)$, then, by linearity of $A_\Gamma$ and concavity of $K_\infty$, $K_\infty(A_\Gamma(N))>\kappa(p)$. This is absurd, and the result follows.\qed
\medskip
\noindent Let $N_0$ be the lower end-point of $I$. We shall see presently that when $K_\infty(A_\Gamma(N_0))<\kappa(p)$ it is relatively straightforward to construct comparison hypersurfaces of $K$-curvature strictly less than $\kappa$. The case of equality is more subtle. We therefore assume that $N_0=N$ and that $K_\infty(A_\Gamma(N))=\kappa(p)$.
\medskip
\noindent Let $\Sigma$ be the closure of the intersection of $\partial X$ with the interior of $X_{p,r}(\hat{\Sigma})$. For $\rho,\epsilon>0$, we define $U_{\rho,\epsilon}$ by:
$$
U_{\rho,\epsilon} = \left\{ q\in M\ |\ d_p(q) < \epsilon,\ d_H(q)<\rho\epsilon^2\right\}.
$$
\proclaim{Proposition \nextprocno}
\noindent There exists $\rho>0$ such that for all sufficiently small $\epsilon$, and for all $p\in\Sigma\minter\partial U_{\rho,\epsilon}$, $d_p(q)<\epsilon$ and $d_H(q) = \rho\epsilon^2$.
\endproclaim
\proclabel{PropASuitableValueOfRhoRestrictsTheBoundary}
\proof Indeed, by definition of $H$, there exists $\rho>0$ such that for sufficiently small $\epsilon$ and for all $p\in X_r(\Gamma,V)\minter\partial U_{\rho,\epsilon}$, $d_p(q)<\epsilon/2$ and $d_H(q)=\rho\epsilon^2$. Since $\Sigma$ is contained in $X_r(\Gamma,V)$, the result follows.\qed
\medskip
\noindent We choose $\rho$ as in Proposition \procref{PropASuitableValueOfRhoRestrictsTheBoundary}, and henceforth keep it fixed. Since $\Gamma$ is LSC and lies strictly above $H$, we may define the function $h$ near $p$ such that $h(p)=0$; $\nabla(h)(p)=0$; the restriction of $\opHess(h)(p)$ to $T_pH$ is positive definite; and $d_H-h=O(d_p^3)$ along $\Gamma$. We choose $x$ as in Proposition \procref{PropFirstBarrierEstimate}. We recall that for $M>0$, the function $\Phi$ is given by:
$$
\Phi = \Phi_0 + x(d_H - h) + Md_H^2.
$$
\noindent Observe that the only parameters that remain to be determined at this stage are $\epsilon$ and $M$. $M$ is chosen as a function of $\epsilon$ using the following result:
\eqnlabel{EqnDefinitionOfPhi}
\proclaim{Proposition \nextprocno}
\noindent There exists a continuous function $B:[0,\infty[\rightarrow[0,\infty[$ such that $B(0)=0$ and for all sufficiently small $\epsilon$, if $M\geqslant B(\epsilon)\epsilon^{-2}$, then $\Phi\geqslant 0$ along $\partial(\Sigma\minter U_{\rho,\epsilon})$.
\endproclaim
\proclabel{PropHowToChooseM}
\proof Let $m:[0,\infty[\rightarrow[0,\infty[$ be as in Proposition \procref{PropUniformModulusOfContinuity}. Recall that, by definition, $-\nabla d_H(p)=V$ lies strictly above $\hat{N}(p)$. This vector therefore also lies strictly above every element of $I$. By Proposition \procref{PropUniformModulusOfContinuity}, there therefore exists $c>0$ such that, upon reducing $r$ if necessary, for all sufficiently large $n$ and for all $q\in\Sigma_n$:
$$
\|\pi_{n,q}(\nabla d_H)\| \geqslant c,
$$
\noindent where $\pi_{n,q}$ is the orthogonal projection of $T_pM$ onto $T_q\Sigma_n$. Furthermore, by definition $\nabla\Phi_0(p)=N_0$. Thus, upon modifying $m$ if necessary, we may suppose that, for all sufficiently large $n$, and for all $q\in\Sigma_n$:
$$
\langle\pi_{n,q}(\nabla\Phi_0),\pi_{n,q}(\nabla d_H)\rangle \geqslant -m(d_p(q)).
$$
\noindent Now consider $q\in\Sigma_n\minter\partial U_{\rho,\epsilon}$. By Proposition \procref{PropASuitableValueOfRhoRestrictsTheBoundary}, for sufficiently large $n$, $d_p(q)<\epsilon$ and $d_H(q)=\rho\epsilon^2$. Let $\gamma$ be an integral curve of $-\pi_{n,q}(\nabla d_H)$ in $\Sigma_n$ starting at $q$. Since $\Sigma_n$ is compact and since $-\pi_{n,q}(\nabla d_H)$ never vanishes, $\gamma$ meets another boundary point $q'$ of $\Sigma_n$ after finite time. Furthermore, since $d_H$ is strictly decreasing along $\gamma$, $q'$ is not an element of $\partial U_{\rho,\epsilon}$, and is therefore an element of $\Gamma$. Thus, since $d_H(q)=\rho\epsilon^2$ and $d_H=0$ along $\Gamma$, $\opLength(\gamma)\leqslant c^{-1}\rho\epsilon^2$. Since $\Phi_0$ vanishes over $\Gamma$, integration over $\gamma$ yields:
$$
\Phi_0(q) \geqslant -m(\epsilon)c^{-2}\rho\epsilon^2.
$$
\noindent Now choose $q\in\Sigma\minter\partial U_{\rho,\epsilon}$. Taking limits, we obtain:
$$
\Phi_0(q)\geqslant -m(\epsilon)c^{-2}\rho\epsilon^2.
$$
\noindent Thus, since $x(q)=O(d_p)=O(\epsilon)$, $d_H(q)=\rho\epsilon^2$ and $h(q)=O(d_p^2)=O(\epsilon^2)$, upon modifying $m$ if necessary:
$$
\Phi(q) = (\Phi_0 + x(d_H - h) + Md_H^2)(q) \geqslant Md_H^2(q) - m(\epsilon)O(\epsilon^2).\eqnum{\nexteqnno}
$$
\noindent However, $(\partial\Sigma)\minter U_{\rho,\epsilon}$ is contained in $\Gamma$. Thus, since $\Phi_0$ vanishes along $\Gamma$, by definition of $h$, along $(\partial\Sigma)\minter U_{\rho,\epsilon}$:
\eqnlabel{EqnLowerBoundsAlongTheBoundaryI}%
$$
\Phi(q) \geqslant Md_H^2 - O(d_p^4).\eqnum{\nexteqnno}%
$$
\noindent However, by definition of $H$, along $\Gamma$, $d_H\geqslant O(d_p^2)$, and along $\Sigma\minter (\partial U_{\rho,\epsilon})$, $d_H=\rho\epsilon^2$. Thus, by \eqnref{EqnLowerBoundsAlongTheBoundaryI} and \eqnlabel{EqnLowerBoundsAlongTheBoundaryII}\eqnref{EqnLowerBoundsAlongTheBoundaryII}, there exists a continuous function $B:[0,\infty[\rightarrow[0,\infty[$ such that $B(0)=0$ and, for sufficiently small $\epsilon>0$, if $M\geqslant B(\epsilon)\epsilon^{-2}$, then $\Phi\geqslant 0$ along $\partial(\Sigma\minter U_{\rho,\epsilon})$, as desired.\qed
\medskip
\noindent For all $\epsilon$ small, we henceforth choose $M:=M(\epsilon)=B(\epsilon)\epsilon^{-2}$.
\proclaim{Proposition \nextprocno}
\noindent For sufficiently small $\epsilon>0$, $\nabla\Phi$ is non-vanishing over $U_{\rho,\epsilon}$.
\endproclaim
\proclabel{PropLevelHypersurfacesAreManifolds}
\proof As in the proof of Proposition \procref{PropDistributionsAreClose}:
$$\matrix
\nabla\Phi \hfill&= \nabla\Phi_0 + (d_H - h)\nabla x + x\nabla(d_H - h) + 2 M d_H\nabla d_H\hfill\cr
&= \nabla\Phi_0 + 2M d_H\nabla d_H + O(d_p) + O(d_H).\hfill\cr
\endmatrix$$
\noindent However, by definition, $\nabla\Phi_0$ and $\nabla d_H$ are non-colinear, and the result follows.\qed
\medskip
\noindent In particular, the level sets of $\Phi$ are smooth hypersurfaces of $U_{\rho,\epsilon}$. We orient these hypersurfaces such that $\nabla\Phi/\|\nabla\Phi\|$ is the upward-pointing unit normal vector.
\proclaim{Proposition \nextprocno}
\noindent For all sufficiently small $\epsilon>0$, there exists $\delta>0$ such that if $L$ is a level subset of $\Phi$ in $U_{\rho,\epsilon}\minter X_{p,r}(\hat{\Sigma})$, if $q$ is a point of $L$, and if $A$ is the shape operator of $L$ at $q$, then, either $A$ is not positive-definite, or $K(A)\leqslant\kappa(q)-\delta$.
\endproclaim
\proclabel{PropCurvatureIsSmallerThanKappa}
\proof Let $q$ be a point in the closure of $U_{\rho,\epsilon}\minter X_{p,r}(\hat{\Sigma})$. Let $A$ be the shape operator of the level hypersurface of $\Phi$ passing through $q$. By compactness, we may suppose that $q$ maximises $K(A)$ over the closure of $U_{\rho,\epsilon}$. It suffices to show that $K(A)<\kappa(q)$ when $A$ is positive definite. However, by Lemma \procref{PropHessianOfRestriction}, $A$ is the restriction of $\|\nabla\Phi\|^{-1}\opHess(\Phi)$ to $\nabla\Phi^\perp$. Thus, if $A$ is positive definite, then bearing in mind Proposition \procref{PropCharacterisationOfFInfinity}:
$$
K(A) = K(\|\nabla\Phi\|^{-1}\opHess(\Phi)|_{(\nabla\Phi)^\perp}) < K_\infty(\|\nabla\Phi\|^{-1}\opHess(\Phi)|_{E(\Phi,d_H)}) = K_{\infty,E}(\Phi).
$$
\noindent However, $d_H\geqslant O(d_p^2)$ over $U_{\rho,\epsilon}\minter X_{p,r}(\hat{\Sigma})$. Thus, by Proposition \procref{PropFirstBarrierEstimateII}, for sufficiently small $\epsilon>0$, $K_{\infty,E}(\Phi)\leqslant\kappa$. In particular, $K(A)<\kappa(q)$, as desired.\qed
\medskip
\noindent We now complete the proof of Proposition \procref{PropFirstOrderLowerEstimatesBoundedCase}:
\medskip
{\bf\noindent Proof of Proposition \procref{PropFirstOrderLowerEstimatesBoundedCase}:\ }Recall that $N_0$ is the lower end-point of $I$. We first suppose that $K_\infty(A_\Gamma(N_0))=\kappa(p)$. Choose $\epsilon>0$ sufficiently small so that Propositions \procref{PropHowToChooseM}, \procref{PropLevelHypersurfacesAreManifolds} and \procref{PropCurvatureIsSmallerThanKappa} hold. In particular, $\Phi\geqslant 0$ along $\partial(\Sigma\minter U_{\rho,\epsilon})$ and there exists $\delta>0$ such that for all $q\in U_{\rho,\epsilon}\minter X_{p,r}(\hat{\Sigma})$, if $A$ is the shape operator of the level hypersurface of $\Phi$ passing through $q$, then either $A$ is not positive-definite, or $K(A)\leqslant\kappa(q)-\delta$. Let $X\in\opN_p\Gamma$ be the unit vector normal to $N_0$ pointing into $\Sigma$. Since $\nabla\Phi(p)=0$, we may perturb $\Phi$ to a function $\Phi'$ such $\Phi'\geqslant 0$ over $\partial(\Sigma\minter U_{\rho,\epsilon})$; for all $q\in U_{\rho,\epsilon}\minter X_{p,r}(\hat{\Sigma})$, if $A$ is the shape operator of the level hypersurface of $\Phi$ passing through $q$, then either $A$ is not positive-definite, or $K(A)<\kappa(q)$; and $\nabla\Phi'(p)=\lambda X$ for some $\lambda<0$. In particular, the restriction of $\Phi$ to $\Sigma\minter U_{\rho,\epsilon}$ achieves its minimum value at some interior point $q\in\Sigma\minter U_{\rho,\epsilon}$. If $L$ is the level hypersurface of $\Phi$ passing through $q$, then $L$ is an interior tangent to $\Sigma$ at $q$. However, since $\Sigma$ is a limit of smooth hypersurfaces of $K$-curvature prescribed by $\kappa$, the $K$-curvature of $\Sigma$ is prescribed by $\kappa$ in the viscosity sense. That is, the $K$-curvature of $L$ at $p$ is no less than $\kappa(p)$. This yields a contradiction in the case where $K_\infty(A_\Gamma(N_0))=\kappa(p)$. When $K_\infty(A_\Gamma(N_0))<\kappa(p)$, choosing $x=0$ in the above construction yields a function $\Phi$ with the required properties, and we thus also obtain a contradiction in this case. This completes the proof.\qed
\newsubhead{The unbounded case} We now consider general $K$. We outline the modifications required in the proof of Proposition \procref{PropFirstOrderLowerEstimatesBoundedCase} to obtain Proposition \procref{PropFirstOrderEstimatesUnboundedCase}. Choose $(\hat{\Sigma},\kappa)\in\Cal{F}$. Denote $\Gamma=\partial\hat{\Sigma}$. Choose $p\in\Gamma$. We define $V$ and $H$ as before. Let $N$ be the lower end-point of $\opCN_p\Gamma$. Using this $N$, we define $\Phi_0$ as before. For appropriate functions $x$ and $h$ defined near $p$, we define $\Phi_1$ as before. Define $E$ as before. For a smooth function, $f$, we define $\lambda_{1,E}(f)$ to be the least eigenvalue of the restriction of $\opHess(f)$ to $E$. Proposition \procref{PropFirstBarrierEstimate} becomes:
\subheadlabel{SubheadTheUnboundedCase}
\proclaim{Proposition \nextprocno}
\noindent For all $h$ defined near $p$ satisfying the hypotheses of Proposition \procref{PropFirstBarrierEstimate}, there exists a function $x$, defined near $p$ such that $x(p)=0$, $\opHess(x)(p)=0$, and:
$$
\lambda_{1,E}(\Phi_1) \leqslant O(d_p^2).
$$
\endproclaim
\proclabel{PropControlOfFirstEigenvalue}
\remark Importantly, in contrast to Proposition \procref{PropFirstBarrierEstimate}, it is not necessary to assume that $N$ is not an end-point of $\opCN_p\Gamma$. This is because the function $T$ used in the proof of Proposition \procref{PropControlOfFirstEigenvalue} is differentiable at the point $(0,...,0,\lambda_{m+1},...,\lambda_{n-1})$, which is not necessarily the case for general $K_\infty$.\qed
\medskip
\proof As in the proof of Proposition \procref{PropFirstBarrierEstimate}, the restrictions of $\opHess(\Phi_0)$ and $\opHess(\Phi_1)$ to $E$ are both equal to $A_\Gamma(N)$. Let $\lambda_1\leqslant\lambda_2\leqslant...\leqslant\lambda_{n-1}$ be the eigenvalues of this restriction. By definition of $N$, $\lambda_1=0$. Let $1\leqslant m\leqslant n-1$ be such that $\lambda_m=0$ but $\lambda_{m+1}>0$. For any function $f$, we define $T(f)$ to be the sum of the first $m$ eigenvalues of the restriction of $\opHess(f)$ to $E$. In particular $T(\Phi_1)(p)=T(\Phi_0)(p)=0$. By \eqnref{EqnFirstVariationOfPrincipleCurvature}, for all $k$, at $p$:
$$
\partial_k T(\Phi_1) = \partial_k T(\Phi_0) - (M\nabla x)_k,
$$
\noindent where, for any vector $U$:
$$
(M U)_k = -\sum_{\alpha=1}^m \left(2U_\alpha h_{;\alpha k} + U_kA_\Gamma(V)_{\alpha\alpha}\right).
$$
\noindent We claim that $M$ is invertible. Indeed, suppose that $MU=0$ for some non-trivial $U$. Taking the inner product with $(0,U_1,...,U_m,0,...,0)$ yields:
$$
\sum_{\alpha,\beta=1}^m U_\alpha U_\beta2 h_{;\alpha\beta} + \sum_{\alpha,\beta=1}^m U_\alpha^2 A_\Gamma(V)_{\beta\beta} = 0.
$$
\noindent Since the restriction of $\opHess(h)$ to $E$ is positive definite, and since $A_\Gamma(V)$ is positive definite, it follows that $U_1=...=U_m=0$. Thus, for all $k$:
$$
U_k = -(MU)_k/\sum_{\alpha=1}^mA_\Gamma(V)_{\alpha\alpha} = 0.
$$
\noindent It follows that the kernel of $M$ is trivial. $M$ is therefore invertible, as asserted. There therefore exists $x$ such that, at $p$, for all $k$, $\partial_k T(\Phi_1)=0$. For such a choice of $x$, $T(\Phi_1)=O(d_p^2)$. Since $\lambda_{1,E}(\Phi_1)\leqslant T(\Phi_1)/m$, the result follows.\qed
\medskip
\noindent For $M>0$, we define $\Phi$ as before. Proposition \procref{PropFirstBarrierEstimateII} now becomes:
\proclaim{Proposition \nextprocno}
\noindent There exists $\delta>0$ such that for all $M$:
$$
\lambda_{1,E(\Phi,d_H)} \leqslant -M\delta d_H + O(d_p^2) + O(d_H).
$$
\endproclaim
\proclabel{PropNonPositiveEigenvalue}
\proof By Proposition \procref{PropControlOfFirstEigenvalue}:
$$
\lambda_{1,E}(\Phi_1) \leqslant O(d_p^2).
$$
\noindent Since $\opHess(\Phi_1) = O(1)$, by Proposition \procref{PropDistributionsAreClose}:
$$
\lambda_{1,E(\Phi,d_H)}(\Phi_1) \leqslant O(d_p^2) + O(d_H).
$$
\noindent Differentiating $Md_H^2$ yields:
$$
\opHess(Md_H^2) = 2M\nabla d_H\otimes\nabla d_H + 2M d_H\opHess(d_H).
$$
\noindent The first term vanishes along $(\nabla d_H)^\perp$. On the other hand, for sufficiently small $\epsilon$, there exists $\delta_1>0$ such that for all $X\in(\nabla d_H)^\perp$, $\opHess(d_H)(X,X)\leqslant -\delta\|X\|^2/2$. Thus:
$$
\lambda_{1,E(\Phi,d_H)}(\Phi) \leqslant  - M\delta d_H + O(d_p^2) + O(d_H),
$$
\noindent as desired.\qed
\medskip
\noindent We repeat the geometric construction of Section \subheadref{SubheadGeometricPropertiesOfFirstOrderBarrier}. Observe that in this case, it is not necessary to assume that $N$ is an interior point of $\opCN_p\Gamma$ (c.f. the remark following Proposition \procref{PropControlOfFirstEigenvalue}). Proposition \procref{PropCurvatureIsSmallerThanKappa} becomes:
\proclaim{Proposition \nextprocno}
\noindent For all sufficiently small $\epsilon>0$, if $L$ is a level subset of $\Phi$ in $U_{\rho,\epsilon}\minter X_{p,r}(\hat{\Sigma})$, then $L$ has at least one non-positive principal curvature at every point.
\endproclaim
\proclabel{PropNonPositivePrincipleCurvature}
\proof Indeed, by Proposition \procref{PropNonPositiveEigenvalue}, for sufficiently small $\epsilon$, $\lambda_{1,E(\Phi,d_H)}(\Phi)\leqslant 0$ throughout $U_{\rho,\epsilon}\minter X_{p,r}(\hat{\Sigma})$. In other words, at every point $p\in U_{\rho,\epsilon}\minter X_{p,r}(\hat{\Sigma})$, the restriction of $\opHess(\Phi)$ to $E(\Phi,d_H)$ is not positive definite. In particular, for all $p\in U_{\rho,\epsilon}\minter X_{p,r}(\hat{\Sigma})$, the restriction of $\opHess(\Phi)$ to $(\nabla\Phi)^\perp$ is also not positive definite. The result now follows by Lemma \procref{PropHessianOfRestriction}.\qed
\medskip
\noindent We now prove of Proposition \procref{PropFirstOrderEstimatesUnboundedCase}:
\medskip
{\bf\noindent Proof of Proposition \procref{PropFirstOrderEstimatesUnboundedCase}:\ }In order to obtain a contradiction, we suppose that $N_0$ is the lower end-point of $\opCN_p\Gamma$. Choose $\epsilon>0$ sufficiently small so that Propositions \procref{PropHowToChooseM}, \procref{PropLevelHypersurfacesAreManifolds} and \procref{PropNonPositivePrincipleCurvature} hold. In particular, $\Phi\geqslant 0$ over $\partial(\Sigma\minter U_{\rho,\epsilon})$. By Proposition \procref{PropNonPositivePrincipleCurvature}, for all $q\in\Sigma\minter U_{\rho,\epsilon}$, the level hypersurface of $\Phi$ passing through $q$ has at least one non-positive principal curvature at $q$. Let $X\in\opN_p\Gamma$ be the unit vector normal to $N_0$ pointing into $\Sigma$. Since $\nabla\Phi(p)=0$, we may perturb $\Phi$ to a function $\Phi'$ such that $\Phi'\geqslant 0$ over $\partial(\Sigma\minter U_{p,\epsilon})$; for all $q\in\Sigma\minter U_{\rho,\epsilon}$, if $A$ is the shape operator of the level hypersurface of $\Phi$ passing through $q$, then either $A$ is not positive-definite or $K(A)<\kappa(q)$; and $\nabla\Phi'(p)=\lambda X$ for some $\lambda<0$. In particular, the restriction of $\Phi$ to $\Sigma\minter U_{\rho,\epsilon}$ achieves its minimum value at some interior point $q\in\Sigma\minter U_{\rho,\epsilon}$. If $L$ is the level hypersurface of $\Phi$ passing through $q$, then $L$ is an interior tangent to $\Sigma$ at $q$. However, since $\Sigma$ is a limit of smooth hypersurfaces of $K$-curvature prescribed by $\kappa$, the $K$-curvature of $\Sigma$ is prescribed by $\kappa$ in the viscosity sense. That is, the $K$-curvature of $L$ at $p$ is no less than $\kappa(p)$. This is absurd, and the result follows.\qed
\newhead{Second Order Boundary Estimates}%
\headlabel{HeadSecondOrderBoundaryEstimates}%
\newsubhead{Main result} Let $M:=M^{n+1}$ be an $(n+1)$-dimensional Riemannian manifold. Let $K$ be a convex curvature function. Let $\Cal{F}$ be a family of pairs $(\hat{\Sigma},\kappa)$ where $\hat{\Sigma}$ is a smooth compact LSC immersed hypersurface with generic boundary in $M$; $\kappa$ is a smooth positive function over $M$; and $K(\hat{\Sigma})>\kappa$. We furnish $\Cal{F}$ with the product topology of the $C^\infty$ topology in the first component and the $C^\infty_\oploc$ topology in the second. We obtain a-priori estimates along the boundary for the shape operators of LSC hypersurfaces of prescribed $K$-curvature.
\subheadlabel{SubheadSecondOrderBoundaryEstimatesMainResult}
\proclaim{Proposition \nextprocno}
\noindent If $\Cal{F}$ is compact, then, for $\theta>0$ small, there exists $B>0$ such that for all smooth compact LSC immersed hypersurfaces $\Sigma$ in $M$ such that $K(\Sigma)=\kappa$; $\Sigma<\hat{\Sigma}$; and $\Sigma$ makes an angle of at least $\theta$ with $\hat{\Sigma}$ along their common boundary, and for all $p\in\partial\Sigma$, if $A$ is the shape operator of $\Sigma$ at $p$, then $\|A\| \leqslant B$.
\endproclaim
\proclabel{PropSecondOrderBoundaryEstimates}
\newsubhead{Preliminary Results} The required barrier will be constructed in the following sections. Upon reducing $\theta$ if necessary, we may suppose that for all $(\hat{\Sigma},\kappa)\in\Cal{F}$ and for every boundary point $p$ of $\Sigma$, the upward-pointing unit normal vector to $\hat{\Sigma}$ at $p$ makes an angle of at least $\theta$ with the upper end-point of $\opCN_p\partial\hat{\Sigma}$. Denote $\Cal{G}=\left\{\partial\hat{\Sigma}\ |\ (\hat{\Sigma},\kappa)\in\Cal{F}\right\}$. Let $r>s>0$ be as in Theorem \procref{ThmMinimalEmbeddingRadius} for the family $\Cal{G}$ and the angle $\theta$.
\subheadlabel{SubheadSecondOrderBoundaryEstimatesPreliminaryResults}
\medskip
\noindent Choose $(\hat{\Sigma},\kappa)\in\Cal{F}$. Let $\Gamma$ be the boundary of $\hat{\Sigma}$. Let $p$ be a point of $\Gamma$. We henceforth identify $\Gamma$ and $\hat{\Sigma}$ with $\Gamma_{p,r}$ and $\hat{\Sigma}_{p,r}$ respectively. For brevity, in the sequel, we denote by $\Cal{B}$ the family of all quantities which only depend upon $\hat{\Sigma}$ and the restrictions to $B_r(p)$ of $\kappa$ and the metric of $M$. For any supplementary object, $X$, we denote by $\Cal{B}(X)$ those terms that also depend on the restriction of $X$ to $B_r(p)$, if $X$ is a submanifold or function or so on, or just on $X$ otherwise.
\medskip
\noindent Let $I$ be a closed subinterval of the closure of $\opCN_p^-\hat{\Sigma}$. Let $m:[0,\infty[\rightarrow[0,\infty[$ be a continuous function such that $m(0)=0$. Both $I$ and $m$ will be determined presently. Let $N_0$ be the lower end-point if $I$. Let $\Sigma$ be a smooth compact LSC immersed hypersurface in $M$ such that $K(\Sigma)=\kappa$; $\Sigma<\hat{\Sigma}$; and $\Sigma$ makes an angle of at least $\theta$ with $\hat{\Sigma}$ along their common boundary. By Theorem \procref{ThmMinimalEmbeddingRadius}, $\Sigma_{p,r}$ is embedded; meets $\partial B_r(p)$ transversally; is contained within $X_r(\hat{\Sigma})$; and divides $X_r(\hat{\Sigma})$ into two connected components, one of which is convex and which we denote by $X_r(\Sigma,\hat{\Sigma})$. We identify $\Sigma$ with $\Sigma_{p,r}$. Let $N$ and $A$ be the upward-pointing unit normal vector field and shape operator of $\Sigma$ respectively. Bearing in mind Proposition \procref{PropUniformModulusOfContinuity}, we suppose that for all $q\in\Sigma$:
$$
D(I,N(q)) \leqslant m(d(p,q)),\eqnum{\nexteqnno}
$$
\noindent where $D$ and $d$ are the distance functions in the total space of $UM$ and in $M$ respectively.
\medskip
\noindent For $q\in\Sigma$, let $e_1,...,e_n$ be an orthonormal basis diagonalising $DK_A$ and let $\mu_1,...,\mu_n$ be its corresponding eigenvalues. We denote:
$$
\mu := \sum_{i=1}^n\mu_i=DK_A(\opId)\geqslant 1.\eqnum{\nexteqnno}
$$
\noindent We define the operator $\Delta^K$ on functions over $\Sigma$ by:
$$
\Delta^K f = \sum_{i=1}^n\mu_i\opHess^\Sigma(f)_{ii},\eqnum{\nexteqnno}
$$
\noindent where $\opHess^\Sigma$ is the Hessian of the Levi-Civita covariant derivative of $\Sigma$. Throughout the sequel, we aim to construct barrier functions that are superharmonic with respect to this generalised Laplacian. The following result plays a fundamental role. It yields a general construction of barrier functions in the non-linear setting.
\proclaim{Lemma \nextprocno}
\noindent Let $\phi:M\rightarrow\Bbb{R}$ be a smooth function such that:
\medskip
\myitem{$(1)$} $\|\nabla\phi\|=1$; and
\medskip
\myitem{$(2)$} the level sets of $\phi$ are LSC with $K$-curvature greater than $\kappa$.
\medskip
\noindent Then the restriction of $\phi$ to $\Sigma$ satisfies:
$$
\Delta^K\phi \geqslant -\|\opHess(\phi)\|\sum_{i=1}^n\mu_i\phi_{;i}\phi_{;i}.
$$
\endproclaim
\proclabel{LemmaSuperHarmonicity}
\proof Choose $q\in\Sigma$. We construct two orthonormal bases for $T_qM$. Let $L_q$ be the level set of $\phi$ passing through $q$. Suppose first that $\nabla\phi$ and $N$ are not colinear at $q$. Then $L_q$ and $\Sigma$ meet transversally at this point. Let $f_1,...,f_{n-1}$ be an orthonormal basis of $T_qL_q\minter T_q\Sigma$ and complete this to an orthonormal basis $f_1,...,f_n$ of $T_q\Sigma$. For $1\leqslant i\leqslant n-1$, denote $f'_i=f_i$, and complete $f'_1,...,f'_{n-1}$ to an orthonormal basis $f'_1,...,f'_{n+1}$ of $T_qM$ such that $f'_n$ is tangent to $L_q$; $f'_n$ makes an angle of at most $\pi/2$ with $f_n$; and $f'_{n+1}$ is normal to $L_q$.
\medskip
\noindent Let $\theta\in]0,\pi/2]$ be the angle between $f_n$ and $f'_n$. Then:
$$
f_n = \opCos(\theta)f'_n \pm \opSin(\theta)f'_{n+1}.
$$
\noindent Let $m_{ij}$ and $m'_{ij}$ be the matrices of the restrictions of $\opHess(\phi)$ to $T_q\Sigma$ and $T_qL_q$ respectively with respect to these bases. Since $\|\nabla\phi\|=1$ and $f'_{n+1}=\pm\nabla\phi$:
$$
\opHess(\phi)(f'_{n+1},\cdot) = 0.
$$
\noindent Consequently:
$$
(m_{ij}) = \pmatrix (m'_{ij})\hfill& \opCos(\theta)(m'_{in})\hfill\cr \opCos(\theta)(m'_{ni})\hfill& \opCos^2(\theta)m'_{nn}\hfill\cr\endpmatrix.
$$
\noindent That is:
$$
(m_{ij}) = \opCos(\theta)(m'_{ij}) + (1-\opCos(\theta))\pmatrix(m'_{ij})\hfill&0\hfill\cr0\hfill&m'_{nn}\hfill\cr\endpmatrix -\opSin^2(\theta)\pmatrix 0\hfill&0\hfill\cr0\hfill&m'_{nn}\hfill\cr\endpmatrix.
$$
\noindent Thus, since $(m'_{ij})$ is positive definite:
$$
(m_{ij}) \geqslant \opCos(\theta)(m'_{ij}) - \opSin^2(\theta)\pmatrix 0\hfill&0\hfill\cr0\hfill&m'_{nn}\hfill\cr\endpmatrix.
$$
\noindent Let $B^{ij}$ be the matrix of $DK_A$ with respect to $f_1,...,f_n$. Then, since $B^{ij}$ is positive definite:
$$
\sum_{i,j=1}^nB^{ij}m_{ij} \geqslant \opCos(\theta)\sum_{i,j=1}^nB^{ij}m'_{ij} - \opSin^2(\theta)B^{nn}m'_{nn}.
$$
\noindent By concavity of $K$, bearing in mind Proposition \procref{PropCurvFnsII}, $(4)$:
$$\matrix
&\sum_{i,j=1}^nB^{ij}(m'_{ij}-A_{ij}) \hfill&\geqslant K(m'_{ij}) - K(A_{ij})\hfill\cr
\Rightarrow\hfill&\sum_{i,j=1}^nB^{ij}m'_{ij}\hfill&\geqslant K(m'_{ij}).\hfill\cr
\endmatrix$$
\noindent However, since $K(L_q)>\kappa(q)$, by Lemma \procref{PropHessianOfRestriction}:
$$
K(m_{ij}') > \kappa(q).
$$
\noindent Thus:
$$
\sum_{i,j=1}^nB^{ij}m_{ij}\geqslant\opCos(\theta)\kappa(q) - \opSin^2(\theta)B^{nn}m'_{nn}.
$$
\noindent However, by Lemma \procref{PropHessianOfRestriction} again:
$$
\opHess^\Sigma(\phi) = \opHess(\phi)|_{T\Sigma} - \langle\nabla\phi,N\rangle A.
$$
\noindent Thus, bearing in mind Proposition \procref{PropCurvFnsII}, $(4)$ again:
$$
\Delta^K\phi\geqslant \opCos(\theta)\kappa(q) - \opCos(\theta)\sum_{i=1}^n\mu_i\lambda_i - \opSin^2(\theta)B^{nn}m'_{nn}=-\opSin^2(\theta)B^{nn}m'_{nn}.
$$
\noindent Finally, $\opSin(\theta)f_n$ is the orthogonal projection onto $T\Sigma$ of $\pm\nabla\phi$. Thus:
$$
\left|m_{nn}'B^{nn}\opSin^2(\theta)\right| \leqslant \|\opHess(\phi)\|\sum_{i=1}^n\mu_i\phi_{;i}\phi_{;i}.
$$
\noindent The result follows in the case where $\nabla\phi$ and $N$ are non-colinear. The case where $\nabla\phi$ and $N$ are colinear follows directly from the concavity of $K$, and this completes the proof.\qed
\medskip
\noindent We require the following modification of this result. For $\phi\in C^\infty(M)$, we define the operator $\Cal{D}_\phi$ over $\Sigma$ by:
$$
\Cal{D}_\phi f = \sum_{i=1}^n\mu_i\phi_{;i}f_{;i}.
$$
\proclaim{Corollary \nextprocno}
\noindent Let $\phi$ be as in Lemma \procref{LemmaSuperHarmonicity}. There exists $\delta, C>0$ in $\Cal{B}(\phi)$ such that:
$$
(\Delta^K + C\Cal{D}_\phi)\phi \geqslant \delta\sum_{i=1}^n\mu_i\opHess(\phi)_{ii}.
$$
\endproclaim
\proclabel{CorSuperHarmonicityI}
\proof In the proof of Lemma \procref{LemmaSuperHarmonicity}, there exists $\delta>0$ in $\Cal{B}$ such that:
$$
\sum_{i,j=1}^n B^{ij}m'_{ij} > \frac{1}{1-\delta}\kappa(q).
$$
\noindent Thus, since $\opSin^2(\theta)B^{nn}m_{nn}'\geqslant 0$:
$$\matrix
&\sum_{i,j=1}^n B^{ij}m_{ij} \hfill&\geqslant \frac{1}{1-\delta}(\opCos(\theta)\kappa(q) - \opSin^2(\theta)B^{nn}m_{nn}')\hfill\cr
\Rightarrow\hfill&\sum_{i,j=1}^n B^{ij}m_{ij}\hfill&\geqslant \delta\sum_{i,j=1}^n B^{ij}m_{ij} + \opCos(\theta)\kappa(q) - \opSin^2(\theta)B^{nn}m'_{nn}.\hfill\cr
\endmatrix$$
\noindent Thus:
$$
\Delta^K\phi \geqslant \delta\sum_{i,j=1}^nB^{ij}\opHess(\phi)_{ij} - \|\opHess(\phi)\|\sum_{i=1}^n\mu_{i}\phi_{;i}\phi_{;i}.
$$
\noindent The result now follows for $C\geqslant \|\opHess(\phi)\|$.\qed
\medskip
\noindent We also require the following straightforward relations. Let $R^M$ be the Riemann curvature tensor of $M$. Let $R^\Sigma$ be the Riemann curvature tensor of $\Sigma$. Let the subscript $;$ denote covariant differentiation with respect to the Levi-Civita covariant derivative of $\Sigma$. We recall the commutation rules of covariant differentiation in a Riemannian manifold:
\proclaim{Proposition \nextprocno}
\myitem{$(1)$} For all $i,j,k$:
$$
A_{ij;k} = A_{kj;i} + R^M_{ki\nu j},
$$
\noindent where $\nu$ represents the direction normal to $\Sigma$; and
\medskip
\myitem{$(2)$} for all $i,j,k,l$:
$$
A_{ij;kl} = A_{ij;lk} + {R^\Sigma_{kli}}^pA_{pj} + {R^\Sigma_{klj}}^pA_{pi}.
$$
\endproclaim
\proclabel{LemmaCommutationRelations}
\proof This is an elementary calculation (c.f. Lemma $6.3$ of \cite{SmiCGC}).\qed
\proclaim{Proposition \nextprocno}
\myitem{$(1)$} For all $p$:
$$
\sum_{i=1}^n\mu_iA_{ii;p} = \kappa_{;p}.
$$
\myitem{$(2)$} For all $p,q$:
$$
\sum_{i=1}^n\mu_iA_{ii;pq} = -(D^2K)^{ij,mn}A_{ij;p}A_{mn;q} + \kappa_{;pq}.
$$
\endproclaim
\proclabel{PropDerivativeOfCurvatureEquation}
\proof This follows by differentiating the equation $K(A)=\kappa$.\qed
\medskip
\noindent Let the subscript $:$ denote covariant differentiation with respect to the Levi-Civita covariant derivative of $M$:
\proclaim{Proposition \nextprocno}
\noindent Let $f$ be the signed distance function to $\Sigma$, let $\nu:=(n+1)$ denote the upward pointing normal direction to $\Sigma$.
\medskip
\myitem{$(1)$} Along $\Sigma$, for all $1\leqslant i,j\leqslant n$:
$$
f_{:ij} = A_{ij},\qquad f_{:i\nu} = f_{:\nu i} = f_{:\nu\nu} = 0,
$$
\noindent where $A$ is the shape operator of $\Sigma$; and
\medskip
\myitem{$(2)$} along $\Sigma$, for all $1\leqslant i,j,k\leqslant n$:
$$
f_{:ijk} = (\nabla^\Sigma A)_{ijk},\qquad f_{:\nu ij} = -A^2_{ij},
$$
\noindent where $\nabla^\Sigma$ is Levi-Civita covariant derivative of $\Sigma$.
\endproclaim
\proclabel{PropSignedDistanceFunction}
\proof This is an elementary calculation (c.f. Lemma $3.16$ of \cite{SmiNLD}).\qed
\newsubhead{Constructing the barrier - part I} The barrier function required to prove Proposition \procref{PropSecondOrderBoundaryEstimates} consists of three components. We construct the first component as follows. Let $X$ be a vector field defined near $p$. Let $f$ be the signed distance to $\Sigma$ in $M$ with sign chosen so that it is positive above $\Sigma$. We define the function $\phi_X$ near $p$ by:
\subheadlabel{SubheadConstructingTheBarrierPartI}
$$
\phi_X = \langle X,\nabla f\rangle.\eqnum{\nexteqnno}
$$
\proclaim{Proposition \nextprocno}
\noindent The restriction of $\phi_X$ to $\Sigma$ satisfies:
$$
\Delta^K\phi_X = O(1)(1+\mu) - \phi_X\sum_{i=1}^n{\mu_i\lambda_i^2},
$$
\noindent where $\lambda_1\geqslant...\geqslant\lambda_n$ and $\mu_1\leqslant...\leqslant\mu_n$ are the eigenvalues of $A$ and $DK_A$ respectively, and $O(1)$ represents terms controlled by $B$, for some $B\in\Cal{B}(X)$.
\endproclaim
\proclabel{LemmaDerOfInnerProduct}
\proof Choose $q\in\Sigma$ and let $e_1,...,e_n$ be an orthonormal basis of $T\Sigma$ at $q$ with respect to which $A$ and $DF_A$ are diagonalised. Let $\lambda_1\geqslant...\geqslant\lambda_n$ and $\mu_1\leqslant...\leqslant\mu_n$ be the corresponding eigenvalues of $A$ and $DF_A$ respectively. We extend $e_1,...,e_n$ to an orthonormal basis of $M$ at $q$ by defining $e_{n+1}=N$. In the sequel, $\nu:=(n+1)$ denotes the upward-pointing normal direction.
\medskip
\noindent We recall that $:$ denotes covariant differentiation with respect to the Levi-Civita covariant derivative of $M$. By Propositions \procref{LemmaCommutationRelations}, $(1)$, and \procref{PropSignedDistanceFunction}, $(2)$, for all $1\leqslant j\leqslant n$:
$$
\sum_{i=1}^n\mu_i f_{:jii} = \sum_{i=1}^n\mu_i(f_{:iij}+R^M_{ij\nu i}) = \kappa_{:j} + O(\mu) = O(1) + O(\mu).
$$
\noindent By Proposition \procref{PropSignedDistanceFunction}, $(2)$, again:
$$
\sum_{i=1}^n\mu_i f_{:\nu ii} = -\sum_{i=1}^n\mu_i\lambda_i^2.
$$
\noindent By Proposition \procref{PropCurvFnsII}, $(4)$:
$$
\sum_{i=1}^n\mu_i\lambda_i = \kappa.
$$
\noindent In particular, for all $i$, since $\mu_i\lambda_i\geqslant 0$, $\lambda_i\mu_i=O(1)$. By Proposition \procref{PropSignedDistanceFunction}, $(1)$, $f_{:ij}=\delta_{ij}\lambda_i$. Thus, recalling that $\|\nabla f\|=1$:
$$\matrix
\sum_{i=1}^n\mu_i\opHess^M(\phi_X)_{ii} \hfill&= \sum_{i=1}^n \mu_i ({X^j}_{:ii}f_{:j} + 2{X^j}_{:i}f_{:ji} + X^jf_{:jii})\hfill\cr
&= O(1)(1+\mu) - X^\nu\sum_{i=1}^n\mu_i\lambda_i^2\hfill\cr
&= O(1)(1+\mu) - \phi_X \sum_{i=1}^n\mu_i\lambda_i^2.\hfill\cr
\endmatrix$$
\noindent Finally, by Lemma \procref{PropHessianOfRestriction}, for any function $h$:
$$
\opHess^\Sigma(h) = \opHess(h)|_{T\Sigma} - \langle\msf{N},\nabla h\rangle A.
$$
\noindent Moreover:
$$
\langle\msf{N},\nabla\phi_X\rangle = {X^k}_{:\nu}f_{:k} + {X^k}f_{:k\nu}.
$$
\noindent By Proposition \procref{PropSignedDistanceFunction}, $(1)$, the second term on the right hand side vanishes along $\Sigma$, and so, using Proposition \procref{PropCurvFnsII}, $(4)$ again:
$$
\sum_{i=1}^n\langle\msf{N},\nabla\phi_X\rangle\mu_i\lambda_i = O(1).
$$
\noindent Thus:
$$
\Delta^K\phi_X =\sum_{i=1}^n\mu_i\opHess^\Sigma(\phi_X)_{;ii}= O(1)(1+\mu) - \phi_X \sum_{i=1}^n\mu_i\lambda_i^2,
$$
\noindent as desired.\qed
\medskip
\noindent The final term on the right-hand side on Proposition \procref{LemmaDerOfInnerProduct} presents an obstacle to the direct use of $\phi_X$ in the construction of the barrier. We remove it by modifying $\phi_X$ as follows. Let $A_\Gamma$ be the shape operator of $\Gamma$, as defined in Section \subheadref{SubheadBoundednessAndMinimalEmbeddingRadii}. Let $N_1$ be the unit vector in $\opN_p\Gamma$ normal to $N_0$ chosen such that $\langle\hat{N},N_1\rangle>0$. Define the interval $J\subseteq\opN_p\Gamma$ by:
$$
J = \left\{ U\in I\ |\ \langle\hat{N},U\rangle>0,\langle N_0,U\rangle>0,\langle N_1, U\rangle>(\kappa(p)/\hat{\kappa})\langle\hat{N},N_1\rangle\right\},
$$
\noindent where $\hat{\kappa}$ is the $K$-curvature of $\hat{\Sigma}$ at $p$. Observe that since $\kappa(p)<\hat{\kappa}$ and since $I$ has length strictly less than $\pi$, $J$ is non-empty. Let $V_p$ be the mid-point of $J$. The reason for such a meticulous choice of $V_p$ will become clear in Propositions \procref{PropPositivityOfPhiV}, \procref{PropLowerCurvBdOfBarrier} and \procref{PropFirstBarrier} below. Let $V$ be a vector field defined near $p$ such that $V(p)=V_p$ and $(\nabla V)(p)=0$.
\proclaim{Proposition \nextprocno}
\noindent There exist $\rho,c>0$ in $B(V,I,m)$ such that, throughout $\Sigma\minter B_\rho(p)$, $\phi_V>c$.
\endproclaim
\proclabel{PropPositivityOfPhiV}
\proof By construction, there exists $\eta>0$ such that every $N'\in I$ makes an angle of at most $\pi/2-2\eta$ with $V_p$. Let $\rho>0$ be such that $m(\rho)<\eta/2$. Upon reducing $\rho$ if necessary, we may assume that $D(V(p),V(q))<\eta/2$ for all $q\in B_\rho(p)$. Then, by definition of $m$ and the triangle inequality, for all $q\in B_\rho(p)$, $D(V(q),N(q))>\pi/2-\eta$, and the result now follows with $c_0:=\opCos(\pi/2-\eta)$.\qed
\medskip
\noindent We define the first order operator $\Cal{D}_1$ on functions over $\Sigma$ by:
$$
\Cal{D}_1h = \frac{2}{\phi_V}\sum_{i=1}^n\mu_i\phi_{V;i}h_{;i},
$$
\noindent and we define the operator $\Cal{L}_1$ by:
$$
\Cal{L}_1 = \Delta^K + \Cal{D}_1.\eqnum{\nexteqnno}
$$
\proclaim{Proposition \nextprocno}
\noindent Using the notation of Proposition \procref{LemmaDerOfInnerProduct}, the restriction of $\phi_X\phi_V^{-1}$ to $\Sigma\minter B_\rho(p)$ satisfies:
$$
\Cal{L}_1(\phi_X\phi_V^{-1}) = O(1)(1+\mu),
$$
\noindent where $O(1)$ represents terms controlled by $B$, for some $B\in\Cal{B}(X,V,I,m)$.
\endproclaim
\proclabel{PropLaplacian}
\proof We use the notation of the proof of Proposition \procref{LemmaDerOfInnerProduct}. By the product rule and the chain rule:
$$\matrix
\Delta^K(\phi_X\phi_V^{-1}) \hfill&=
\sum_{i=1}^n\mu_i(\phi_V^{-1}\opHess^\Sigma(\phi_X)_{ii} + 2(\nabla^\Sigma\phi_X)_i(\nabla^\Sigma\phi_V^{-1})_i + \phi_X\opHess^\Sigma(\phi_V^{-1})_{ii})\hfill\cr
&=\sum_{i=1}^n\mu_i(\phi_V^{-1}\opHess^\Sigma(\phi_X)_{ii} - 2\phi_V^{-2}(\nabla^\Sigma\phi_X)_i(\nabla^\Sigma\phi_V)_i\hfill\cr
&\qquad-\phi_X\phi_V^{-2}\opHess^\Sigma(\phi_V)_{ii} + 2\phi_X\phi_V^{-3}(\nabla^\Sigma\phi_V)_i(\nabla^\Sigma\phi_V)_i)\hfill\cr
&=\phi_V^{-1}(\Delta^K\phi_X)
- 2\phi_V^{-2}\sum_{i=1}^n\mu_i(\phi_X)_{;i}(\phi_V)_{;i}\hfill\cr
&\qquad + 2\phi_X\phi_V^{-3}\sum_{i=1}^n\mu_i(\phi_V)_{;i}(\phi_V)_{;i}
- \phi_X\phi_V^{-2}(\Delta^K\phi_V).\hfill\cr
\endmatrix$$
\noindent Thus, by Proposition \procref{LemmaDerOfInnerProduct}, bearing in mind that $\left|\phi_X\right|,\left|\phi_V\right|,\left|\phi_V^{-1}\right|=O(1)$ over $\Sigma\minter B_\rho(p)$:
$$\matrix
\Delta^K(\phi_X\phi_V^{-1}) \hfill&= O(1)(1 + \mu)-
2\phi_V^{-1}\sum_{i=1}^n\mu_i\opLn(\phi_V)_{;i}(\phi_X)_{;i}\hfill\cr
&\qquad\qquad -2\phi_X\sum_{i=1}^n\mu_i\opLn(\phi_V)_{;i}(\phi_V^{-1})_{;i}\hfill\cr
&=O(1)(1+\mu) - 2\sum_{i=1}^n\mu_i\opLn(\phi_V)_{;i}(\phi_X\phi_V^{-1})_{;i}\hfill\cr
&=O(1)(1+\mu) - \Cal{D}_1(\phi_X\phi_V^{-1}).\hfill\cr
\endmatrix$$
\noindent This completes the proof.\qed
\newsubhead{Constructing the Barrier - Part II}The second component of the barrier function is constructed by taking the signed distance function to a carefully chosen locally strictly concave embedded hypersurface. This hypersurface is constructed as follows. Let $P\subseteq T_pM$ be the hyperplane normal to $N_0$. We identify $P$ with its image under the exponential map in $M$. Observe that $P$ is totally geodesic at $p$. We orient $P$ such that $N_0$ points upward. Observe that $P$ is transverse to $\hat{\Sigma}$ at $P$. Thus, upon reducing $r$ if necessary, we may assume that $P\minter_p B_r(p)$ meets $\hat{\Sigma}$ along a smooth codimension-$2$ submanifold, $\Gamma'$, say of $M$. Let $A_{\Gamma'}$ be the shape operator of $\Gamma'$, as defined in Section \subheadref{SubheadBoundednessAndMinimalEmbeddingRadii}.
\subheadlabel{SubheadConstructingTheBarrierPartII}
\proclaim{Proposition \nextprocno}
\noindent $K_\infty(A_{\Gamma'}(V)) > \kappa(p)$.
\endproclaim
\proclabel{PropLowerCurvBdOfBarrier}
\proof Recall the construction of $V_p$ in the preceeding section. Let $N_1$ be the unit vector in $\opN_p\Gamma$ normal to $N_0$ such that $\langle\hat{N},N_1\rangle>0$. By linearity:
$$
A_{\Gamma'}(\hat{N}) = \langle\hat{N},N_1\rangle A_{\Gamma'}(N_1) + \langle\hat{N},N_0\rangle A_{\Gamma'}(N_0).
$$
\noindent Since $P$ is totally geodesic at $p$, $A_{\Gamma'}(N_0)=0$. Thus, $A_{\Gamma'}(\hat{N})=\langle\hat{N},N_1\rangle A_{\Gamma'}(N_1)$. Likewise, $A_{\Gamma'}(V)=\langle V,N_1\rangle A_{\Gamma'}(N_1)$. Combining these relations yields:
$$
A_{\Gamma'}(V) = \frac{\langle V,N_1\rangle}{\langle\hat{N},N_1\rangle}A_{\Gamma'}(\hat{N}).
$$
\noindent Thus, bearing in mind Proposition \procref{PropCharacterisationOfFInfinity}:
$$
K_\infty(A_{\Gamma'}(V)) = \frac{\langle V,N_1\rangle}{\langle\hat{N},N_1\rangle}K_\infty(A_{\Gamma'}(\hat{N})) > (\kappa(p)/\hat{\kappa})\hat{\kappa} = \kappa(p),
$$
\noindent as desired.\qed
\medskip
\noindent By Proposition \procref{PropLowerCurvBdOfBarrier}, we may extend $\Gamma'$ to a smooth LSC embedded hypersurface $H_0$ passing through $p$ such that the upward pointing normal to $H_0$ at $p$ is $V(p)$ and  $K(H_0)(p)>\kappa(p)$. Upon perturbing $H_0$ slightly, and reversing the orientation, we obtain a smooth locally strictly {\sl concave} embedded hypersurface $H$ passing through $p$ such that the upward pointing normal to $H$ at $p$ is $-V(p)$; if $A$ is the shape operator of $H$ at $p$, then $K(-A)>\kappa(p)$; and $H_0\minter_p B_r(p)$ lies above the graph of $\epsilon d_p^2$ over $H$, for some $\epsilon>0$.
\medskip
\noindent We will show that $H$ is the desired hypersurface. Let $d_H$ be the signed distance to $H$ in $M$ with sign chosen so that $d_H$ is positive above $H$. In this section, we prove superharmonicity of $d_H$. Later we prove non-negativity of the restriction of $d_H$ to $\partial\Sigma$. This will depend upon a suitable choice of $N_0$. For $C>0$, define the first order operator $\Cal{D}_2$ on functions over $\Sigma$ such that:
$$
\Cal{D}_2h = -C\sum_{i=1}^n\mu_id_{H;i}h_{;i}.
$$
\noindent We define $\Cal{L}_2$ by:
$$
\Cal{L}_2 = \Delta^K + \Cal{D}_1 + \Cal{D}_2.
$$
\noindent $C$ is determined by the following result:
\proclaim{Proposition \nextprocno}
\noindent Using the notation of Proposition \procref{LemmaDerOfInnerProduct}, there exists $\rho,C\in\Cal{B}(V,I,m,H)$ such that the restriction of $d_H$ to $\Sigma\minter B_\rho(p)$ satisfies:
$$
\Cal{L}_2 d_H\leqslant 0.
$$
\endproclaim
\proclabel{PropFirstBarrier}
\proof Choose $\rho>0$ such that the $K$-curvature of every level set of $d_H$ in $B_\rho(p)$ is strictly greater than $\kappa$. Bearing in mind that $-d_H$ is convex, by Corollary \procref{CorSuperHarmonicityI}, there exists $C,\delta_a>0$ in $\Cal{B}(H)$ such that:
$$
(\Delta^K + \Cal{D}_2)d_H \leqslant \delta_a\sum_{i=1}^n\mu_i\opHess(d_H)_{ii}.
$$
\noindent However, since $\nabla d_H(p)=-V(p)$, as in Proposition \procref{PropPositivityOfPhiV}, there exists $\delta_b>0$ in $\Cal{B}(I,m,H)$ such that, upon reducing $\rho$ if necessary, $\langle\nabla d_H,N\rangle<-\delta_b$ throughout $\Sigma\minter B_\rho(p)$. Since $H$ is strictly convex, there therefore exists $\delta_c>0$ in $\Cal{B}(I,m,H)$ such that throughout $\Sigma\minter B_\rho(p)$ and for all $1\leqslant i\leqslant n$:
$$
\opHess(d_H)_{ii} \leqslant -\delta_b.
$$
\noindent Thus, denoting $\delta_d=\delta_a\delta_c$:
$$
(\Delta^K + \Cal{D}_2)d_H \leqslant -\delta_d\mu.\eqnum{\nexteqnno}
$$
\noindent It remains to consider the contribution from $\Cal{D}_1$. As in the proof of Proposition \procref{LemmaDerOfInnerProduct}:
$$
(\phi_V)_{;i} = {V^k}_{:i}f_k + V^kf_{:ki} = {V^\nu}_{:i} + \lambda_iV^i.
$$
\noindent Since $\nabla V(p)=0$, bearing in mind Proposition \procref{PropPositivityOfPhiV}, upon reducing $\rho$ if necessary, we may assume that:
$$
\left|2\phi_V^{-1}\sum_{i=1}^n\mu_i{V^\nu}_{:i}d_{H;i}\right| \leqslant \frac{\delta_d\mu}{2}.\eqnum{\nexteqnno}
$$
\noindent Moreover, since $(\nabla d_H + V)(p)=0$ at $p$, after reducing $\rho$ further if necessary, we may assume that:
$$
\|\nabla d_H + V\|\leqslant \frac{\delta_d\phi_V\mu}{4\kappa\|V\|}.
$$
\noindent Thus, bearing in mind Proposition \procref{PropCurvFnsII}, $(4)$, and the fact that $\mu_i\lambda_i\geqslant 0$ for all $i$:
$$\matrix
2\phi_V^{-1}\sum_{i=1}^n\mu_i\lambda_iV^id_{H;i} \hfill&= -2\phi_V^{-1}\sum_{i=1}^n\mu_i\lambda_iV^iV^i\hfill\cr
&\qquad + 2\phi_V^{-1}\sum_{i=1}^n\mu_i\lambda_i V^i(d_{H;i} + V^i)\hfill\cr
&\leqslant \frac{\delta_d\mu}{2}.\hfill\cr
\endmatrix\eqnum{\nexteqnno}$$
\noindent Combining these relations yields:
$$
\Cal{L}_2d_H\leqslant 0,
$$
\noindent as desired.\qed
\medskip
\noindent We now verify that the addition of the term $\Cal{D}_2$ does not affect the conclusion of Proposition \procref{PropLaplacian}. Indeed:
\proclaim{Proposition \nextprocno}
\noindent Using the notation of Proposition \procref{PropLaplacian}, the restriction of $\phi_X\phi_V^{-1}$ to $\Sigma\minter B_\rho(p)$ satisfies:
$$
\Cal{L}_2(\phi_X\phi_V^{-1}) = O(1)(1+\mu),
$$
\noindent where $O(1)$ represents terms controlled by $B$ for some $B$ in $\Cal{B}(X,V,I,m,H)$.
\endproclaim
\proclabel{PropYetAnotherSuperharmonicityRelation}
\proof Indeed, for any vector field, $Y$:
$$
\phi_{Y;i} = {Y^\nu}_{;i} + Y_i\lambda_i.
$$
\noindent It follows by Proposition \procref{PropCurvFnsII}, $(4)$, and the fact that $\lambda_i\mu_i\geqslant 0$ for all $i$ that:
$$
\Cal{D}_2(\phi_X\phi_V^{-1}) = O(1)(1+\mu).
$$
\noindent The result now follows by Proposition \procref{PropLaplacian}.\qed
\newsubhead{Constructing the Barrier - Part III} The third component of the barrier function is simply the squared distance to $p$ in $M$:
\subheadlabel{ConstrucingTheBarrierPartIII}
\proclaim{Proposition \nextprocno}
\noindent There exists $\epsilon>0$ in $\Cal{B}(V,I,m,H)$ such that, after reducing $\rho$ if necessary, the restriction of $d_p^2$ to $\Sigma\minter B_\rho(p)$ satisfies:
$$
\Cal{L}_2 d_p^2 \geqslant \frac{1}{2}(1+\mu).
$$
\endproclaim
\proclabel{PropSecondBarrier}
\proof We continue to use the notation of the proof of Proposition \procref{LemmaDerOfInnerProduct}. Since $\mu\geqslant 1$, by Proposition \procref{PropCurvFnsII}, $(4)$:
$$
\Delta^K d_p^2 \geqslant (1 + \mu) - 2d_p\langle\nabla d_p,N\rangle\kappa.
$$
\noindent For $\rho<1/8$, throughout $\Sigma\minter B_\rho(p)$:
$$
2d_p\langle\nabla d_p,\msf{N} \rangle\kappa < 1/4.
$$
\noindent Thus, throughout $B_\rho(p)$:
$$
\Delta^K(d_p^2) \geqslant \frac{3}{4}(1+\mu).
$$
\noindent Reducing $\rho$ further if necessary, by Proposition \procref{PropPositivityOfPhiV}:
$$
\left|4d_p\phi_V^{-1}\sum_{i=1}^n\mu_i{V^\nu}_{;i}d_{p;i}\right| \leqslant \frac{1}{16}\mu.
$$
\noindent Moreover, upon reducing $\rho$ further if necessary, for all $i$, by Proposition \procref{PropCurvFnsII}, $(4)$, bearing in mind that $\mu_i\lambda_i\geqslant 0$ for all $i$:
$$
\left|4d_p\phi_V^{-1}\mu_i\lambda_iV^id_{p;i}\right| \leqslant \frac{\mu}{16n}
$$
\noindent Combining these relations yields:
$$
\left|\Cal{D}_1d_p^2\right| \leqslant \frac{1}{8}(1+\mu).
$$
\noindent In like manner, after reducing $\rho$ yet further if necessary:
$$
\left|\Cal{D}_2d_p^2\right| \leqslant \frac{1}{8}(1+\mu).
$$
\noindent Thus:
$$
\Cal{L}_2(d_p^2) \geqslant \frac{1}{2}(1+\mu),
$$
\noindent as desired.\qed
\medskip
\noindent We now prove Proposition \procref{PropSecondOrderBoundaryEstimates}:
\medskip
{\noindent\bf Proof of Proposition \procref{PropSecondOrderBoundaryEstimates}:\ }We assume the contrary. Let $(\hat{\Sigma}_n,\kappa_n)_\ninn$ be a sequence in $\Cal{F}$ converging to $(\hat{\Sigma},\kappa)$. For convenience, we suppose that $\hat{\Sigma}_n=\hat{\Sigma}$ and $\kappa_n=\kappa$ for all $n$. Let $\Gamma$ be the boundary of $\hat{\Sigma}$ and let $\hat{N}$ be its upward-pointing unit normal vector field. Let $(\Sigma_n)_\ninn$ be a sequence of smooth compact LSC immersed hypersurfaces such that for all $n$, $K(\Sigma_n)=\kappa$; $\Sigma_n<\hat{\Sigma}$; and $\Sigma$ makes an angle of at least $\theta$ with $\hat{\Sigma}$ along their common boundary. For all $n$, let $N_n$ be the upward-pointing unit normal vector field over $\Sigma_n$ and let $A_n$ be its shape operator. Let $(p_n)_\ninn$ be a sequence of points of $\Gamma$ converging to $p$ and suppose that $(\|A_n(p_n)\|)_\ninn$ converges to $+\infty$. For convenience, we suppose that $p_n=p$ for all $n$.
\medskip
\noindent We henceforth identify $\Sigma_n$, $\hat{\Sigma}$ and $\Gamma$ with $\Sigma_{n,p,r}$, $\hat{\Sigma}_{p,r}$ and $\Gamma_{p,r}$ respectively. For all $n$, we denote $X_n=X_{p,r}(\Sigma_n,\hat{\Sigma})$. Observe that, for all $n$, $X_n$ is contained in a ball of radius $r$. Furthermore, by Theorem \procref{ThmMinimalEmbeddingRadius} there exists $s>0$ such that for all $n$, $X_n$ contains a ball of radius $s$. Thus, upon extracting a subsequence, there exists a convex subset $X$ with non-trivial interior towards which $(X_n)_\ninn$ converges in the Hausdorff sense.
\medskip
\noindent Observe that $p$ is a boundary point of $X$. Let $I$ be the set of supporting normal vectors to $X$ at $p$. By Proposition \procref{PropSupportingNormalIsContainedBelowOuterBarrier}, $I$ is a subset of the closure of $\opCN_p^-\hat{\Sigma}$. By Proposition \procref{PropUniformModulusOfContinuity}, upon extracting a subsequence, we may suppose that there exists a sequence $(I_n)_\ninn$ of closed subintervals of $\opCN_p^-\hat{\Sigma}$ and a continuous function $m:[0,\infty[\rightarrow[0,\infty[$ such that $(I_n)_\ninn$ converges to $I$ in the Hausdorff sense; $m(0)=0$; and for all $n$ and for all $q\in\Sigma_n$, $D(N_n(q),I_n) \leqslant m(d(q,p))$. In particular, for all $n$, $N_n(p)$ is an element of $I_n$. For all $n$, let $N_{n,0}$ be the lower end-point of $I_n$, and with this choice of $I_n$, define $V$, $P$, $H$ and $d_H$ as above.
\medskip
\noindent We claim that, upon reducing $r$ if necessary, there exists $\delta_a>0$ such that, for all sufficiently large $n$, $d_H\geqslant\delta_a d_p^2$ along $(\partial\Sigma_n)\minter B_\rho(p)=\Gamma\minter B_\rho(p)$ and $d_H\geqslant\delta_a$ along $\Sigma_n\minter(\partial B_\rho(p))$. Indeed, let $Y$ be the convex component of the complement of $H_0\minter_p B_\rho(p)$ in $B_\rho(p)$. Observe that $P\minter\hat{\Sigma}$ divides $\hat{\Sigma}$ into two connected components. We denote by $\hat{\Sigma}^-$ the component lying to the left of this hypersurface. Since $V$ lies below $\hat{N}_p$, $\hat{\Sigma}^-$ is contained in $Y$. Likewise, $P\minter\hat{\Sigma}$ divides $P$ into two connected components. We denote by $P^+$ the component lying to the right of this hypersurface. Since $V$ lies above $N_0$, $P^+$ is also contained in $Y$. We denote by $Y_1$ the closure of the convex component of the complement of $\hat{\Sigma}$ in $B_\rho(p)$, and by $Y_2$ the closure of the connected component of the complement of $P$ in $B_\rho(p)$ lying below $P$. Since the boundary of $Y_1\minter Y_2$ coincides with the union of $\hat{\Sigma}^-$ with $P^+$, $Y_1\minter Y_2\subseteq Y$. However, by construction, $X\subseteq Y_1$. Furthermore, since $N_0$ is a supporting normal to $X$ at $p$, by convexity $X\subseteq Y_2$. It follows that $X$ lies in $Y$. In particular, $\Gamma$ is contained in $Y$, and so, by definition of $H$, for sufficiently small $\delta_a>0$, $d_H\geqslant\delta_a d_p^2$ along $\Gamma\minter B_\rho(p)$ as asserted. Likewise, upon reducing $\delta$ further if necessary, $d_H(q)\geqslant 2\delta_a$ for all $q\in X\minter (\partial B_\rho(p))$. Thus, for sufficiently large $n$, $d_H(q)\geqslant\delta_a$ for all $q\in\Sigma_n\minter(\partial B_\rho(p))\subseteq X_n\minter(\partial B_\rho(p))$, as asserted.
\medskip
\noindent Now let $X$ be any vector field over $B_r(p)$ which is tangent along $\Gamma$. Denote $\phi=\phi_X\phi_V^{-1}$. Observe that $\phi$ vanishes along $\Gamma$. By Propositions \procref{PropYetAnotherSuperharmonicityRelation} and \procref{PropSecondBarrier}, there exists $A_->0$ in $\Cal{B}(X,V,I,m,H)$ such that, throughout $\Sigma_n\minter B_\rho(p)$:
$$
\Cal{L}_2(\phi-A_-d_p^2) < 0.
$$
\noindent Bearing in mind Proposition \procref{PropFirstBarrier} and the preceeding paragraph, upon reducing $\rho$ if necessary, there therefore exists $B_->0$ in $\Cal{B}(X,V,I,m,H)$ such that:
\medskip
\myitem{$(1)$} $\Cal{L}_2(\phi+B_-d_{H} - A_-d_p^2)< 0$ throughout $\Sigma_n\minter B_\rho(p)$; and
\medskip
\myitem{$(2)$} $\phi + B_-d_H - A_-d_p^2 \geqslant 0$ along $\partial(\Sigma_n\minter B_\rho(p))$.
\medskip
\noindent It thus follows by the maximum principle that, throughout $\Sigma_n\minter B_\rho(p)$:
$$
\phi \geqslant A_-d_p^2 - B_-d_{H}.
$$
\noindent Likewise, reducing $\rho$ further if necessary, there exists $A_+$ and $B_+$ in $\Cal{B}(X,V,I,m,H)$ such that, throughout $\Sigma_n\minter B_\rho(p)$:
$$
\phi\leqslant  - A_+d_p^2 + B_+ d_{H}.
$$
\noindent We thus obtain a-priori bounds on $d\phi$ at $p$. Let $f$ be the signed distance function to $\Sigma_n$. For all $Y$, since $\phi_X(p)=0$:
$$
\opHess(f)(X,Y) = \langle\nabla\phi,Y\rangle\phi_V(p) - \langle\nabla_Y X,\msf{N}\rangle.
$$
\noindent Thus, since $X$ is arbitrary, we obtain a-priori bounds on $\opHess(f)(X,Y)$ for all pairs of vectors $X,Y\in T_p\Sigma_n$ where at least one of $X$ or $Y$ is tangent to $\partial\Sigma_n$. By Lemma \procref{PropHessianOfRestriction}, the second fundamental form of $\Sigma_n$ is the restriction to $T\Sigma_n$ of the hessian of $f$, we deduce that there exists $B$ in $\Cal{B}(X,V,I,m,H)$ such that:
$$
\|A_n(X,Y)\|(p)\leqslant B\|X\|\|Y\|,
$$
\noindent for all $n$ and for all such pairs of vectors. Since, by hypotheses, $\|A_n(p)\|\rightarrow+\infty$, it follows that $\|A_n(X_n,X_n)\|\rightarrow+\infty$ where, for all $n$, $X_n$ is the unit vector normal to $\partial\Sigma_n$ in $T_p\Sigma_n$.
\medskip
\noindent However, we may assume that $(X_n)_\ninn$ converges to $X_\infty$, say, which is normal to $\Gamma$ at $p$. Let $\lambda_1'\leqslant...\leqslant\lambda_{n-1}'$ be the eigenvalues of $A_\Gamma(X_\infty)$. For all $m$, let $\lambda_{1,m}\leqslant...\leqslant\lambda_{n,m}$ be the eigenvalues of $A_n$. By the above discussion, $(\lambda_{n,m})_{m\in\Bbb{N}}\rightarrow+\infty$. By Lemma $1.2$ of \cite{CaffNirSprIII} and the bounds already obtained, for all $1\leqslant i\leqslant n-1$:
$$
(\lambda_{i,m})_{m\in\Bbb{N}}\rightarrow\lambda_i'.
$$
\noindent Suppose first that $K$ is of bounded type. By Proposition \procref{PropFirstOrderEstimatesUnboundedCase}, $\lambda_i'>0$ for all $i$. By Proposition \procref{PropFirstOrderLowerEstimatesBoundedCase}, $K_\infty(\lambda_1',...,\lambda_{n-1}')>\kappa(p)$. However, by concavity, $K(x_1,...,x_{n-1},t)$ converges locally uniformly to $K_\infty(x_1,...,x_{n-1})$ in $(x_1,...,x_{n-1})$ as $t\rightarrow+\infty$. That is:
$$
\mlim_{m\rightarrow+\infty}K(\lambda_{1,m},...,\lambda_{n,m})
=K_\infty(\lambda'_1,...,\lambda'_{n-1})>\kappa(p_0),
$$
\noindent which is absurd.
\medskip
\noindent Suppose now that $K$ is of unbounded type. By Proposition \procref{PropFirstOrderEstimatesUnboundedCase}, $\lambda_i'>0$ for all $i$. In the same manner, we obtain:
$$
\mlim_{m\rightarrow+\infty}K(\lambda_{1,m},...,\lambda_{n,m})=+\infty>\kappa(p),
$$
\noindent which is likewise absurd. There therefore exists $B_2\geqslant 0$ such that, for all $n$:
$$
\|A_n(p_n)\|\leqslant B_2,
$$
\noindent which is absurd, and this completes the proof.\qed
\newhead{Global Second Order Estimates}%
\headlabel{HeadGlobalSecondOrderEstimates}%
\newsubhead{Main results} Let $M:=M^{n+1}$ be an $(n+1)$-dimensional Riemannian manifold. Let $K$ be a convex curvature function. Let $\kappa:M\rightarrow]0,\infty[$ be a smooth positive function. In this section we obtain global a-priori estimates for the shape operators of smooth compact LSC immersed hypersurfaces in $M$ of $K$-curvature prescribed by $\kappa$. Let $X$ be a compact geodesically convex subset of $M$. As in Section \subheadref{SubheadSecondOrderBoundaryEstimatesPreliminaryResults}, we denote by $\Cal{B}$ the family of all quantities which only depend on the restrictions to $X$ of $\kappa$ and the metric of $M$.
\subheadlabel{SubheadGlobalSecondOrderEstimatesMainResults}
\proclaim{Proposition \nextprocno}
\noindent Choose $R>0$. Suppose that $\kappa(p)<\mu_\infty(K)/R$ for all $p\in X$. There exists $B>0$ in $\Cal{B}(R)$ with the property that if $\Sigma$ is a smooth compact LSC immersed hypersurface in $M$ such that $K(\Sigma)=\kappa$ and $\Sigma\subseteq X\minter B_R(q)$, for some $q\in M$, then, for all $p\in\Sigma$:
$$
\|A(p)\| \leqslant B(1 + \msup_{p'\in\partial\Sigma}\|A(p')\|).
$$
\endproclaim
\proclabel{PropSecondOrderIntBoundsAlt}
\noindent Using a slightly different approach, we obtain the following complementary result:
\proclaim{Proposition \nextprocno}
\noindent Suppose that the sectional curvature of $M$ is bounded above by $-1$ and that $\kappa(p)\in]0,1[$ for all $p\in X$. There exists $B>0$ in $\Cal{B}$ with the property that if $\Sigma$ is a smooth compact LSC immersed hypersurface in $M$ such that $K(\Sigma)=\kappa$ and $\Sigma\subseteq X$, then, for all $p\in\Sigma$:
$$
\|A(p)\| \leqslant B(1 + \msup_{p'\in\partial\Sigma}\|A(p')\|).
$$
\endproclaim
\proclabel{PropSecondOrderIntBounds}
\newsubhead{Asymptotic behaviour of curvature functions} We briefly clarify the hypotheses of Proposition \procref{PropSecondOrderIntBoundsAlt}. Observe that for Plateau problems in $\Bbb{R}^n$, we would want open subsets of the sphere of radius $R$ to serve as barriers for hypersurfaces of constant curvature equal to $k$ for all $k\in]0,1/R]$, and in particular for $k=1/R$. This is guaranteed by the following result.
\proclaim{Proposition \nextprocno}
\noindent $\mu_\infty(K)>1$.
\endproclaim
\proclabel{PropMuInfinityIsGreaterThanOne}
\remark This follows from the fact that $K$ is $C^1$ over the interior of $\Gamma$.\qed
\medskip
\proof Denote $\Bbb{I}:=(1,...,1)$. Let $X\subseteq\Gamma$ be the set of all points $(x_1,...,x_n)$ such that $0\leqslant x_i\leqslant 1$ for all $i$, and $x_i=1$ for at least one $i$. Since $K$ is homogeneous of order $1$, $DK_x(\Bbb{I})$ is homogeneous of order $0$. It thus suffices to show that:
$$
\mliminf_{x\in X,x\rightarrow\partial\Gamma}DK_x(\Bbb{I}) > 1.
$$
\noindent Suppose the contrary. By Proposition \procref{PropCurvFnsII}, $(6)$, $DK_x(\Bbb{I})\geqslant 1$ for all $x$. We suppose therefore that there exists a sequence $(x_n)_\ninn$ in $\partial X$ converging to $x_\infty\in\Gamma$ such that for all $n$, $DK_{x_n}(\Bbb{I})=1$. We claim that $DK_\Bbb{I}(x_n)=K(x_n)$ for all $n$. Indeed, by Proposition \procref{PropCurvFnsII} $(4)$, for all $x\in X$, $DK_x(x)=K(x)$. Thus, if $DK_x(\Bbb{I})=1$, then $DK_x(\Bbb{I}-x) = 1 - K(x)$. However, by concavity, for all $t\in[0,1]$, $DK_{t\Bbb{I}+(1-t)x}(\Bbb{I}-x)\leqslant DK_x(\Bbb{I}-x)$. Thus, for all $t\in[0,1]$, $DK_{t\Bbb{I} + (1-t)x}(\Bbb{I}-x)\leqslant 1 - K(x)$. Since the integral equals $1-K(x)$, it follows that for all $t\in[0,1]$, $DK_{t\Bbb{I}+(1-t)x}(\Bbb{I}-x)=1-K(x)$. In particular, since $K$ is $C^1$, $DK_\Bbb{I}(\Bbb{I} - x)=1-K(x)$. Thus, if $DK_x(\Bbb{I})=1$, then $DK_\Bbb{I}(x)=K(x)$. In particular, $DK_\Bbb{I}(x_n)=K(x_n)$ for all $n$, as asserted. By compatibility, taking limits yields $DK_\Bbb{I}(x_\infty)=K(x_\infty)=0$. However, since $x_\infty\in\partial X$, $DK_\Bbb{I}(x_\infty)=x_{\infty,1}/n+...+x_{\infty,n}/n\geqslant 1/n$. This is absurd, and the result follows.\qed
\medskip
\noindent When $K$ is of unbounded type, the hypotheses of Proposition \procref{PropSecondOrderIntBoundsAlt} are trivially satisfied for sufficiently large $R$. Indeed:
\proclaim{Proposition \nextprocno}
\noindent If $K$ is of unbounded type, then $\mu_\infty(K)=\infty$.
\endproclaim
\proclabel{PropMuInfinityIsInfinite}
\proof Indeed, choose $B>0$. Since $K$ is of unbounded type, there exists $C>0$ such that $K(1,...,1,C)\geqslant B\sqrt{n}$. For $c\geqslant 0$, let $X_c$ be the set of all points $(x_1,...,x_n)$ such that $c\leqslant x_i\leqslant C$ for all $i$ and $x_i=C$ for at least one $i$. By homogeneity, it suffices to show that there exists a neighbourhood $U$ of $\partial X_0$ in $X_0$ such that for all $x\in U$, $DK_x(\Bbb{I})\geqslant B$. By invariance and ellipticity of $K$, $K\geqslant B\sqrt{n}$ throughout $X_1$. However, by compatibility, $K$ vanishes over $\partial X_0$. For every point $x\in\partial X_0$, let $y_x$ be such that $x+y_x$ is the closest point to $x$ in $\partial X_1$. Observe that for all $x\in\partial X_0$, $y_x\in[0,1]^n$ and $y_{x,i}=1$ for at least one $i$. By the intermediate value theorem, for all $x\in\partial X_1$, there exists $t_x\in[0,1]$ such that $DK_{t_xy_x}(y_x)\geqslant B\sqrt{n}$. By concavity, for all $t\in]0,t_x[$, $DK_{ty_x}(y_x)\geqslant DK_{t_xy_x}(y_x)\geqslant B\sqrt{n}$. In particular, $\|DK_{ty_x}\|\geqslant B$ for all such $t$. By ellipticity, every component of $DK_{ty_x}$ is positive. Thus $DK_{ty_x}(\Bbb{I})\geqslant\|DK_{ty_x}\|\geqslant B$. Since $y_x$ and $t_x$ may be chosen to vary continuously with $x\in\partial X_0$, there exists a neighbourhood $U$ of $\partial X_0$ in $X_0$ such that for all $x\in U$, $DK_x(\Bbb{I})\geqslant B$, as desired.\qed
\newsubhead{Preliminary results}Let $\Sigma$ be a smooth compact LSC immersed hypersurface in $M$ such that $K(\Sigma)=\kappa$ and $\Sigma\subseteq X$. Let $N$ be the upward-pointing unit normal vector field over $\Sigma$. Let $A$ be the shape operator of $\Sigma$. Observe that, since $\Sigma$ is LSC, $A$ is everywhere positive definite. Let $R^M$ be the Riemann curvature tensor of $M$. Let $R^\Sigma$ be the Riemann curvature tensor of $\Sigma$. Let the subscripts $:$ and $;$ denote covariant differentiation with respect to the Levi-Civita covariant derivatives of $M$ and $\Sigma$ respectively. Let the subscript $\nu:=n+1$ denote the upward-pointing normal direction.
\subheadlabel{SubheadGlobalSecondOrderBoundsPreliminaryResults}
\proclaim{Proposition \nextprocno}
\noindent For all $i,j,k$ and $l$:
$$
A_{ij;kl} = A_{kl;ij} + R^M_{kj\nu i;l} + R^M_{li\nu k;j} + {R^\Sigma_{jlk}}^pA_{pi} + {R^\Sigma_{jli}}^pA_{pk}.
$$
\endproclaim
\proclabel{CorSecondDerivativesOfA}
\proof This follows immediately from Proposition \procref{LemmaCommutationRelations} (c.f. Corollary $6.4$ of \cite{SmiCGC}).\qed
\medskip
\noindent Choose $p\in\Sigma$. Let $e_1,...,e_n$ be an orthonormal basis of $T_p\Sigma$ diagonalising $A$ and $DK_A$. Let $\lambda_1\geqslant...\geqslant\lambda_n$ and $\mu_1\leqslant...\leqslant\mu_n$ be the corresponding eigenvalues of $A$ and $DK_A$ respectively. We extend $\lambda_1$ to a continuous function near $p$ such that for all $q$, $\lambda_1(q)$ is the greatest eigenvalue of $A(q)$. We define the operator $\Delta^K$ as in Section \subheadref{SubheadSecondOrderBoundaryEstimatesPreliminaryResults}. For $\phi\in C^\infty(M)$, we define the homogeneous first order operator $\Cal{D}_\phi$ on functions over $\Sigma$ by:
$$
\Cal{D}_\phi f=\sum_{i=1}^n\mu_i\phi_{;i}f_{;i}.
$$
\noindent We define $I,J\subseteq\left\{1,...,n\right\}$ by:
$$
I=\left\{1\leqslant i\leqslant n\ |\ \mu_i\leqslant 4\mu_1\right\},\qquad
J=\left\{1\leqslant i\leqslant n\ |\ \mu_i>4\mu_1\right\}.
$$
\noindent We recall that a function $f$ is said to satisfy $(\Delta^K + C\Cal{D}_\phi)f\geqslant g$ in the weak sense at $p$ if and only if there exists a smooth function $a$, defined near $p$ such that $f\geqslant a$ near $p$; $f(p)=a(p)$ at $p$; and $(\Delta^K + CD_\phi)a\geqslant g$ at $p$. We obtain:
\proclaim{Proposition \nextprocno}
\noindent For all $\phi\in C^\infty(M)$, $C\geqslant 0$, there exists $K\geqslant 0$ in $\Cal{B}(\phi,C)$ such that if $\lambda_1(p)\geqslant 1$, then, at $p$:
$$
(\Delta^K + C\Cal{D}_\phi)\opLog(\lambda_1) \geqslant -K(1+\mu) - \sum_{i\in I}\frac{\mu_i}{\lambda_1^2}A_{11;i}^2,
$$
\noindent in the weak sense.
\endproclaim
\proclabel{CorSuperharmonic}
\noindent We extend $e_1,...,e_n$ to a frame near $p$ by parallel transport along geodesics. We define the function $a$ near $p$ by $a=A(e_1,e_1)$. Observe that $\lambda_1\geqslant a$ and $\lambda_1(p)=a(p)$.
\proclaim{Proposition \nextprocno}
\noindent For all $i$, at $p$, $a_{;i} = A_{11;i}$ and $a_{;ii} = A_{11;ii}$.
\endproclaim
\proclabel{PropSecondDerOfFirstEval}
\proof This is an elementary calculation (c.f. Proposition $6.5$ of \cite{SmiCGC}).\qed
\proclaim{Proposition \nextprocno}
\noindent For all $\phi\in C^\infty(M)$, $C\geqslant 0$, there exists $K>0$ in $\Cal{B}(\phi,C)$ such that, if $a(p)\geqslant 1$, then, at $p$:
$$
(\Delta^K + C\Cal{D}_\phi)\opLog(a)(p) \geqslant - K(1 + \mu) - \sum_{i\in I}\frac{\mu_i}{\lambda_1^2}A_{11;i}^2.
$$
\endproclaim
\proclabel{PropSuperharmonic}
\proof By Propositions \procref{CorSecondDerivativesOfA} and \procref{PropSecondDerOfFirstEval}:
$$
a_{;ii} = A_{11;ii}= A_{ii;11} + R^M_{i1\nu 1;i} + R^M_{i1\nu i;1} + {R^\Sigma_{1ii}}^pA_{p1} + {R^\Sigma_{1i1}}^pA_{pi}.
$$
\noindent However, at $p$, by Proposition \procref{PropDerivativeOfCurvatureEquation}, $(2)$:
$$
\sum_{i=1}^n\frac{\mu_i}{\lambda_1}A_{ii;11} = -\frac{1}{\lambda_1}(D^2K)^{ij,mn}A_{ij;1}A_{mn;1} + \frac{1}{\lambda_1}\kappa_{;11}.
$$
\noindent Thus, at $p$:
$$\matrix
\Delta^K\opLog(a) \hfill&= \frac{1}{\lambda_1}\kappa_{;11}
-\frac{1}{\lambda_1}(D^2K)^{ij,mn}A_{ij;1}A_{mn;1} - \sum_{i=1}^n\frac{\mu_i}{\lambda_1^2}A_{11;i}A_{11;i}\hfill\cr
&\qquad + \sum_{i=1}^n\frac{\mu_i}{\lambda_1}(R^M_{i1\nu 1;i}+R^M_{i1\nu i;1}) + \sum_{i,j=1}^n\frac{\mu_i}{\lambda_1}({R^\Sigma_{1ii}}^pA_{p1} + {R^\Sigma_{1i1}}^pA_{pi}).\hfill\cr
\endmatrix$$
\noindent We consider each contribution seperately. Since, for all $a,b\in\Bbb{R}$ and for all $\eta>0$, $(a+b)^2\leqslant (1+\eta)a^2+(1+\eta^{-1})b^2$, by Lemma \procref{LemmaCommutationRelations}, $(1)$, there exists $K_1$ in $\Cal{B}$ such that for all $i\in J$:
$$
\frac{9}{8}A_{11;i}^2 = \frac{9}{8}(A_{i1;1} + R^M_{i1\nu 1})^2 \leqslant \frac{5}{4}A_{i1;1}^2 + K_1.
$$
\noindent Thus, by Proposition \procref{PropCurvFnsII}, $(7)$, bearing in mind the definition of $J$ and the fact that $\lambda_1\geqslant 1$:
$$\matrix
-\frac{1}{\lambda_1}(D^2K)^{ij,mn}A_{ij;1}A_{mn;1} - \frac{9}{8}\sum_{i\in J}\frac{\mu_i}{\lambda_1^2}A_{11;i}A_{11;i}\hfill\cr
\qquad\qquad\qquad\geqslant \sum_{i\in J}(\frac{2(\mu_i-\mu_1)}{\lambda_1(\lambda_1-\lambda_i)} - \frac{5}{4}\frac{\mu_i}{\lambda_1^2})A_{i1;1}^2 -K_1\mu\hfill\cr
\qquad\qquad\qquad\geqslant \sum_{i\in J}\frac{\mu_1}{(\lambda_1-\lambda_i)\lambda_1}A_{i1;1}^2 -K_1\mu\hfill\cr
\qquad\qquad\qquad\geqslant -K_1\mu.\hfill\cr
\endmatrix$$
\noindent Thus:
$$
-\frac{1}{\lambda_1}(D^2K)^{ij,mn}A_{ij;1}A_{mn;1} - \sum_{i\in J}\frac{\mu_i}{\lambda_1^2}A_{11;i}A_{11;i}\geqslant \frac{1}{8}\sum_{i\in J}\frac{\mu_i}{\lambda_1^2}A_{11;i}A_{11;i} - K_1\mu.\eqnum{\nexteqnno}
$$
\noindent For all $\xi$, $X$ and $Y$:
$$\matrix
\nabla^\Sigma\xi(Y;X) \hfill&= \nabla^M\xi(Y;X) - A(X,Y)\xi(N);\text{ and }\hfill\cr
X\xi(N) \hfill&= \nabla^M\xi(N;X) + \xi(A X).\hfill\cr
\endmatrix$$
\noindent Thus:
$$\matrix
R^M_{i1\nu 1;i} \hfill&= R^M_{i1\nu 1:i} + \lambda_i(1 - \delta_{i1}) R^M_{1 \nu\nu 1} + \lambda_i R^M_{i1i1},\hfill\cr
R^M_{i1\nu i;1} \hfill&= R^M_{i1\nu i:1} - \lambda_1(1 - \delta_{i1}) R^M_{i \nu\nu i} - \lambda_1 R^M_{i1i1}.\hfill\cr
\endmatrix$$
\noindent Bearing in mind that $\lambda_1\geqslant 1$, by Proposition \procref{PropCurvFnsII}, $(4)$, there exists $K_2$ in $\Cal{B}$ such that:
$$
\sum_{i=1}^n\frac{\mu_i}{\lambda_1}(R^M_{i1\nu 1;i}+R^M_{i1\nu i;1})\geqslant -K_2(1 + \mu).\eqnum{\nexteqnno}
$$
\noindent Next:
$$
{R^\Sigma_{1ii}}^pA_{p1} + {R^\Sigma_{1i1}}^pA_{pi} = R^M_{1ii1}(\lambda_1-\lambda_i) + \lambda_1\lambda_i(\lambda_1 - \lambda_i).
$$
\noindent Thus, bearing in mind that $\lambda_1\geqslant 1$ and that $\lambda_1\geqslant\lambda_i$ for all $i$, there exists $K_3$ in $\Cal{B}$ such that:
$$
\sum_{i,j=1}^n\frac{\mu_i}{\lambda_1}({R^\Sigma_{1ii}}^pA_{p1} + {R^\Sigma_{1i1}}^pA_{pi}) \geqslant -K_3(1 + \mu).\eqnum{\nexteqnno}
$$
\noindent Finally, bearing in mind Lemma \procref{PropHessianOfRestriction}:
$$\matrix
\kappa_{;11} \hfill&= \opHess^\Sigma(\kappa)(e_1,e_1)\hfill\cr
&= \opHess^M(\kappa)(e_1,e_1) - \langle\nabla\kappa,N\rangle A_{11}\hfill\cr
&= \opHess^M(\kappa)(e_1,e_1) - \lambda_1d\kappa(N).\hfill\cr
\endmatrix$$
\noindent Bearing in mind that $\lambda_1\geqslant 1$, there thus exists $K_4$ in $\Cal{B}$ such that:
$$
\frac{1}{\lambda_1}\kappa_{;11} \geqslant -K_4.\eqnum{\nexteqnno}
$$
\noindent Combining the above relations, there exists $K_5$ in $\Cal{B}$ such that:
$$
\Delta^K\opLog(a) \geqslant -K_5(1+\mu) - \sum_{i\in I}\frac{\mu_i}{\lambda_1^2}A_{11;i}^2 + \frac{1}{8}\sum_{i\in J}\frac{\mu_i}{\lambda_1^2}A_{11;i}^2.
$$
\noindent Finally, bearing in mind that $\lambda_1=a\geqslant 1$, there exists $K_6$ in $\Cal{B}(\phi,C)$ such that:
$$\matrix
C\Cal{D}_\phi\opLog(a)\hfill&= C\sum_{i=1}^n\frac{\mu_i}{\lambda_1}\phi_{;i}A_{11;i}\hfill\cr
&\geqslant -\frac{1}{8}\sum_{i=1}^n\frac{\mu_i}{\lambda_1^2}A_{11;i}^2 - K_6\mu.\hfill\cr
\endmatrix$$
\noindent The result now follows by combining the above relations.\qed
\medskip
\noindent We now prove Proposition \procref{CorSuperharmonic}:
\medskip
{\bf\noindent Proof of Proposition \procref{CorSuperharmonic}:\ }This follows immediately from Proposition \procref{PropSuperharmonic}.\qed
\newsubhead{Global second order a priori estimates - part I}Let $p$ be a point in $M$. Let $d_p$ be the distance in $M$ to $p$.
\proclaim{Proposition \nextprocno}
\noindent Let $R>0$ be such that $\kappa(q)<\mu_\infty(K)/R$ for all $q\in X\minter B_R(p)$. There exist $c,\epsilon>0$ in $\Cal{B}(R)$ such that if $\Sigma\subseteq X\minter B_R(p)$, then, over $\Sigma$:
$$
\lambda_1 \geqslant c\ \Rightarrow \Delta^K d_p^2 \geqslant \epsilon(1 + \mu).
$$
\endproclaim
\proclabel{PropLapOfDistanceFunctionAlt}
\proof Observe that, by homogeneity:
$$
\mliminf_{K(x)=1,x\rightarrow\infty}DK_x(\Bbb{I}) = \mu_\infty(K).
$$
\noindent Since $M$ has non-positive curvature, bearing in mind Proposition \procref{PropCurvFnsII} and Lemma \procref{PropHessianOfRestriction}:
$$\matrix
&\opHess^M(\frac{1}{2}d_p^2)\hfill&\geqslant \opId\hfill\cr
\Rightarrow\hfill&\opHess^\Sigma(\frac{1}{2} d_p^2)\hfill&\geqslant \opId - d_p\langle\msf{N},\nabla d_p\rangle A\hfill\cr
\Rightarrow\hfill&\Delta^K\frac{1}{2}d_p^2\hfill&\geqslant \mu - \kappa d_p\hfill\cr
& &\geqslant \mu - \kappa R.\hfill\cr
\endmatrix$$
\noindent Thus, by definition of $R$, there exists $c,\epsilon>0$ in $\Cal{B}(R)$ such that, for $\lambda_1\geqslant c$:
$$
\Delta^K d_p^2 \geqslant 2\epsilon\mu \geqslant \epsilon(1 + \mu),
$$
\noindent as desired.\qed
\proclaim{Proposition \nextprocno}
\noindent Let $R>0$ be such that $\kappa(q)<\mu_\infty(K)/R$ for all $q\in X\minter B_R(p)$. There exist $C,c>0$ in $\Cal{B}(R)$ and a homogeneous first order operator $\Cal{D}$ such that if $\Sigma\subseteq X\minter B_R(p)$ and if $\lambda_1(q)\geqslant c$, then, at $q$:
$$
(\Delta^K + \Cal{D})(\opLog(\lambda_1) + C\delta_p^2) > 0,
$$
\noindent in the weak sense.
\endproclaim
\proclabel{CorMaximumPrincipalAlt}
\proof Choose $q\in\Sigma$. Define $a$ near $q$ as in Section \subheadref{SubheadGlobalSecondOrderBoundsPreliminaryResults}. By Propositions \procref{PropSuperharmonic} and \procref{PropLapOfDistanceFunctionAlt}, there exists $c_1,C>0$ in $\Cal{B}(R)$ such that, for $\lambda_1>c_1$, at $q$:
$$
\Delta^K(\opLog(a) + Cd_p^2) \geqslant 1 - \sum_{i\in I}\frac{\mu_i}{\lambda_1^2}A_{11;i}^2.
$$
\noindent Denote $\Phi:=\opLog(a) + Cd_p^2$. For all $k$, at $q$:
$$
\frac{1}{\lambda_1}A_{11;k} = \opLog(a)_{;k} = -2Cd_p d_{p;k} + \Phi_{;k}.
$$
\noindent There therefore exist smooth functions $(D_i)_{i\in I}$ such that:
$$
\sum_{i\in I}\frac{\mu_i}{\lambda_1^2}A^2_{11;i} = 4C^2d_p^2\sum_{i\in I}\mu_i d_{p;k}^2 - \sum_{i\in I} D_i\Phi_{;k}. $$
\noindent By Proposition \procref{PropCurvFnsII}, $(4)$, $\lambda_1\mu_1\leqslant\kappa$. Thus $\mu_1\leqslant\kappa/\lambda_1$, and so $\mu_i\leqslant 4\kappa/\lambda_1$ for all $i\in I$. There therefore exists $c_2\in\Cal{B}(R)$ such that:
$$
\Delta^K\Phi + \sum_{i\in I}D_i\Phi_{;i} \geqslant 1 - c_2/\lambda_1.
$$
\noindent The result follows with $c=\opMax(c_1,c_2)$.\qed
\medskip
{\bf\noindent Proof of Proposition \procref{PropSecondOrderIntBoundsAlt}:\ }Let $C,c>0$ be as in Proposition \procref{CorMaximumPrincipalAlt}. Consider the function $\|A\|e^{Cd_p^2}=\lambda_1 e^{Cd_p^2}$. If this function acheives its maximum along $\partial\Sigma$, then the result follows since $e^{Cd_p^2}$ is uniformly bounded above and below. Otherwise, it acheives its maximum in the interior of $\Sigma$, in which case, by Proposition \procref{CorMaximumPrincipalAlt} and the maximum principle, at this point $\|A\|\leqslant\lambda_1\leqslant c$. The result follows.\qed
\newsubhead{Global second order a priori bounds - part II} We require the following variant of Lemma \procref{LemmaSuperHarmonicity}:
\proclaim{Proposition \nextprocno}
\noindent Let $\phi$ be as in Lemma \procref{LemmaSuperHarmonicity}. Suppose that $\phi\geqslant\delta>0$. There exists $\epsilon,\alpha,C$ in $\Cal{B}(\phi,\delta)$ such that:
$$
(\Delta^K + CD_\phi) \phi^{1+\alpha} \geqslant \epsilon(1 + \mu).
$$
\endproclaim
\proclabel{CorSuperHarmonicityII}
\proof Let $\pi$ be the orthogonal projection along $\nabla\phi$. Observe that for all $\alpha>0$ and for all $X\in T_pM$:
$$
\opHess(\phi^{1+\alpha})(X,X) \geqslant \opHess(\phi^{1+\alpha})(\pi(X),\pi(X)).
$$
\noindent Defining $m$ and $m'$ as in the proof of Lemma \procref{LemmaSuperHarmonicity}, we thus obtain:
$$
(m_{ij}) \geqslant \pmatrix (m'_{ij})\hfill& \opCos(\theta)(m'_{in})\hfill\cr \opCos(\theta)(m'_{ni})\hfill&\opCos^2(\theta)m'_{nn}\hfill\cr\endpmatrix.
$$
\noindent As in the proof of Corollary \procref{CorSuperHarmonicityI}, there exists $\epsilon_1>0$ in $\Cal{B}$ such that:
$$
\Delta^K\phi^{1+\alpha} \geqslant \epsilon_1\sum_{i=1}^n\mu_i\opHess(\phi^{1+\alpha})_{ii} - \|\opHess(\phi^{1+\alpha})\|\sum_{i=1}^n\mu_i\phi_{;i}\phi_{;i}.
$$
\noindent However, there exists $\epsilon>0$ in $\Cal{B}(\phi)$ such that:
$$
\epsilon_1\sum_{i=1}^n\mu_i\opHess(\phi^{1+\alpha})_{ii} \geqslant \epsilon\mu.
$$
\noindent Finally:
$$
D_\phi\phi^{1+\alpha} = (1+\alpha)\phi^\alpha\sum_{i=1}^n\mu_i\phi_{;i}\phi_{;i}.
$$
\noindent Thus, for $(1+\alpha)\delta^\alpha C \geqslant \|\opHess(\phi^{1+\alpha})\|$, the result follows.\qed
\medskip
\noindent Let $p$ be a point in $M$. Let $d_p$ be the distance in $M$ to $p$.
\proclaim{Proposition \nextprocno}
\noindent If $\kappa<1$ and if the sectional curvature of $M$ is bounded above by $-1$, then there exists $\epsilon,\alpha, C>0$ in $\Cal{B}$ such that, over $\Sigma$ and away from $p$:
$$
(\Delta^K + C\Cal{D}_{d_p})d_p^{1+\alpha} \geqslant \epsilon(1 + \mu).
$$
\endproclaim
\proclabel{PropLapOfDistanceFunction}
\proof Trivially, $\|\nabla d_p\|=1$. Moreover, since the level sets of $d_p$ are geodesic spheres and since $M$ is a Hadamard manifold, they are strictly convex. Since the sectional curvature of $M$ is bounded above by $-1$, by Properties $(3)$ and $(5)$ of $K$, the level sets have $K$-curvature greater than $1$. The result now follows by Proposition \procref{CorSuperHarmonicityII}.\qed
\proclaim{Proposition \nextprocno}
\noindent If $\kappa<1$ and if the section curvature of $M$ is bounded above by $-1$, then there exists $C,c>0$ in $\Cal{B}$ and a first order homogeneous operator $\Cal{D}$ such that over $\Sigma$:
$$
\lambda_1\geqslant c\ \Rightarrow (\Delta^K + \Cal{D})(\opLog(\lambda_1) + C\delta) > 0,
$$
\noindent in the weak sense.\qed
\endproclaim
\proclabel{CorMaximumPrincipal}
\proof By Propositions \procref{PropLapOfDistanceFunction} and \procref{CorSuperharmonic}, there exist $C_1,C_2>0$ and $\alpha\in]0,1[$ in $\Cal{B}$ such that, over $\Sigma$ and away from $p$:
$$
(\Delta + C\Cal{D}_{d_p})(\opLog(a) + C d_p^{1+\alpha}) \geqslant 1-\sum_{i\in I}\frac{\mu_i}{\lambda_1^2}A_{11;i}^2.
$$
\noindent The result now follows as in the proof of Proposition \procref{CorMaximumPrincipalAlt}.\qed
\medskip
{\bf\noindent Proof of Proposition \procref{PropSecondOrderIntBounds}:\ }This is identical to the proof of Proposition \procref{PropSecondOrderIntBoundsAlt}, with Proposition \procref{CorMaximumPrincipal} used instead of Proposition \procref{CorMaximumPrincipalAlt}.\qed
\newhead{Existence}%
\headlabel{HeadExistence}%
\newsubhead{Compactness} Let $K$ be a convex curvature function. Bearing in mind the two complementary results of Section \subheadref{SubheadGlobalSecondOrderEstimatesMainResults}, we obtain the following two complementary compactness results:
\subheadlabel{SubheadExistenceCompactness}
\proclaim{Proposition \nextprocno}
\noindent Let $K$ be a convex curvature function; let $M$ be an $(n+1)$-dimensional Hadamard manifold; let $\Cal{F}$ be a family of triplets $(\hat{\Sigma},\kappa,R)$ where $\hat{\Sigma}$ is a smooth compact LSC immersed hypersurface with generic boundary in $M$; $\kappa$ is a smooth positive function on $M$ such that $\kappa<\mu_\infty(K)/R$; $\hat{\Sigma}$ is contained in $B_R(p)$ for some $p\in M$; and $K(\hat{\Sigma})>\kappa$. Let $\Cal{G}$ be the family of smooth compact LSC immersed hypersurfaces $\Sigma$ in $M$ such that $K(\Sigma)=\kappa$ and $\Sigma<\hat{\Sigma}$ for some triplet $(\hat{\Sigma},\kappa,R)\in\Cal{F}$. If $\Cal{F}$ is compact, then so too is $\Cal{G}$.
\endproclaim
\proclabel{PropCompactnessI}
\proof By compactness, upon perturbing every element of $\Cal{F}$ if necessary, we may assume that there exists $\theta>0$ such that for all $\Sigma\in\Cal{G}$, there exists $(\hat{\Sigma},\kappa,R)\in\Cal{F}$ such that $K(\Sigma)=\kappa$; $\Sigma<\hat{\Sigma}$; and $\Sigma$ makes an angle of at least $\theta$ with $\hat{\Sigma}$ along their common boundary. Let $p\in M$ be such that $\hat{\Sigma}$ is contained in $B_R(p)$. We first show that $\Sigma$ is also contained in $B_R(p)$. Suppose the contrary. Let $q\in\Sigma$ be the point lying furthest from $p$. Let $L\subseteq M$ be the geodesic ray leaving $q$ in the outward-pointing normal direction from $\Sigma$. Trivially $L$ does not intersect $B_R(p)$. However, since $\hat{\Sigma}$ bounds $\Sigma$, $L$ meets $\hat{\Sigma}$ at some point. This is absurd, and it follows that $\Sigma$ is contained in $B_R(p)$ as asserted.
\medskip
\noindent Proposition \procref{PropSecondOrderBoundaryEstimates} yields a priori $C^2$ bounds for elements of $\Cal{G}$ along the boundary. By the preceeding paragraph, Proposition \procref{PropSecondOrderIntBoundsAlt} yields global a priori $C^2$-bounds for elements of $\Cal{G}$. Theorem $1$ of \cite{CaffNirSprII} yields a priori $C^{2,\alpha}$-bounds for elements of $\Cal{G}$. The Schauder estimates (c.f. \cite{GilbTrud}) yield a priori $C^k$-bounds for elements of $\Cal{G}$ for all $k\in\Bbb{N}$. Proposition \procref{PropDiameterBounds} yields a priori diameter bounds for elements of $\Cal{G}$. It follows from the Arzela-Ascoli Theorem for immersed hypersurfaces (c.f. \cite{SmiAAT}) that every sequence in $\Cal{G}$ contains a convergent subsequence.
\medskip
\noindent It remains to show that $\Cal{G}$ is closed. Let $\Sigma$ by a limit point of $\Cal{G}$. Let $(\Sigma_n)_\ninn$ be a sequence in $\Cal{G}$ converging to $\Sigma$. Let $(\hat{\Sigma}_n,\kappa_n,R_n)$ be a sequence in $\Cal{F}$ such that, for all $n$, $\Sigma_n<\hat{\Sigma}_n$ and $K(\Sigma_n)=\kappa_n$. Since $\Cal{F}$ is compact, we may suppose that $(\hat{\Sigma}_n,\kappa_n,R_n)$ converges to $(\hat{\Sigma},\kappa,R)$, say. Taking limits, $K(\Sigma)=\kappa$. Furthermore, by Proposition \procref{PropBoundednessIsAClosedProperty}, $\Sigma<\hat{\Sigma}$. It follows that $\Sigma$ is an element of $\Cal{G}$, and this completes the proof.\qed
\proclaim{Proposition \nextprocno}
\noindent Let $K$ be a convex curvature function; let $M$ be an $(n+1)$-dimensional Hadamard manifold with sectional curvature bounded above by $-1$; let $\Cal{F}$ be a family of pairs $(\hat{\Sigma},\kappa)$ where $\hat{\Sigma}$ is a smooth compact LSC immersed hypersurface with generic boundary in $M$; $\kappa:M\rightarrow]0,1[$ is a smooth positive function; and $K(\hat{\Sigma})>\kappa$. Let $\Cal{G}$ be the family of smooth compact LSC immersed hypersurfaces $\Sigma$ in $M$ such that $\Sigma<\hat{\Sigma}$ and $K(\Sigma)=\kappa$ for some pair $(\hat{\Sigma},\kappa)\in\Cal{F}$. If $\Cal{F}$ is compact, then so too is $\Cal{G}$.
\endproclaim
\proclabel{PropCompactnessII}
\proof This is proven in the same way as Proposition \procref{PropCompactnessI}, using Proposition \procref{PropSecondOrderIntBounds} instead of Proposition \procref{PropSecondOrderIntBoundsAlt}.\qed
\newsubhead{Proof of Main Results}We first show:
\proclaim{Proposition \nextprocno}
\noindent Let $M$ be a Hadamard manifold. Let $\Sigma$ be a smooth compact LSC immersed hypersurface in $M$. $\Sigma$ is isotopic through smooth compact LSC immersed hypersurfaces to an immersed submanifold lying on a geodesic sphere.
\endproclaim
\proclabel{PropIsotopies}
\proof Let $N_{\Sigma}$ be the outward pointing, unit, normal vector field over $\Sigma$ and define $I:\Sigma\times[0,\infty[\rightarrow M$ by:
$$
I(p,t) = \opExp(t\msf{N}_{\Sigma}(p)).
$$
\noindent Let $d:\Sigma\times[0,\infty[\rightarrow\Bbb{R}$ be the distance in $\Sigma\times[0,\infty[$ to $\Sigma\times\left\{0\right\}$. $d$ is a locally convex function. For all $t\in[0,\infty[$, let $\Sigma_t=d^{-1}(\left\{t\right\})$ be the level set of $d$ at height $t$. Choose $p\in M$ and let $d_p$ be the distance to $p$ in $M$. $d_p$ is also a convex function and we identify it with $d_p\circ I$.
\medskip
\noindent Observe that $\partial(\Sigma\times[0,\infty[)$ consists of $2$ components, being $\Sigma$ and $(\partial\Sigma)\times[0,\infty[$. Choose $R\geqslant 0$ such that, for all $q\in\Sigma$:
$$
d_p(q)<R.
$$
\noindent Observe that, since $M$ is non-positively curved, as $d(q)$ tends to $+\infty$, the angle between $\nabla d$ and $\nabla d_p$ at $q$ tends to $0$. Thus, increasing $R$ if necessary, we may assume that, for $d_p(q)\geqslant R$:
$$
\langle \nabla d,\nabla d_p\rangle(q) > 0.
$$
\noindent For $s\in [0,1]$ define $d_s$ and $\tilde{\Sigma}_s$ by:
$$
d_s = sd_p + (1-s)d,\qquad \tilde{\Sigma}_s = d_s^{-1}(\left\{R\right\}).
$$
\noindent For all $s$, and for all $q$ such that $d_p(q)\geqslant R$:
$$
\langle\nabla d,\nabla d_s\rangle(q) >0.
$$
\noindent Thus, for all $s\in[0,1]$, $\tilde{\Sigma}_s\minter\Sigma_0=\emptyset$ and $\tilde{\Sigma}_s$ is transverse to $\partial\Sigma\times[0,\infty[$. Moreover, for all $s$, $d_s$ is convex, and so $(\Sigma_s)_{s\in[0,1]}$ defines an isotopy through smooth compact LSC immersed hypersurfaces from $\Sigma_R=\tilde{\Sigma}_0$ to $\tilde{\Sigma}_1\subseteq\partial B_R(p)$. Since $\Sigma_0$ is isotopic to $\Sigma_R$, the result follows.\qed
\medskip
\noindent This allows us to prove Theorems \procref{ThmFirstExistenceTheorem} and \procref{ThmSecondExistenceTheorem}:
\medskip
{\bf\noindent Proof of Theorems \procref{ThmFirstExistenceTheorem} and \procref{ThmSecondExistenceTheorem}:\ }By Proposition \procref{PropBoundednessIsAnOpenProperty}, boundedness is an open property and so differential topological techniques may be applied. The compactness result of Propositions \procref{PropCompactnessI} and \procref{PropCompactnessII} together with the isotopy result of Proposition \procref{PropIsotopies} then allow us to apply the differential topological degree theory developed for the case of extrinsic curvature in \cite{SmiPPH}. Existence follows.\qed
\medskip
\noindent We prove Theorem \procref{ThmDegreeUniqueness}:
\medskip
{\noindent\bf Proof of Theorem \procref{ThmDegreeUniqueness}:\ }Let $\kappa:M\rightarrow]0,1[$ be a smooth function. Let $\Sigma=(i,(S,\partial S))$ be a smooth compact LSC immersed hypersurface $K(\Sigma)=\kappa$. Let $N$ be the upward- pointing unit normal vector field over $\Sigma$. Let $J$ be the Jacobi operator of $K$-curvature over $\Sigma$. As in  Proposition $3.1.1$ of \cite{Lab}, for all $f\in C^\infty(S)$:
$$
\Cal{L}f = (DK_A(W) - DK_A(A^2))f - DK_A(\opHess(f)),
$$
\noindent where $W:TS\rightarrow TS$ is given by:
$$
W\cdot X = R_{\msf{N}X}\msf{N}.
$$
\noindent Since the sectional curvature of $M$ is bounded above by $-1$, for all $X\in TS$:
$$\matrix
&\langle W\cdot X,X\rangle \hfill&= \langle R_{\msf{N}X}\msf{N},X\rangle\hfill\cr
& &\geqslant \|X\|^2\hfill\cr
\Rightarrow\hfill&W\hfill&\geqslant \opId\hfill\cr
\Rightarrow\hfill&DK_A(W)\hfill&\geqslant DK_A(\opId).\hfill\cr
\endmatrix$$
\noindent It thus follows from the hypotheses on $K$ that:
$$
DK_A(W) - DK_A(A^2) > 0.
$$
\noindent Thus, for all $f\in C^\infty(S)$:
$$
\langle\Cal{L}f,f\rangle \geqslant 0,
$$
\noindent with equality if and only if $f=0$.
\medskip
\noindent We may conclude in one of two different ways. First, we may interpret this in terms of the degree theory of \cite{SmiDTI}, in which case we see that the contribution of any solution to the degree of $\kappa$ is equal to $+1$, and that there is therefore only one solution. Alternatively, we reason more directly, as in the proof of Lemma $3.0.2$ of \cite{Lab}, to reach the same conclusion. This completes the proof.\qed
\goodbreak
\newhead{Bibliography}
{\leftskip = 5ex \parindent = -5ex
\leavevmode\hbox to 4ex{\hfil \cite{Caffarelli}}\hskip 1ex{Caffarelli L., Monge-Amp\`ere equation, div-curl theorems in Lagrangian coordinates, compression and rotation, Lecture Notes, 1997}
\medskip
\leavevmode\hbox to 4ex{\hfil \cite{CaffNirSprI}}\hskip 1ex{Caffarelli L., Nirenberg L., Spruck J., The Dirichlet problem for nonlinear second-order elliptic equations. I. Monge Amp\`ere equation, {\sl Comm. Pure Appl. Math.} {\bf 37} (1984), no. 3, 369--402}
\medskip
\leavevmode\hbox to 4ex{\hfil \cite{CaffNirSprII}}\hskip 1ex{Caffarelli L., Kohn J. J., Nirenberg L., Spruck J., The Dirichlet problem for nonlinear second-order elliptic equations. II. Complex Monge Amp\`ere, and uniformly elliptic, equations, {\sl Comm. Pure Appl. Math.} {\bf 38} (1985), no. 2, 209--252}
\medskip
\leavevmode\hbox to 4ex{\hfil \cite{CaffNirSprIII}}\hskip 1ex{Caffarelli L., Nirenberg L., Spruck J., The Dirichlet problem for nonlinear second-order elliptic equations. III. Functions of the eigenvalues of the Hessian, {\sl Acta Math.} {\bf 155} (1985), no. 3-4, 261--301}
\medskip
\leavevmode\hbox to 4ex{\hfil \cite{CaffNirSprV}}\hskip 1ex{Caffarelli L., Nirenberg L., Spruck J., Nonlinear second-order elliptic equations. V. The Dirichlet problem for Weingarten hypersurfaces, {\sl Comm. Pure Appl. Math.} {\bf 41} (1988), no. 1, 47--70}
\medskip
\leavevmode\hbox to 4ex{\hfil \cite{CalI}}\hskip 1ex{Calabi E., Improper affine hyperspheres of convex type and a generalization of a theorem by K. J\"orgens, {\sl Michigan Math. J.} {\bf 5} (1958), 105--126}
\medskip
\leavevmode\hbox to 4ex{\hfil \cite{CrandallIshiiLions}}\hskip 1ex{Crandall M. G., Ishii H., Lions P. L., User's guide to viscosity solutions of second order partial differential equations, {\sl Bull. Amer. Math. Soc.} {\bf 27} (1992), no. 1, 1--67}
\medskip
\leavevmode\hbox to 4ex{\hfil \cite{GilbTrud}}\hskip 1ex{Gilbarg D., Trudinger N. S., {\sl Elliptic partial differential equations of second order}, Classics in Mathematics. Springer-Verlag, Berlin, 2001}
\medskip
\leavevmode\hbox to 4ex{\hfil \cite{GuanSpruckA}}\hskip 1ex{Guan B., Spruck J., Boundary value problems on $S^n$ for surfaces of constant Gauss curvature, {\sl Ann. of Math.} {\bf 138} (1993), 601--624}
\medskip
\leavevmode\hbox to 4ex{\hfil \cite{Guan}}\hskip 1ex{Guan B., The Dirichlet problem for Monge-Amp\`ere equations in non-convex domains and spacelike hypersurfaces of constant Gauss curvature, {\sl Trans. Amer. Math. Soc.} {\bf 350} (1998), 4955--4971}
\medskip
\leavevmode\hbox to 4ex{\hfil \cite{GuanSpruck}}\hskip 1ex{Guan B., Spruck J., The existence of hypersurfaces of constant Gauss curvature with prescribed boundary, {\sl J. Differential Geom.} {\bf 62} (2002), no. 2, 259--287}%
\medskip
\leavevmode\hbox to 4ex{\hfil \cite{GuanSpruckI}}\hskip 1ex{Guan B., Spruck J., Locally convex hypersurfaces of constant curvature with boundary, {\sl Comm. Pure Appl. Math.} {\bf 57} (2004), no. 10, 1311--1331}%
\medskip
\leavevmode\hbox to 4ex{\hfil \cite{GuanSpruckII}}\hskip 1ex{Guan B., Spruck J., Szapiel M., Hypersurfaces of constant curvature in hyperbolic space. I., {\sl J. Geom. Anal.} {\bf 19} (2009), no. 4, 772--795}%
\medskip
\leavevmode\hbox to 4ex{\hfil \cite{GuanSpruckIII}}\hskip 1ex{Guan B., Spruck J., Hypersurfaces of constant curvature in hyperbolic space. II., {\sl J. Eur. Math. Soc. (JEMS)} {\bf 12} (2010), no. 3, 797--817}
\medskip
\leavevmode\hbox to 4ex{\hfil \cite{GuanSpruckIV}}\hskip 1ex{Guan B., Spruck J., Convex hypersurfaces of constant curvature in hyperbolic space, {\sl in preparation}}
\medskip
\leavevmode\hbox to 4ex{\hfil \cite{IvochTomi}}\hskip 1ex{Ivochkina N. M., Tomi F., Locally convex hypersurfaces of prescribed curvature and boundary, {\sl Calc. Var. Partial Differential Equations} {\bf 7} (1998), no. 4, 293--314}
\medskip
\leavevmode\hbox to 4ex{\hfil \cite{Jorgens}}\hskip 1ex{J\"orgens K., \"Uber die L\"osungen der Differentialgleichung $rt-s^2=1$  (German), {\sl Math. Ann.} {\bf 127} (1954), 130--134}%
\medskip
\leavevmode\hbox to 4ex{\hfil \cite{Lab}}\hskip 1ex{Labourie F., Un lemme de Morse pour les surfaces convexes (French), {\sl Invent. Math.} {\bf 141} (2000), no. 2, 239--297}
\medskip
\leavevmode\hbox to 4ex{\hfil \cite{Pog}}\hskip 1ex{Pogorelov A. V., On the improper convex affine hyperspheres, {\sl Geometriae Dedicata} {\bf 1} (1972), no. 1, 33--46}
\medskip
\leavevmode\hbox to 4ex{\hfil \cite{SmiDTI}}\hskip 1ex{Rosenberg H., Smith G., Degree Theory of Immersed Hypersurfaces, arXiv:1010.1879}
\medskip
\leavevmode\hbox to 4ex{\hfil \cite{RosSpruck}}\hskip 1ex{Rosenberg H., Spruck J., On the existence of convex hypersurfaces of constant Gauss curvature in hyperbolic space, {\sl J. Differential Geom.} {\bf 40} (1994), no. 2, 379--409}
\medskip
\leavevmode\hbox to 4ex{\hfil \cite{ShengUrbasWang}}\hskip 1ex{Sheng W., Urbas J., Wang X., Interior curvature bounds for a class of curvature equations. (English summary), {\sl Duke Math. J.} {\bf 123} (2004), no. 2, 235--264}
\medskip
\leavevmode\hbox to 4ex{\hfil \cite{Smale}}\hskip 1ex{Smale S., An infinite dimensional version of Sard's theorem, {\sl Amer. J. Math.} {\bf 87} (1965), 861--866}
\medskip
\leavevmode\hbox to 4ex{\hfil \cite{SmiAAT}}\hskip 1ex{Smith G., An Arzela-Ascoli Theorem for Immersed Submanifolds, {\sl Ann. Fac. Sci. Toulouse Math.} {\bf 16} (2007), no. 4, 817--866}
\medskip
\leavevmode\hbox to 4ex{\hfil \cite{SmiCGC}}\hskip 1ex{Smith G., Compactness results for immersions of prescribed Gaussian curvature I - analytic aspects, to appear in {\sl Adv. Math.}}
\medskip
\leavevmode\hbox to 4ex{\hfil \cite{SmiSLC}}\hskip 1ex{Smith G., Special Lagrangian curvature, to appear in {\sl Math. Ann.}}
\medskip
\leavevmode\hbox to 4ex{\hfil \cite{SmiNLD}}\hskip 1ex{Smith G., The non-linear Dirichlet problem in Hadamard manifolds, arXiv:0908.3590}
\medskip
\leavevmode\hbox to 4ex{\hfil \cite{SmiFCS}}\hskip 1ex{Smith G., Moduli of Flat Conformal Structures of Hyperbolic Type, {\sl Geom. Dedicata} {\bf 154} (2011), no. 1, 47--80}
\medskip
\leavevmode\hbox to 4ex{\hfil \cite{SmiPPH}}\hskip 1ex{Smith G., Compactness results for immersions of prescribed Gaussian curvature II - geometric aspects, to appear in {\sl Geom. Dedicata}}
\medskip
\leavevmode\hbox to 4ex{\hfil \cite{SmiNLP}}\hskip 1ex{Smith G., The non-linear Plateau problem in non-positively curved manifolds, to appear in {\sl Trans. Amer. Math. Soc}}
\medskip
\leavevmode\hbox to 4ex{\hfil \cite{SmiPPG}}\hskip 1ex{Smith G., {\sl The Plateau problem for Gaussian curvature}, arXiv:1206.5544}
\medskip
\leavevmode\hbox to 4ex{\hfil \cite{Spruck}}\hskip 1ex{Spruck J., Fully nonlinear elliptic equations and applications to geometry, {\sl Proceedings of the International Congress of Mathematicians}, (Z\"urich, 1994), 1145--1152, Birkh\"a user, Basel, 1995.}
\medskip
\leavevmode\hbox to 4ex{\hfil \cite{Trud}}\hskip 1ex{Trudinger N. S., On the Dirichlet problem for Hessian equations, {\sl Acta Math.}, {\bf 175}, (1995), 151--164}
\medskip
\leavevmode\hbox to 4ex{\hfil \cite{TrudWang}}\hskip 1ex{Trudinger N. S., Wang X., On locally locally convex hypersurfaces with boundary, {\sl J. Reine Angew. Math.} {\bf 551} (2002), 11--32}
\medskip
\leavevmode\hbox to 4ex{\hfil \cite{WhiteI}}\hskip 1ex{White B., The space of $m$-dimensional hypersurfaces that are stationary for a parametric elliptic functional, {\sl Indiana Univ. Math. J.}, {\bf 36}, (1987), no. 3, 567--602}
\medskip
\leavevmode\hbox to 4ex{\hfil \cite{WhiteII}}\hskip 1ex{White B., The space of minimal submanifolds for varying Riemannian metrics, {\sl Indiana Univ. Math. J.}, {\bf 40}, (1991), no. 1, 161--200}
\par}
\enddocument
%
%%%%%%%%%%%%%%%%%%%%%%%%%%%%%%%%%%%%%%%%%%%%%%%%%%%%%%%%%%%%%%%%%%%%%%%%%%%%%%%%%%%%%%%%%%%%%%%%%%%%%%%%%%%%%%%%%%%%%%%
%
% 3: Closing commands.
%
%%%%%%%%%%%%%%%%%%%%%%%%%%%%%%%%%%%%%%%%%%%%%%%%%%%%%%%%%%%%%%%%%%%%%%%%%%%%%%%%%%%%%%%%%%%%%%%%%%%%%%%%%%%%%%%%%%%%%%%
%
\enddocument

%% file: preamble.tex
%%%%%%%%%%%%%%%%%%%%%%%%%%%%%%%%%%%%%%%%%%%%%%%%%%%%%%%%%%%%%%%%%%%%%%%%%%%%%%%%%%%%%%%%%%%%%%%%%%%%%%%%%%%%%%%%%%%%%%%
%
% 0: Preliminaries.
%
%%%%%%%%%%%%%%%%%%%%%%%%%%%%%%%%%%%%%%%%%%%%%%%%%%%%%%%%%%%%%%%%%%%%%%%%%%%%%%%%%%%%%%%%%%%%%%%%%%%%%%%%%%%%%%%%%%%%%%%
%
%
\let\myfrac=\frac%
\input eplain %
\let\frac=\myfrac%
\let\myfootnote=\footnote%
\input amstex \input epsf %
\let\footnote=\myfootnote%
%
% Here we load the functions permitting us to use "amsmath" without using "amsppt".
%
\loadeufm\loadmsam\loadmsbm\message{symbol names}\UseAMSsymbols\message{,}%
\magnification 1200 %
\font\myfontdefault=cmr10%
\newif\ifmakebiblio%
\newif\ifinappendices%
\newif\ifundefinedreferences%
\newif\ifchangedreferences%
\makebibliofalse%
\undefinedreferencesfalse%
\changedreferencesfalse%
%
%%%%%%%%%%%%%%%%%%%%%%%%%%%%%%%%%%%%%%%%%%%%%%%%%%%%%%%%%%%%%%%%%%%%%%%%%%%%%%%%%%%%%%%%%%%%%%%%%%%%%%%%%%%%%%%%%%%%%%%
%
% 1: Abstract machinery.
%
% Here we define the macro "makecounter", which gives functionality to the counters defined presently. In order
% to implement it, it is necessary to provisionally change the catcodes.
%
%%%%%%%%%%%%%%%%%%%%%%%%%%%%%%%%%%%%%%%%%%%%%%%%%%%%%%%%%%%%%%%%%%%%%%%%%%%%%%%%%%%%%%%%%%%%%%%%%%%%%%%%%%%%%%%%%%%%%%%
%
\def\setcatcodes{\catcode`\!=0 \catcode`\\=11}%
{\global\let\noe=\noexpand%
\catcode`\@=11 \catcode`\_=11 \setcatcodes%
!global!def!_@@internal@@makeref#1{%
!global!expandafter!def!csname #1ref!endcsname##1{%
!csname _@#1@##1!endcsname%
!expandafter!ifx!csname _@#1@##1!endcsname!relax%
    !write16{#1 ##1 not defined - run saving references}%
    !undefinedreferencestrue%
!fi}}%
!global!def!_@@internal@@makelabel#1{%
!global!expandafter!def!csname #1label!endcsname##1{%
!edef!temptoken{!csname #1info!endcsname}%
!ifloadreferences%
    !expandafter!ifx!csname _@#1@##1!endcsname!relax%
        !write16{#1 ##1 not hitherto defined - rerun saving references}%
        !changedreferencestrue%
    !else%
        !expandafter!ifx!csname _@#1@##1!endcsname!temptoken%
        !else%
            !write16{#1 ##1 reference has changed - rerun saving references}%
            !changedreferencestrue%
        !fi%
    !fi%
!else%
    !expandafter!edef!csname _@#1@##1!endcsname{!temptoken}%
    !edef!textoutput{!write!references{\global\def\_@#1@##1{!temptoken}}}%
    !textoutput%
!fi}}%
!global!def!makecounter#1{!_@@internal@@makelabel{#1}!_@@internal@@makeref{#1}}%
!unsetcatcodes%
}
%
%%%%%%%%%%%%%%%%%%%%%%%%%%%%%%%%%%%%%%%%%%%%%%%%%%%%%%%%%%%%%%%%%%%%%%%%%%%%%%%%%%%%%%%%%%%%%%%%%%%%%%%%%%%%%%%%%%%%%%%
%
% 2: Counters.
%
% Here we define the various counters.
%
%%%%%%%%%%%%%%%%%%%%%%%%%%%%%%%%%%%%%%%%%%%%%%%%%%%%%%%%%%%%%%%%%%%%%%%%%%%%%%%%%%%%%%%%%%%%%%%%%%%%%%%%%%%%%%%%%%%%%%%
%
\def\turnintolatin#1{\ifcase #1 _\or i\or ii\or iii\or iv\or v\or vi\or vii\or viii\or ix\or x\or xi\or xii\or xiii\or xiv\or xv\or xvi\or xvii\or xviii\or xix\or xx\or xxi\or xxii\or xxiii\or xxiv\or xxv\or xxvi\fi}%
\def\alphanum#1{\ifcase #1 _\or A\or B\or C\or D\or E\or F\or G\or H\or I\or J\or K\or L\or M\or N\or O\or P\or Q\or R\or S\or T\or U\or V\or W\or X\or Y\or Z\fi}%
\newwrite\references%
\ifloadreferences{\catcode`\@=11 \catcode`\_=11 \input references.tex }%
\else{\openout\references=references.tex }%
\fi%
%
% A: Headings.
%
\newcount\headno%
\global\headno=0%
\def\headinfo{\ifinappendices\alphanum\headno\else\the\headno\fi}%
\def\nextheadno{\global\advance\headno by 1 \global\subheadno=0 \global\procno=0 \headinfo}%
\makecounter{head}%
%
% B: Subheadings.
%
\newcount\subheadno%
\global\subheadno=0%
\def\subheadinfo{\headinfo.\the\subheadno}%
\def\nextsubheadno{\global\advance\subheadno by 1 \global\procno=0 \subheadinfo}%
\makecounter{subhead}%
%
% C: Proclaims (Theorems, Propositions, Lemmas, Corollories, Definitions).
%
\newcount\procno%
\global\procno=0%
\def\procinfo{\subheadinfo.\the\procno}%
\def\nextprocno{\global\advance\procno by 1 \procinfo}%
\makecounter{proc}%
%
% D: Figures.
%
\newcount\figno%
\global\figno=0%
\def\figinfo{\subheadinfo.\the\figno}%
\def\nextfigno{\global\advance\figno by 1 \figinfo}%
\makecounter{fig}%
%
% E: Equations.
%
\newcount\eqnno%
\global\eqnno=0%
\def\eqninfo{\text{(\the\eqnno)}}%
\def\nexteqnno{\global\advance\eqnno by 1 \eqninfo}%
\makecounter{eqn}%
%
%%%%%%%%%%%%%%%%%%%%%%%%%%%%%%%%%%%%%%%%%%%%%%%%%%%%%%%%%%%%%%%%%%%%%%%%%%%%%%%%%%%%%%%%%%%%%%%%%%%%%%%%%%%%%%%%%%%%%%%
%
% 3: Citations.
%
% Citations are treated as a special type of counter.
%
%%%%%%%%%%%%%%%%%%%%%%%%%%%%%%%%%%%%%%%%%%%%%%%%%%%%%%%%%%%%%%%%%%%%%%%%%%%%%%%%%%%%%%%%%%%%%%%%%%%%%%%%%%%%%%%%%%%%%%%
%
\def\gobbleeight#1#2#3#4#5#6#7#8{}%
\newcount\citationno%
\global\citationno=0%
\def\citationinfo{\the\citationno}%
\makecounter{citation}%
\newwrite\biblio%
\def\newref#1#2{%
\def\temptext{#2}%
\edef\bibliotextoutput{\expandafter\gobbleeight\meaning\temptext}%
\global\advance\citationno by 1\citationlabel{#1}%
\ifmakebiblio%
    \edef\fileoutput{\write\biblio{\noindent\hbox to 0pt{\hss$[\the\citationno]$}\hskip 0.2em\bibliotextoutput\medskip}}%
    \fileoutput%
\fi}%
\def\cite#1{%
$[\citationref{#1}]$%
\ifmakebiblio%
    \edef\fileoutput{\write\biblio{#1}}%
    \fileoutput%
\fi%
}%
%
%%%%%%%%%%%%%%%%%%%%%%%%%%%%%%%%%%%%%%%%%%%%%%%%%%%%%%%%%%%%%%%%%%%%%%%%%%%%%%%%%%%%%%%%%%%%%%%%%%%%%%%%%%%%%%%%%%%%%%%
%
% 4: Formatting.
%
%%%%%%%%%%%%%%%%%%%%%%%%%%%%%%%%%%%%%%%%%%%%%%%%%%%%%%%%%%%%%%%%%%%%%%%%%%%%%%%%%%%%%%%%%%%%%%%%%%%%%%%%%%%%%%%%%%%%%%%
%
\let\mypar=\par%
\edef\Pagetitle={Blank}\headline={\hfil\Pagetitle\hfil}%
\edef\Pagefooter={Blank}\footline={\hfil\Pagefooter\hfil}%
%
% A: Move to next odd page.
%
\newcount\showpagenumflag%
\global\showpagenumflag=0 %
\def\nextoddpage%
{\newpage\ifodd\pageno%
\else\global\showpagenumflag=0 %
\null\vfil\eject%
\global\showpagenumflag=1 %
\fi}%
%
% B: Headings.
%
\font\headfont=cmb12%
\def\newhead#1%
{\ifhmode\mypar\fi%
\ifnum\headno=0 \else\goodbreak\bigskip\fi%
{\headfont\noindent\nextheadno\ - #1.}
\nobreak\medskip}%
%
% C: Subheadings.
%
\def\newsubhead#1%
{\ifhmode\mypar\fi%
\ifnum\subheadno=0 \else\goodbreak\medskip\fi%
{\bf\noindent\nextsubheadno\ - #1.\ }}%
%
% D: Proclaims.
%
\newif\ifinproclaim%
\global\inproclaimfalse%
\def\proclaim#1{%
\goodbreak\medskip
\bgroup\inproclaimtrue%
\noindent{\bf #1}%
\nobreak\medskip\sl}%
\def\noskipproclaim#1{%
\goodbreak\medskip%
\bgroup\inproclaimtrue%
\noindent{\bf #1}\nobreak\sl}%
\def\endproclaim{\mypar\egroup\nobreak\medskip\ignorespaces}%
%
% E: Figures.
%
% The macro "\makelabelgrid" is a useful utility for guiding the positioning of labels is figures.
%
\newcount\xpos\newcount\ypos
\def\makelabelgrid{%
\xpos=-5 \ypos=-5 %
\loop\ifnum\xpos<6 %
{\loop\ifnum\ypos<6 %
\def\labeltext{x}%
\ifnum\xpos=0\def\labeltext{+}\fi%
\ifnum\ypos=0\def\labeltext{+}\fi%
\placelabel[\xpos][\ypos]{\labeltext}%
\advance\ypos by 1 %
\repeat}%
\advance\xpos by 1 %
\repeat}%
\def\placelabel[#1][#2]#3{{%
\setbox10=\hbox{\raise #2cm \hbox{\hskip #1cm #3}}%
\ht10=0pt \dp10=0pt \wd10=0pt \box10}}%
%
%
%
% F: Items.
%
\def\myitem#1{\noindent\hbox to .5cm{\hfill#1\hss}}%
%
% G: Right justification.
%
%
%
%%%%%%%%%%%%%%%%%%%%%%%%%%%%%%%%%%%%%%%%%%%%%%%%%%%%%%%%%%%%%%%%%%%%%%%%%%%%%%%%%%%%%%%%%%%%%%%%%%%%%%%%%%%%%%%%%%%%%%%
%
% 5: Special fonts.
%
%%%%%%%%%%%%%%%%%%%%%%%%%%%%%%%%%%%%%%%%%%%%%%%%%%%%%%%%%%%%%%%%%%%%%%%%%%%%%%%%%%%%%%%%%%%%%%%%%%%%%%%%%%%%%%%%%%%%%%%
%
% A: The "mathsf" font is not defined in Plain.
%
%
\font\sansseriften=cmss10%
\font\sansserifseven=cmss7%
\font\sansseriffive=cmss5%
\newfam\sansseriffam%
\textfont\sansseriffam=\sansseriften%
\scriptfont\sansseriffam=\sansserifseven%
\scriptscriptfont\sansseriffam=\sansseriffive%
\def\mathsf{\fam\sansseriffam}%
\def\msf#1{{\mathsf #1}}%
%
% B: The "mathbf" font is not defined in Plain.
%
\font\boldten=cmb10%
\font\boldseven=cmb7%
\font\boldfive=cmb5%
\newfam\mathboldfam%
\textfont\mathboldfam=\boldten%
\scriptfont\mathboldfam=\boldseven%
\scriptscriptfont\mathboldfam=\boldfive%
\def\mathbf{\fam\mathboldfam}%
%
%
% C: Here we define the macros "mathi" and "mathj". This is not really necessary, since the macros "imath" and
% "jmath" perform the same function. Still, it makes for an interesting exercise.
%
\font\mycmmiten=cmmi10%
\font\mycmmiseven=cmmi7%
\font\mycmmifive=cmmi5%
\newfam\mycmmifam%
\textfont\mycmmifam=\mycmmiten%
\scriptfont\mycmmifam=\mycmmiseven%
\scriptscriptfont\mycmmifam=\mycmmifive%
\def\hexa#1{\ifcase #1 0\or 1\or 2\or 3\or 4\or 5\or 6\or 7\or 8\or 9\or A\or B\or C\or D\or E\or F\fi}%
\mathchardef\mathi="7\hexa\mycmmifam7B%
\mathchardef\mathj="7\hexa\mycmmifam7C%
%
% D: Here we define a few Hebrew letters: "mybeth", "mygimmel" and "mydaleth".
%
\font\mymsbmten=msbm10 at 8pt%
\font\mymsbmseven=msbm7 at 5.6pt%6
\font\mymsbmfive=msbm5 at 4pt%
\newfam\mymsbmfam%
\textfont\mymsbmfam=\mymsbmten%
\scriptfont\mymsbmfam=\mymsbmseven%
\scriptscriptfont\mymsbmfam=\mymsbmfive%
\mathchardef\mybeth="7\hexa\mymsbmfam69%
\mathchardef\mygimmel="7\hexa\mymsbmfam6A%
\mathchardef\mydaleth="7\hexa\mymsbmfam6B%
%
%%%%%%%%%%%%%%%%%%%%%%%%%%%%%%%%%%%%%%%%%%%%%%%%%%%%%%%%%%%%%%%%%%%%%%%%%%%%%%%%%%%%%%%%%%%%%%%%%%%%%%%%%%%%%%%%%%%%%%%
%
% 6: Mathematics operators and symbols.
%
%%%%%%%%%%%%%%%%%%%%%%%%%%%%%%%%%%%%%%%%%%%%%%%%%%%%%%%%%%%%%%%%%%%%%%%%%%%%%%%%%%%%%%%%%%%%%%%%%%%%%%%%%%%%%%%%%%%%%%%
%
\def\proof{{\noindent\bf Proof:\ }}%
\def\remark{{\noindent\bf Remark:\ }}%
\def\qed{~$\square$}%
\def\makeop#1{\global\expandafter\def\csname op#1\endcsname{{\text{#1}}}}%
\def\makeopsmall#1{\global\expandafter\def\csname op#1\endcsname{{\text{\lowercase{#1}}}}}%
%
% A: Set Theory.
%
\def\munion{\mathop{\cup}}%
\def\minter{\mathop{\cap}}%
%
% B: Point set topology.
%
\makeop{Ext}%
\makeop{Int}%
\makeop{Dist}%
\makeop{Diam}%
\makeop{Length}%
%
% C: Sequences.
%
%
\def\ninn{{n\in\Bbb{N}}}%
\def\mlim{\mathop{{\text{Lim}}}}%
\def\mliminf{\mathop{{\text{LimInf}}}}%
\def\msup{\mathop{{\text{Sup}}}}%
\def\minf{\mathop{{\text{Inf}}}}%
%
% D: Linear Algebra.
%
\makeop{Dim}%
\makeop{Ker}%
\makeop{Coker}%
\makeop{Tr}%
\makeop{Adj}%
\makeop{Det}%
\makeop{End}%
\makeop{Lin}%
\makeop{Symm}%
\makeop{Mult}%
%
% E: Basic calculus.
%
\makeop{dx}%
\makeop{dy}%
\makeop{dz}%
\makeop{dt}%
\makeop{dVol}%
\makeop{dArea}%
\makeop{Supp}%
\makeop{Hess}%
\makeop{Lip}%
%
% F: Complex Numbers.
%
\makeop{Re}%
\makeop{Im}%
\makeop{Arg}%
\makeop{Log}%
\makeop{Exp}%
%
% G: Trigonometry.
%
\makeopsmall{Cos}%
\makeopsmall{Sin}%
\makeopsmall{Tan}%
\makeopsmall{Sec}%
\makeopsmall{Cosec}%
\makeopsmall{Cot}%
\makeopsmall{ArcCos}%
\makeopsmall{ArcSin}%
\makeopsmall{ArcTan}%
\makeopsmall{ArcSec}%
\makeopsmall{ArcCosec}%
\makeopsmall{ArcCot}%
%
% H: Hyperbolic Trigonometry.
%
\makeopsmall{Cosh}%
\makeopsmall{Sinh}%
\makeopsmall{Tanh}%
\makeopsmall{ArcCosh}%
\makeopsmall{ArcSinh}%
\makeopsmall{ArcTanh}%
%
% I: Differential and Riemannian Geometry.
%
\makeop{Vol}%
\makeop{Area}%
\makeop{Riem}%
\makeop{Ric}%
\makeop{Scal}%
\makeop{Euc}%
\makeop{Imm}%
\makeop{Emb}%
%
% J: Lie Groups.
%
\makeop{Id}%
\makeop{Ad}%
\makeop{O}%
\makeop{SO}%
\makeop{SL}%
\makeop{GL}%
\makeop{Conf}%
\makeop{Homeo}%
\makeop{Diff}%
\makeop{Isom}%
%
% K: Functional Analysis.
%
\makeop{Ind}%
\makeop{Sig}%
\makeop{Spec}%
%
% L: Other.
%
\makeop{Conv}%
\makeop{Max}%
\makeop{Min}%
\makeop{Mod}%
\makeop{Deg}%
\makeop{loc}%
%
%%%%%%%%%%%%%%%%%%%%%%%%%%%%%%%%%%%%%%%%%%%%%%%%%%%%%%%%%%%%%%%%%%%%%%%%%%%%%%%%%%%%%%%%%%%%%%%%%%%%%%%%%%%%%%%%%%%%%%%
%
% 7: Redundant Material.
%
%%%%%%%%%%%%%%%%%%%%%%%%%%%%%%%%%%%%%%%%%%%%%%%%%%%%%%%%%%%%%%%%%%%%%%%%%%%%%%%%%%%%%%%%%%%%%%%%%%%%%%%%%%%%%%%%%%%%%%%
%
% This file contains redundant functionality for constructing a bibliography. Although it is not used, it may one day
% prove useful, and so I leave it here.
%
% A: Before the citations, the file should contain:
%
% \newif\ifmakebiblio
%
% followed by either \makebibliotrue or \makebibliofalse.
%
% (note that, for conveniance, the commands \newif\ifmakebiblio and \makebibliofalse have been included at the
% beginning of this preamble. These commands should be removed before making the changes outlined here.
%
% B: Immediately before the first citation, the file should contain the following instructions:
%
% \ifmakebiblio%
% \openout\biblio=biblio.tex %
% \edef\fileoutput{\write\biblio{\bgroup\leftskip=2em}}%
% \fileoutput%
% \fi%
%
% C: Immediately after the last citation, the file should contain the following instruction:
%
% \ifmakebiblio%
% {\edef\fileoutput{\write\biblio{\egroup}}%
% \fileoutput}%
% \fi%
%

%% file: references.tex
\global\def\_@citation@Caffarelli{1}
\global\def\_@citation@CaffNirSprI{2}
\global\def\_@citation@CaffNirSprII{3}
\global\def\_@citation@CaffNirSprIII{4}
\global\def\_@citation@CaffNirSprV{5}
\global\def\_@citation@CalI{6}
\global\def\_@citation@CrandallIshiiLions{7}
\global\def\_@citation@GilbTrud{8}
\global\def\_@citation@GuanSpruckA{9}
\global\def\_@citation@Guan{10}
\global\def\_@citation@GuanSpruck{11}
\global\def\_@citation@GuanSpruckI{12}
\global\def\_@citation@GuanSpruckII{13}
\global\def\_@citation@GuanSpruckIII{14}
\global\def\_@citation@GuanSpruckIV{15}
\global\def\_@citation@IvochTomi{16}
\global\def\_@citation@Jorgens{17}
\global\def\_@citation@Lab{18}
\global\def\_@citation@Pog{19}
\global\def\_@citation@SmiDTI{20}
\global\def\_@citation@RosSpruck{21}
\global\def\_@citation@ShengUrbasWang{22}
\global\def\_@citation@Smale{23}
\global\def\_@citation@SmiAAT{24}
\global\def\_@citation@SmiCGC{25}
\global\def\_@citation@SmiSLC{26}
\global\def\_@citation@SmiNLD{27}
\global\def\_@citation@SmiFCS{28}
\global\def\_@citation@SmiPPH{29}
\global\def\_@citation@SmiNLP{30}
\global\def\_@citation@SmiPPG{31}
\global\def\_@citation@Spruck{32}
\global\def\_@citation@Trud{33}
\global\def\_@citation@TrudWang{34}
\global\def\_@citation@WhiteI{35}
\global\def\_@citation@WhiteII{36}
\global\def\_@subhead@SubheadNonLinearCurvatureFunctions{1.1}
\global\def\_@subhead@SubheadTheNonLinearPlateauProblem{1.2}
\global\def\_@proc@ThmFirstExistenceTheorem{1.3.1}
\global\def\_@proc@ThmSecondExistenceTheorem{1.3.2}
\global\def\_@proc@ThmDegreeUniqueness{1.3.3}
\global\def\_@head@HeadConvexCurvatureFunctions{2}
\global\def\_@proc@PropLimitOfCurvFns{2.1.1}
\global\def\_@proc@PropCurvatureQuotients{2.1.2}
\global\def\_@proc@PropCharacterisationOfFInfinity{2.1.3}
\global\def\_@proc@PropCurvFnsII{2.1.4}
\global\def\_@head@HeadConvexCobordismsAndEmbeddingRadii{3}
\global\def\_@subhead@SubheadDefinitionOfConvexCobordisms{3.1}
\global\def\_@subhead@SubheadGeometryOfConvexPrecobordisms{3.2}
\global\def\_@proc@PropGeometryOfCobordism{3.2.1}
\global\def\_@proc@PropCharacterisationOfConvexCobordisms{3.2.2}
\global\def\_@subhead@SubheadGeometryOfConvexCobordisms{3.3}
\global\def\_@proc@PropUniquenessOfConvexCobordisms{3.3.1}
\global\def\_@proc@PropExistsUniqueMinimiser{3.3.2}
\global\def\_@proc@PropDistInCobordismIsConvex{3.3.3}
\global\def\_@proc@PropDiameterBounds{3.3.4}
\global\def\_@subhead@SubheadGlueingAndExcision{3.4}
\global\def\_@proc@PropGeosDontReturn{3.4.1}
\global\def\_@proc@PropCobIsTrans{3.4.2}
\global\def\_@proc@PropHowTheHypersurfaceMeetsTheFoliation{3.4.3}
\global\def\_@proc@PropExcisionI{3.4.6}
\global\def\_@proc@PropGeoDoesNotReturnII{3.4.7}
\global\def\_@proc@PropExcisionII{3.4.8}
\global\def\_@subhead@SubheadOpennessAndClosedness{3.5}
\global\def\_@proc@PropTransversalityOfInteriorHypersurfaces{3.5.1}
\global\def\_@proc@PropBoundednessIsAnOpenProperty{3.5.2}
\global\def\_@proc@PropBoundednessIsAClosedProperty{3.5.3}
\global\def\_@subhead@SubheadBoundednessAndMinimalEmbeddingRadii{3.6}
\global\def\_@proc@ThmMinimalEmbeddingRadius{3.6.1}
\global\def\_@subhead@SubheadSupportingNormalVectors{3.7}
\global\def\_@proc@PropSupportingNormalIsContainedBelowOuterBarrier{3.7.1}
\global\def\_@proc@PropUniformModulusOfContinuity{3.7.2}
\global\def\_@head@HeadFirstOrderLowerEstimates{4}
\global\def\_@subhead@SubheadFirstOrderLowerEstimatesMainResults{4.1}
\global\def\_@proc@PropFirstOrderEstimatesUnboundedCase{4.1.1}
\global\def\_@proc@PropFirstOrderLowerEstimatesBoundedCase{4.1.2}
\global\def\_@proc@PropHessianOfRestriction{4.1.3}
\global\def\_@subhead@SubheadAnalyticPropertiesOfTheBarrierFunction{4.2}
\global\def\_@eqn@EqnDefinitionOfPhiOne{\relax \unhbox \voidb@x \hbox {(3)}}
\global\def\_@eqn@EqnFirstDerivativeOfTheseVectorFields{\relax \unhbox \voidb@x \hbox {(4)}}
\global\def\_@proc@PropFirstBarrierEstimate{4.2.1}
\global\def\_@eqn@EqnFirstVariationOfPrincipleCurvature{\relax \unhbox \voidb@x \hbox {(5)}}
\global\def\_@eqn@EqnAUsefulTermThatDoesNotVanish{\relax \unhbox \voidb@x \hbox {(6)}}
\global\def\_@proc@PropDistributionsAreClose{4.2.2}
\global\def\_@proc@PropControlOfDerivativeOfPhi{4.2.3}
\global\def\_@eqn@EqnControlOfPhiOne{\relax \unhbox \voidb@x \hbox {(8)}}
\global\def\_@proc@PropFirstBarrierEstimateII{4.2.4}
\global\def\_@subhead@SubheadGeometricPropertiesOfFirstOrderBarrier{4.3}
\global\def\_@proc@PropASuitableValueOfRhoRestrictsTheBoundary{4.3.2}
\global\def\_@eqn@EqnDefinitionOfPhi{\relax \unhbox \voidb@x \hbox {(8)}}
\global\def\_@proc@PropHowToChooseM{4.3.3}
\global\def\_@eqn@EqnLowerBoundsAlongTheBoundaryI{\relax \unhbox \voidb@x \hbox {(9)}}
\global\def\_@eqn@EqnLowerBoundsAlongTheBoundaryII{\relax \unhbox \voidb@x \hbox {(10)}}
\global\def\_@proc@PropLevelHypersurfacesAreManifolds{4.3.4}
\global\def\_@proc@PropCurvatureIsSmallerThanKappa{4.3.5}
\global\def\_@subhead@SubheadTheUnboundedCase{4.4}
\global\def\_@proc@PropControlOfFirstEigenvalue{4.4.1}
\global\def\_@proc@PropNonPositiveEigenvalue{4.4.2}
\global\def\_@proc@PropNonPositivePrincipleCurvature{4.4.3}
\global\def\_@head@HeadSecondOrderBoundaryEstimates{5}
\global\def\_@subhead@SubheadSecondOrderBoundaryEstimatesMainResult{5.1}
\global\def\_@proc@PropSecondOrderBoundaryEstimates{5.1.1}
\global\def\_@subhead@SubheadSecondOrderBoundaryEstimatesPreliminaryResults{5.2}
\global\def\_@proc@LemmaSuperHarmonicity{5.2.1}
\global\def\_@proc@CorSuperHarmonicityI{5.2.2}
\global\def\_@proc@LemmaCommutationRelations{5.2.3}
\global\def\_@proc@PropDerivativeOfCurvatureEquation{5.2.4}
\global\def\_@proc@PropSignedDistanceFunction{5.2.5}
\global\def\_@subhead@SubheadConstructingTheBarrierPartI{5.3}
\global\def\_@proc@LemmaDerOfInnerProduct{5.3.1}
\global\def\_@proc@PropPositivityOfPhiV{5.3.2}
\global\def\_@proc@PropLaplacian{5.3.3}
\global\def\_@subhead@SubheadConstructingTheBarrierPartII{5.4}
\global\def\_@proc@PropLowerCurvBdOfBarrier{5.4.1}
\global\def\_@proc@PropFirstBarrier{5.4.2}
\global\def\_@proc@PropYetAnotherSuperharmonicityRelation{5.4.3}
\global\def\_@subhead@ConstrucingTheBarrierPartIII{5.5}
\global\def\_@proc@PropSecondBarrier{5.5.1}
\global\def\_@head@HeadGlobalSecondOrderEstimates{6}
\global\def\_@subhead@SubheadGlobalSecondOrderEstimatesMainResults{6.1}
\global\def\_@proc@PropSecondOrderIntBoundsAlt{6.1.1}
\global\def\_@proc@PropSecondOrderIntBounds{6.1.2}
\global\def\_@proc@PropMuInfinityIsGreaterThanOne{6.2.1}
\global\def\_@proc@PropMuInfinityIsInfinite{6.2.2}
\global\def\_@subhead@SubheadGlobalSecondOrderBoundsPreliminaryResults{6.3}
\global\def\_@proc@CorSecondDerivativesOfA{6.3.1}
\global\def\_@proc@CorSuperharmonic{6.3.2}
\global\def\_@proc@PropSecondDerOfFirstEval{6.3.3}
\global\def\_@proc@PropSuperharmonic{6.3.4}
\global\def\_@proc@PropLapOfDistanceFunctionAlt{6.4.1}
\global\def\_@proc@CorMaximumPrincipalAlt{6.4.2}
\global\def\_@proc@CorSuperHarmonicityII{6.5.1}
\global\def\_@proc@PropLapOfDistanceFunction{6.5.2}
\global\def\_@proc@CorMaximumPrincipal{6.5.3}
\global\def\_@head@HeadExistence{7}
\global\def\_@subhead@SubheadExistenceCompactness{7.1}
\global\def\_@proc@PropCompactnessI{7.1.1}
\global\def\_@proc@PropCompactnessII{7.1.2}
\global\def\_@proc@PropIsotopies{7.2.1}